\documentclass[12pt]{article}
\usepackage[dvips]{graphicx}
\usepackage{amsmath, amssymb}
\usepackage{amsthm}

\usepackage{enumerate}

\newtheorem{theorem}{Theorem}[section]

\newtheorem{corollary}{Corollary}[section]

\newtheorem{proposition}{Proposition}[section]
\newtheorem{definition}{Definition}[section]

\newcommand{\tp}{^{\mathrm t}}

\newlength{\IndentI}
\newlength{\IndentII}
\newlength{\IndentIII}
\newlength{\WidthI}
\newlength{\WidthII}
\newlength{\WidthIII}
\numberwithin{equation}{section}

\title{Conformal Geometry of Sequential Test in Multidimensional Curved Exponential Family}

\author{
Masayuki Kumon\footnote{
Association for Promoting Quality Assurance in Statistics,\ Tokyo,\ Japan}  
\\
Akimichi Takemura
\footnote{
Graduate School of Information Science and Technology,
University of Tokyo,\ Tokyo,\ Japan}
\\
Kei Takeuchi
\footnote{
Emeritus, Graduate School of Economics, 
University of Tokyo,\ Tokyo,\ Japan}
}
\date{June 2014}

\begin{document}
\maketitle

\begin{abstract}
This article presents a differential geometrical method for analyzing sequential test procedures. 
It is based on the primal result on the conformal geometry of statistical manifold developed in Kumon, Takemura and Takeuchi (2011). 
By introducing curvature-type random variables, the condition is first clarified for a statistical manifold to be an exponential family under an appropriate sequential test procedure. 
This result is further elaborated for investigating the efficient sequential test in a multidimensional curved exponential family. 
The theoretical results are numerically examined by using von Mises-Fisher and hyperboloid models. 
\end{abstract}


\


\section{INTRODUCTION}

Sequential inferential procedure continues observations until the observed sample satisfies a certain prescribed criterion. 
Its properties have been shown to be superior on the average to those of nonsequential inferential procedure in which the number of observations is fixed a priori. 
Among others Takeuchi and Akahira (1988), Akahira and Takeuchi (1989) formulated the higher-order asymptotic theory of sequential estimation procedures rigorously and analyzed the higher-order efficiency of these procedures in the scalar parameter case. 
They showed that the exponential curvature term in the second-order variance can be eliminated by the sequential maximum likelihood estimator with the best stopping rule.

Okamoto, Amari and Takeuchi (1991) generalized these works to the multiparameter case by using geometrical notions, and studied characteristics of more general sequential estimation procedures. 
As indicated there, a sequential estimation procedure with an adequate stopping rule causes a nonuniform expansion of a statistical manifold  
(Riemannian manifold with a dual couple of affine connections) which is called the conformal transformation.
The result of Takeuchi and Akahira can be interpreted such that it is possible to reduce the exponential curvature of a statistical manifold to zero by a suitable conformal transformation. 
The conformal geometry of statistical manifold thus can be an adequate framework for the analysis of the sequential inferential procedures.

Kumon, Takemura and Takeuchi (2011) conducted this scheme for investigating the sequential estimating procedures from the information geometrical viewpoint. 
The dual conformal Weyl-Schouten curvature quantity of a statistical manifold was introduced there, and this quantity was proved to play a 
central role when considering the problem of covariance minimization 
under the sequential estimating procedures.

This article continues the scheme for investigating the sequential test  procedures. 
We focus on the powers of sequential tests, and pursue the possibility of uniformly efficient (most powerful) sequential test.
Different from the previous work, we first introduce the dual Euler-Schouten curvature random variables, and directly study their transformations under the sequential test procedures. 
These enable us to clarify the condition for a statistical manifold to be an exponential family under an adequate sequential test procedure.  
We then elaborate the result for analyzing the structure of the uniformly efficient sequential test in a multidimensional curved exponential family.  
The above approach is motivated by the fact that the exponentiality of a statistical manifold ensures the uniformly efficient inferential procedures in the nonsequential case as well as in the sequential case. 
Information geometry has been steadily producing mathematical 
methodologies for a variety of statistical sciences 
(see e.g., Amari, 1985; Amari et al., 1987; Amari and Nagaoka, 2000; Kumon, 2009, 2010).
The present article also is situated in these developments.

The organization of the article is as follows. 
The established known results are cited as propositions, and the results obtained in this article are stated as theorems. 
In Section 2, we prepare some statistical notations and preliminary 
results on the sequential test procedures which will be relevant in this article. 
In Section 3, we formulate the conformal transformation of statistical manifold, where the dual Euler-Schouten curvature tensor variables are   introduced. 
The related notion of totally exponential umbilic is proved to give the 
criterion for determining that a statistical manifold can be an exponential family under an suitable sequential test procedure. 
Its relation to the dual conformal flatness is also explained.    
In Section 4, we summarize the asymptotics on the nonsequential test in a multidimensional curved exponential family, which are already developed in Kumon and Amari (1983), and Amari (1985, Chapter 6). 
In Section 5, the general result in Section 3 is used to write down the structure of multidimensional curved exponential family under the sequential test procedures, where the meanings of the dual Euler-Schouten curvatures of statistical submanifold are clarified.  
In Section 6, we formulate the asymptotics of the sequential test in a multidimensional curved exponential family in line with the nonsequential case, and give the concrete condition and procedure for the uniformly efficient sequential test. 
In Section 7, the results in Section 6 and also in Section 4 are examined numerically by using two curved exponential families, the von Mises-Fisher model and the hyperboloid model, both of which admit the uniformly efficient sequential test.  
Section 8 is devoted to some additional discussions and a perspective of future work.

\section{PRELIMINARIES}

We denote by $\{X(t) = [X_1(t), \dots, X_k(t)]^{\tp}\ |\ t\in T\}$ a $k$-dimensional random sequence defined on the probability space $[\Omega, {\cal F}, P]$ 
with values in $[E, \cal{E}]$, where $E = \mathbb{R}^k$ and $\cal{E}$ is the $\sigma$-field of all Borel sets in $E$. 
The time parameter $t \in T$ runs over all positive integers 
$T = \{1, 2, \dots \}$. 
Moreover, let ${\cal F}_n, n \in T$ denote the $\sigma$-field generated by the $n$-tuple of random vectors $X^n = [X(1), \dots, X(n)]$.

In this article, we assume that the probability measure $P$ depends on an unknown parameter $\theta = (\theta^1, \dots, \theta^m)^{\tp} \in \Theta$, $P = P_\theta$, where $\Theta$ is homeomorphic to ${\mathbb R}^m$, 
and we shall consider the case where the following conditions are fulfilled:

\begin{enumerate}[(i)]
\item The random sequence $\{X(t)\ |\ t\in T\}$ consists of independent and identically distributed $k$-dimensional random variables 
$X(1), X(2), \dots $. 
\item The probability distribution of $X(1)$ is dominated by a $\sigma$-finite measure $\mu$ and has the density $f(x, \theta)$ with respect to $\mu$ 
\begin{align}
\label{eq:2-1}
\frac{dP_\theta}{d\mu}(x) = f(x, \theta),\quad x\in E,\ \theta \in \Theta,
\end{align}
so that the joint density $f(x^n, n, \theta)$ of 
$X^n = [X(1), \dots, X(n)]$ is written as
\begin{align}
\label{eq:2-2}
f(x^n, n, \theta) = \prod_{t=1}^n f(x(t), \theta),\quad 
x^n\in E^n,\ \theta \in \Theta.
\end{align}
\end{enumerate}
We say that $\{P_\theta\ |\ \theta\in \Theta\}$ is an $m$-dimensional full regular minimally represented exponential family (f.r.m.\ exponential family) when the density (\ref{eq:2-1}) can be written as
\begin{align}
\label{eq:2-3}
f_e(x, \theta) = \exp \{\theta^ix_i - \psi(\theta)\},
\end{align}
where $x = (x_1, \dots, x_m)^{\tp}\in {\mathbb R}^m$ is the canonical statistic, $\theta$ is the natural parameter, and 
$\eta = \partial \psi(\theta)/\partial \theta$ is the expectation parameter with $\psi(\theta)$ a smooth (infinitely differentiable) convex function of $\theta$. 
In the right-hand side of (\ref{eq:2-3}) and hereafter the Einstein summation convention will be assumed, so that summation will be automatically taken over indices repeated twice in the sense, for example, 
$\theta^ix_i = \sum_{i=1}^m \theta^ix_i$.

Sequential inferential procedures are characterized by a random sample size, where stopping times are used to stop the observations of the process. 
We denote by $\tau$ an arbitrary stopping time; that is, a random variable $\tau$ defined on $\Omega$ with values in $T\cup \{\infty\}$ and possessing the property $\{\omega \in \Omega: \tau(\omega) \leqslant n\}\in 
{\cal F}_n,\ \forall n\in T$. 
This article treats the case to which the Sudakov lemma applies, where the stopped sequence  
$\{(X^\tau, \tau)\ |\ \tau \in T\cup \infty\}$ has the density 
${\bar f}(x^n, n, \theta)$, and this  can be regarded as the same function in (\ref{eq:2-2}) 
(cf. e.g., Magiera, 1974).

As a typical statistical problem suppose that we wish to test a simple null hypothesis
\begin{align*}
H_0 : \theta = \theta_0 \in \Theta\quad \textrm{vs}\quad 
H_1 : \theta \neq \theta_0.
\end{align*} 
A sequential test procedure will be based on a test statistic $Y_n$ and a  stopping time $\tau_c \wedge n_1$, where  
\begin{align}
\label{eq:2-4}
\tau_c = \min\{n\ge n_0 : Y_n > c\},
\end{align}
if $\tau_c \le n_1$, $H_0$ should be rejected, whereas if $\tau_c > n_1$, 
$H_0$ should not be rejected. 
Here $n_0$ denotes a given initial sample size and $n_1$ denotes some given maximal sample size.
Specifically we take the stopping time 
$\tau_c \wedge n_1$ with 
\begin{align}
\label{eq:2-5}
\tau_c = \min\{n\ge n_0 : Y_n = n\cdot y(\bar{X}_n) > ca(n)\},\quad 
\bar{X}_n = \frac{1}{n}\sum_{t=1}^n X(t),
\end{align}
where $y$ and $a$ belong to the following classes of functions $\mathcal{Y}$ and $\mathcal{A}$, respectively. 
\begin{definition}
The function $y(\eta) \in {\cal Y}$ if $y(\eta)$ is positive and smooth in 
$\eta \in {\mathbb R}^k$. 
\end{definition}
\begin{definition}
The function $a(u) \in {\cal A}$ if $a(u)$ is nondecreasing, concave, smooth in $u \ge \exists A \ge 0$ and regularly varying at infinity with exponent $\beta\ (0\le \beta < 1)$ (see Gut, 2009, Appendix B).
\end{definition}
Let $\nu_c = \nu_c(\eta)$ be the solution of the equation
\begin{align*}
\nu_cy(\eta) = ca(\nu_c).
\end{align*}
In particular, if $a(u) = u^\beta$ for some $\beta\ (0\le \beta < 1)$, then 
$\nu_c = (c/y(\eta))^{1/(1-\beta)}$. 
For $\beta = 0$, that is, $a(u)\equiv 1$ and $y(\eta) = |c^0 + c^i\eta_i|$, the relation reduces to $\nu_c = c/|c^0 + c^i\eta_i|$.

The above class of stopping times $\tau_c$ is denoted by ${\cal C}$.  
The following is the central limit theorem for a stopping time 
$\tau_c \in {\cal C}$ (see Gut, 2009, Theorem 3.3 in Chapter 6).

\begin{proposition}
If $\sigma_{ij}(\theta) = Cov_\theta[X_i(1),X_j(1)] < \infty$ and $\partial^iy(\eta) \neq 0$ for any 
$\partial^i = \partial/\partial \eta_i$, then  
$\nu_c(\eta)\ \to\ \infty\ \textrm{as}\ c\ \to\ \infty$, and we have
\begin{align}
\label{eq:2-7}
\frac{\tau_c - \nu_c(\eta)}
{\sqrt{\nu_c(\eta)\sigma_{ij}(\theta)s^i_c(\eta)s^j_c(\eta)}}\ 
\xrightarrow{d}\ N(0, 1)\ \ \textrm{as}\ \ c\ \to\ \infty, 
\end{align}
where $s^i_c(\eta) = \partial^i \log \nu_c(\eta)$.
\end{proposition}
\

This result will be utilized for examining the exponentiality of statistical manifold under a sequential statistical procedure in the next section.

\section{CONFORMAL TRANSFORMATION OF STATISTICAL MANIFOLD}

Let $M = \{f(x, \theta)\ |\ \theta \in \Theta\}$ be an original family of probability densities of $X(1)$, where $\Theta$ is homeomorphic to ${\mathbb R}^m$. 
The family $M$ can be regarded as an $m$-dimensional statistical manifold, where the $m$-dimensional vector parameter $\theta$ serves as a coordinate system to specify a point, that is, a density $f(x, \theta)\in M$. 
In order to introduce differential geometrical notions in $M$, 
we associate a linear space $R_\theta$ with each point $\theta \in M$ 
as follows (cf. Amari and Kumon, 1988).
\begin{align*}
R_\theta = \{r(x)\ |\ \langle r(x)\rangle_\theta = 0,\ 
\langle r(x)^2\rangle_\theta < \infty \},\quad 
\langle r(x)\rangle_\theta = E_\theta[r(X)] = \int_{x \in E} r(x)f(x, \theta)d\mu(x).
\end{align*}
The $m$-dimensional linear space
\begin{align*}
T_\theta = \{r(x)\ |\ r(x) = a^i\partial_il_1(x, \theta) \},\quad 
l_1(x, \theta) = \log f(x, \theta),\quad 
\partial_i = \frac{\partial}{\partial \theta^i}
\end{align*} 
is a subspace of $R_\theta$ spanned by $m$ score functions 
$\partial_il_1(x, \theta),\ i = 1, \dots, m$, which is the tangent space at $\theta$ of $M$. 
We have a natural direct sum
\begin{align*}
R_\theta = T_\theta \oplus N_\theta,\quad 
N_\theta = \{r(x)\ |\ r(x)\in R_\theta,\ \langle r(x)s(x)\rangle_\theta = 0,\ 
\forall s(x)\in T_\theta \}, 
\end{align*}
where $N_\theta$ is the normal space at $\theta$ of $M$.  
The aggregates of $R_\theta's,\ T_\theta's,\ N_\theta's$ at all 
$\theta \in \Theta$ and the resulting direct sums are denoted by
\begin{align*}
{\cal R}_M = \underset{\theta \in \Theta}\cup R_\theta,\quad 
{\cal T}_M = \underset{\theta \in \Theta}\cup T_\theta,\quad
{\cal N}_M = \underset{\theta \in \Theta}\cup N_\theta,\quad 
{\cal R}_M = {\cal T}_M \oplus {\cal N}_M,
\end{align*} 
and these are called the  
\textit{Hilbert bundle, tangent bundle, normal bundle}  
over $M$, respectively.

A function $r(x, \theta)$ is called  
\textit{a smooth section} of the Hilbert bundle ${\cal R}_M$ when 
$r(x, \theta)$ is smooth in $\theta$ and 
$r(x, \theta)\in R_\theta,\ \forall \theta \in \Theta$. 
The set of all sections of ${\cal R}_M$ is denoted by
\begin{align*}
S({\cal R}_M) = \{r(x, \theta)\ |\ r(x, \theta)\in R_\theta \},
\end{align*}
which is again a linear space. 
The $S({\cal T}_M)$ and $S({\cal N}_M)$ are similarly defined, which are the subspaces of $S({\cal R}_M)$. 
Corresponding to ${\cal R}_M = {\cal T}_M \oplus {\cal N}_M$, we also have a natural direct sum
\begin{align*}
S({\cal R}_M) = S({\cal T}_M) \oplus S({\cal N}_M).
\end{align*}
The $m$ score functions $\partial_il_1(x, \theta),\ i = 1, \dots, m$ belong to $S({\cal T}_M)$, and their inner products are denoted by
\begin{align}
\label{eq:3-1}
g_{ij}(\theta) = \langle \partial_il_1\partial_jl_1\rangle_\theta,
\end{align}
which constitute the Fisher information metric tensor of $M$.

The one parameter family of $\alpha$-covariant derivative for 
$r(x, \theta)\in S({\cal R}_M)$ in the direction of $\theta^i$-coordinate is defined by
\begin{align}
\label{eq:3-2}
\nabla_i^{(\alpha)}r = \partial_i r - 
\frac{1+\alpha}{2}\langle \partial_ir\rangle_\theta + 
\frac{1-\alpha}{2}r\partial_il_1 
\in S({\cal R}_M).
\end{align}
Specifically for $\partial_jl_1(x, \theta)\in S({\cal T}_M)$, we have the 
direct sum decomposition
\begin{align}
\label{eq:3-3}
&\nabla_i^{(\alpha)}\partial_jl_1 = \partial_i\partial_jl_1 + 
\frac{1+\alpha}{2}g_{ij}(\theta) + \frac{1-\alpha}{2}\partial_il_1\partial_jl_1 = 
\Gamma^{(\alpha)k}_{ij}(\theta)\partial_kl_1(x, \theta) + 
h_{ij}^{(\alpha)}(x, \theta),\\ 
&\Gamma^{(\alpha)k}_{ij}(\theta)\partial_kl_1(x, \theta) \in S({\cal T}_M),
\quad h_{ij}^{(\alpha)}(x, \theta) \in S({\cal N}_M),\nonumber
\end{align}
where
\begin{align}
&\Gamma^{(\alpha)k}_{ij}(\theta) = 
\Gamma^{(\alpha)}_{ijl}(\theta)g^{lk}(\theta),\nonumber \\ 
\label{eq:3-4}
&\Gamma_{ijk}^{(\alpha)}(\theta) = 
\langle \nabla_i^{(\alpha)}\partial_jl_1 \partial_kl_1\rangle_\theta = 
\langle \partial_i\partial_jl_1\partial_kl_1\rangle_\theta + 
\frac{1-\alpha}{2}T_{ijk}(\theta),\\
\label{eq:3-5}
&T_{ijk}(\theta) = 
\langle \partial_il_1\partial_jl_1\partial_kl_1\rangle_\theta,
\end{align}
$T_{ijk}(\theta)$ is called the 
\textit{skewness tensor} of $M$,  $\Gamma_{ijk}^{(\alpha)}(\theta)$ is the one parameter family of affine connections called the 
\textit{$\alpha$-connection} of $M$, and $h_{ij}^{(\alpha)}(x, \theta)$ is called the 
\textit{$\alpha$-Euler-Schouten curvature tensor variable} of $M$ in 
${\cal R}_M$.

The $\alpha$-curvature of the Hilbert bundle ${\cal R}_M$ is measured by the second-order differential operator 
$K_{ij}^{(\alpha)} = \nabla_i^{(\alpha)}\nabla_j^{(\alpha)} - \nabla_j^{(\alpha)}\nabla_i^{(\alpha)}$ in 
$S({\cal R}_M)$, and by direct calculation we have
\begin{align*}
K_{ij}^{(\alpha)} r = 
\frac{1-\alpha^2}{4}(\langle \partial_ir\rangle_\theta \partial_jl_1 - 
\langle \partial_jr\rangle_\theta \partial_il_1).
\end{align*}
Specifically for $\partial_kl_1(x, \theta)\in S({\cal T}_M)$, we define
the $\alpha$-enveloping curvature tensor of $M$ by
\begin{align}
\label{eq:3-6}
K_{ijkl}^{(\alpha)}(\theta) = 
\langle K_{ij}^{(\alpha)}\partial_kl_1 \partial_ll_1\rangle_\theta =  
\frac{1-\alpha^2}{4}(g_{jk}g_{il} - g_{ik}g_{jl}).
\end{align}

On the other hand, the intrinsic curvature of $M$ is measured by the the $\alpha$-Riemann-Christoffel curvature tensor 
\begin{align}
\label{eq:3-7}
R_{ijkl}^{(\alpha)}(\theta) = \partial_i\Gamma_{jkl}^{(\alpha)} - 
\partial_j\Gamma_{ikl}^{(\alpha)} + 
g^{rs}(\Gamma_{ikr}^{(\alpha)}\Gamma_{jsl}^{(\alpha)} - 
\Gamma_{jkr}^{(\alpha)}\Gamma_{isl}^{(\alpha)}).
\end{align}
The $\alpha$- and the $(-\alpha)$-connections are mutually dual
\begin{align*}
\partial_ig_{jk} = \Gamma_{ijk}^{(\alpha)} + 
\Gamma_{ikj}^{(-\alpha)},
\end{align*}
and the $\pm \alpha$-RC curvature tensors are in the dual relation
\begin{align*}
R_{ijkl}^{(\alpha)} = - R_{ijlk}^{(-\alpha)}.
\end{align*}
We can regard (\ref{eq:3-3}) as the partial differential equation for 
$\partial_jl_1(x, \theta)$, and its integrability condition
\begin{align*}
\partial_k \partial_i \partial_j l_1(x, \theta) = 
\partial_i \partial_k \partial_j l_1(x, \theta)
\end{align*}
yields the equation of Gauss 
(cf. Amari, 1985, p.52; Vos, 1989)
\begin{align}
\label{eq:3-8}
R_{ijkl}^{(\alpha)}(\theta) = K_{ijkl}^{(\alpha)}(\theta) + 
\langle h_{jk}^{(\alpha)}h_{il}^{(-\alpha)} - 
h_{ik}^{(\alpha)}h_{jl}^{(-\alpha)}\rangle_\theta,
\end{align}
which connects three kinds of curvature quantities.
Among them, the 1-ES curvature tensor variable $h_{ij}^{(1)}(x, \theta)$ plays an important role in the evaluations of various statistical inferences. 
In fact we have the following criterion for a statistical manifold $M$ to be an f.r.m.\ exponential family $M_e$.

\begin{theorem}
A statistical manifold $M = \{f(x, \theta)\ |\ \theta \in \Theta \}$ is 
an f.r.m.\ exponential family 
$M_e = \{f_e(x, \theta)\ |\ \theta \in \Theta \}$ if and only if 
\begin{align}
\label{eq:3-9}
&h_{ij}^{(1)}(x, \theta) = \partial_i\partial_jl_1(x, \theta) + 
g_{ij}(\theta) - \Gamma_{ij}^{(1)k}(\theta)\partial_kl_1(x, \theta) 
\equiv 0.
\end{align}
\end{theorem}
\

\begin{proof}
Suppose that $M$ is an f.r.m.\ exponential family $M_e$ with the natural parameter $\theta$. 
From (\ref{eq:2-3}), we have 
$l_1(x, \theta) = \log f_e(x, \theta) = \theta^ix_i - \psi(\theta)$.  
Then a direct calculation and (\ref{eq:3-4}) with $\alpha = 1$ yield
\begin{align*}
&\partial_i\partial_jl_1(x, \theta) = -\partial_i\partial_j\psi(\theta) = -g_{ij}(\theta),\quad \Gamma_{ijk}^{(1)}(\theta) \equiv 0\
\Rightarrow\ h_{ij}^{(1)}(x, \theta) \equiv 0.
\end{align*}

Conversely suppose that (\ref{eq:3-9}) holds with a certain parameterization $\theta$. 
Then from the equation of Gauss (\ref{eq:3-8}) with $\alpha = 1$, 
we have $R_{ijkl}^{(1)}(\theta) \equiv 0$, so that there exists a coordinate system $(\xi^\alpha)$ such that $\Gamma_{\alpha\beta}^{(1)\gamma}(\xi) \equiv 0$ holds. 
Since $h_{ij}^{(1)}(x, \theta)$ is the tensor variable, 
$h_{\alpha\beta}^{(1)}(x, \xi) \equiv 0$ also holds, and hence
\begin{align*}
&h_{\alpha\beta}^{(1)}(x, \xi) = 
\partial_\alpha\partial_\beta l_1(x, \xi) + g_{\alpha\beta}(\xi) \equiv 0\\
&\Rightarrow\ \partial_\gamma \partial_\alpha\partial_\beta l_1(x, \xi) = 
-\partial_\gamma g_{\alpha\beta}(\xi) = 
-\partial_\alpha g_{\gamma\beta}(\xi) 
= \partial_\alpha \partial_\gamma\partial_\beta l_1(x, \xi)\\
&\Rightarrow\ \exists \psi(\xi)\ 
\textrm{smooth and convex in $\xi$ such that}\ 
g_{\alpha\beta}(\xi) = \partial_\alpha \partial_\beta \psi(\xi)\\
&\Rightarrow\ \partial_\alpha\partial_\beta l_1(x, \xi) + 
\partial_\alpha \partial_\beta \psi(\xi) = 0\\
&\Rightarrow\ l_1(x, \xi) = l_0(x) + \xi^\alpha t_\alpha(x) - \psi(\xi),
\end{align*}
which shows that $M$ is an f.r.m.\ exponential family $M_e$ with the natural parameter $\xi$.
\end{proof}
\

Related to $h_{ij}^{(1)}(x, \theta)$, we define the mean 1-ES curvature variable $h^{(1)}(x, \theta)$ and the 1-ES umbilic curvature tensor variable $k_{ij}^{(1)}(x, \theta)$ of $M$ in ${\cal R}_M$ by
\begin{align}
\label{eq:3-10}
&h^{(1)}(x, \theta) = \frac{1}{m}h_{ij}^{(1)}(x, \theta)g^{ij}(\theta)
\in S({\cal N}_M),\\
\label{eq:3-11}
&k_{ij}^{(1)}(x, \theta) = h_{ij}^{(1)}(x, \theta) - 
g_{ij}(\theta)h^{(1)}(x, \theta)\in S({\cal N}_M).
\end{align}
We say that a statistical manifold $M$ is 
\textit{totally exponential umbilic (totally $e$-umbilic)} in the Hilbert bundle ${\cal R}_M$ when
\begin{align}
\label{eq:3-12}
k_{ij}^{(1)}(x, \theta) = h_{ij}^{(1)}(x, \theta) - 
g_{ij}(\theta)h^{(1)}(x, \theta) \equiv 0,
\end{align}
and thus nonzero $k_{ij}^{(1)}(x, \theta)$ represents the anisotropy of the $1$-ES curvature of $M$ in ${\cal R}_M$.
By definition and Theorem 3.1, we know that
\begin{align}
\label{eq:3-13}
\textrm{$M$ is an f.r.m. exponential family $M_e$}\
\Rightarrow\ \textrm{$M$ is totally $e$-umbilic in ${\cal R}_M$}.
\end{align}
The squared scalar quantity of $h_{ij}^{(1)}(x, \theta)$ is denoted by 
$\lambda^2$, which can be decomposed into
\begin{align}
\label{eq:3-14}
&\lambda^2 = 
\langle h_{ij}^{(1)}(x, \theta)h_{kl}^{(1)}(x, \theta)\rangle_\theta
g^{ik}(\theta)g^{jl}(\theta) = \kappa^2 + \gamma^2,\\
\label{eq:3-15}
&\kappa^2 = 
\langle k_{ij}^{(1)}(x, \theta)k_{kl}^{(1)}(x, \theta)\rangle_\theta
g^{ik}(\theta)g^{jl}(\theta),\\
\label{eq:3-16}
&\gamma^2 = m\langle h^{(1)}(x, \theta)^2\rangle_\theta.
\end{align}
In the evaluation of power functions of test procedures, two kinds of nonnegative scalar $1$-ES curvatures $\kappa^2$ and $\gamma^2$ 
will be the crucial quantities as is shown in the next section.

We turn to the conformal transformation of the original statistical manifold $M$. 
Let $\bar{M} = \{\bar{f}(x^n, n, \theta)\ |\ \theta \in \Theta\}$ be a 
family of densities $\bar{f}(x^n, n, \theta)$ under a sequential inferential procedure. 
This family $\bar{M}$ is also an $m$-dimensional extended statistical manifold with the coordinate system $(\theta^i)$ specifying 
$\bar{f}(x^n, n, \theta) \in \bar{M}$. 
With each point $\theta \in \bar{M}$ we again associate a linear space
\begin{align*}
&\bar{R}_\theta = \{\bar{r}(x^n, n)\ |\ 
\langle \bar{r}(x^n, n)\rangle_\theta = 0,\ 
\langle \bar{r}(x^n, n)^2\rangle_\theta < \infty \},\\ 
&\langle \bar{r}(x^n, n)\rangle_\theta = E_\theta[\bar{r}(X^\tau, \tau)] = 
\int_{x^n \in E^n} \bar{r}(x^n, n)\bar{f}(x^n, n, \theta)d\mu(x^n, n).
\end{align*}
The $m$-dimensional linear space
\begin{align*}
\bar{T}_\theta = \{\bar{r}(x^n, n)\ |\ \bar{r}(x^n, n) = 
a^i\partial_i\bar{l}_n(x^n, \theta) \},\quad 
\bar{l}_n(x^n, \theta) = \log \bar{f}(x^n, n, \theta)
\end{align*} 
is a subspace of $\Bar{R}_\theta$ spanned by $m$ score functions 
$\partial_i\bar{l}_n(x^n, \theta),\ i = 1, \dots, m$, which is the tangent space at $\theta$ of $\bar{M}$. 
There is a natural direct sum
\begin{align*}
\bar{R}_\theta = \bar{T}_\theta \oplus \bar{N}_\theta,\quad 
\bar{N}_\theta = \{\bar{r}(x^n, n)\ |\ \bar{r}(x^n, n)\in 
\bar{R}_\theta,\ \langle \bar{r}(x^n, n)\bar{s}(x^n, n)\rangle_\theta = 0,\ 
\forall \bar{s}(x^n, n)\in \bar{T}_\theta \}, 
\end{align*}
where $\bar{N}_\theta$ is the normal space at $\theta$ of $\bar{M}$.  
The aggregates of $\bar{R}_\theta's,\ \bar{T}_\theta's,\ \bar{N}_\theta's$ at all $\theta \in \Theta$ and the resulting direct sums are denoted by
\begin{align*}
\bar{{\cal R}}_{\bar{M}} = \underset{\theta \in \Theta}\cup \bar{R}_\theta,\quad 
\bar{{\cal T}}_{\bar{M}} = \underset{\theta \in \Theta}\cup \bar{T}_\theta,\quad
\bar{{\cal N}}_{\bar{M}} = \underset{\theta \in \Theta}\cup \bar{N}_\theta,\quad 
\bar{{\cal R}}_{\bar{M}} = \bar{{\cal T}}_{\bar{M}} \oplus 
\bar{{\cal N}}_{\bar{M}},
\end{align*} 
and these are the Hilbert bundle, tangent bundle, normal bundle over 
$\bar{M}$, respectively.

The set of all sections of the Hilbert bundle $\bar{{\cal R}}_{\bar{M}}$ is denoted by
\begin{align*}
S(\bar{{\cal R}}_{\bar{M}}) = 
\{\bar{r}(x^n, n, \theta)\ |\ \bar{r}(x^n, n, \theta)\in \bar{R}_\theta \},
\end{align*}
which is again a linear space. 
The $S(\bar{{\cal T}}_{\bar{M}})$ and $S(\bar{{\cal N}}_{\bar{M}})$ are similarly defined, which are the subspaces of $S(\bar{{\cal R}}_{\bar{M}})$. 
Corresponding to $\bar{{\cal R}}_{\bar{M}} = 
\bar{{\cal T}}_{\bar{M}} \oplus \bar{{\cal N}}_{\bar{M}}$, there is also a natural direct sum
\begin{align*}
S(\bar{{\cal R}}_{\bar{M}}) = S(\bar{{\cal T}}_{\bar{M}}) \oplus 
S(\bar{{\cal N}}_{\bar{M}}).
\end{align*}
The $m$ score functions $\partial_i\bar{l}_n(x^n, \theta),\ i = 1, \dots, m$ belong to $S(\bar{{\cal T}}_{\bar{M}})$, and their inner products 
constitute the Fisher information metric tensor 
$\bar{g}_{ij}(\theta)$ of $\bar{M}$. 
From the Wald identity this is expressed as (see Akahira and Takeuchi, 1989)
\begin{align}
\label{eq:3-17}
\bar{g}_{ij}(\theta) = 
\langle \partial_i\bar{l}_n\partial_j\bar{l}_n\rangle_\theta = 
\nu(\theta)g_{ij}(\theta),\quad 
\nu(\theta) = \langle n\rangle_\theta = E_\theta[\tau].
\end{align}

The one parameter family of $\alpha$-covariant derivative for 
$\bar{r}(x^n, n, \theta)\in S(\bar{{\cal R}}_{\bar{M}})$ in the direction of $\theta^i$-coordinate is defined by
\begin{align}
\label{eq:3-18}
\bar{\nabla}_i^{(\alpha)}\bar{r} = \partial_i \bar{r} - 
\frac{1+\alpha}{2}\langle \partial_i\bar{r}\rangle_\theta + 
\frac{1-\alpha}{2}\bar{r}\partial_i\bar{l}_n 
\in S(\bar{{\cal R}}_{\bar{M}}).
\end{align}
Specifically for 
$\partial_j\bar{l}_n(x^n, \theta)\in S(\bar{{\cal T}}_{\bar{M}})$, we have 
the direct sum decomposition
\begin{align}
\label{eq:3-19}
&\bar{\nabla}_i^{(\alpha)}\partial_j\bar{l}_n = 
\partial_i\partial_j\bar{l}_n + 
\frac{1+\alpha}{2}\bar{g}_{ij}(\theta) + \frac{1-\alpha}{2}
\partial_i \bar{l}_n\partial_j\bar{l}_n = 
\bar{\Gamma}^{(\alpha)k}_{ij}(\theta)\partial_k\bar{l}_n(x^n, \theta) + 
\bar{h}_{ij}^{(\alpha)}(x^n, n, \theta),\\ 
&\bar{\Gamma}^{(\alpha)k}_{ij}(\theta)\partial_k\bar{l}_n(x^n, \theta) \in 
S(\bar{{\cal T}}_{\bar{M}}),
\quad \bar{h}_{ij}^{(\alpha)}(x^n, n, \theta) \in 
S(\bar{{\cal N}}_{\bar{M}}),\nonumber
\end{align}
where
\begin{align}
&\bar{\Gamma}^{(\alpha)k}_{ij}(\theta) = 
\bar{\Gamma}^{(\alpha)}_{ijl}(\theta)\bar{g}^{lk}(\theta),\nonumber \\ 
\label{eq:3-20}
&\bar{\Gamma}_{ijk}^{(\alpha)}(\theta) = 
\langle \bar{\nabla}_i^{(\alpha)}\partial_j\bar{l}_n \partial_k\bar{l}_n\rangle_\theta = 
\langle \partial_i\partial_j\bar{l}_n\partial_k\bar{l}_n\rangle_\theta + 
\frac{1-\alpha}{2}\bar{T}_{ijk}(\theta),\\
\label{eq:3-21}
&\bar{T}_{ijk}(\theta) = 
\langle \partial_i\bar{l}_n\partial_j\bar{l}_n\partial_k\bar{l}_n\rangle_\theta,
\end{align}
$\bar{T}_{ijk}(\theta)$ is the skewness tensor of $\bar{M}$,    $\bar{\Gamma}_{ijk}^{(\alpha)}(\theta)$ is the $\alpha$-connection of $\bar{M}$, and $\bar{h}_{ij}^{(\alpha)}(x^n, n, \theta)$ is the $\alpha$-ES curvature tensor variable of $\bar{M}$ in $\bar{{\cal R}}_{\bar{M}}$. 
Also from the Wald identity, $\bar{T}_{ijk}(\theta)$ and $\bar{\Gamma}_{ijk}^{(\alpha)}(\theta)$ are expressed as
\begin{align}
\label{eq:3-22}
&\bar{T}_{ijk}(\theta) = \nu[T_{ijk} + 3g_{(ij}s_{k)}],\quad 
3g_{(ij}s_{k)} = g_{ij}s_{k} + g_{jk}s_{i} + g_{ki}s_{j},\\
\label{eq:3-23}
&{\bar \Gamma}_{ijk}^{(\alpha)}(\theta) = 
\nu\left[\Gamma_{ijk}^{(\alpha)} + \frac{1-\alpha}{2}(g_{ki}s_j + g_{kj}s_i)
- \frac{1+\alpha}{2}g_{ij}s_k\right],\quad 
s_k(\theta) = \partial_k \log \nu(\theta).
\end{align}
The relations (\ref{eq:3-17}) and (\ref{eq:3-22}) show that a sequential inferential procedure induces a conformal transformation of statistical manifold $M \mapsto \bar{M}$ by the gauge function $\nu(\theta) > 0$.

The $\alpha$-curvature of the Hilbert bundle $\bar{{\cal R}}_{\bar{M}}$ is measured by the second-order differential operator 
$\bar{K}_{ij}^{(\alpha)} = \bar{\nabla}_i^{(\alpha)}\bar{\nabla}_j^{(\alpha)} - \bar{\nabla}_j^{(\alpha)}\bar{\nabla}_i^{(\alpha)}$ in 
$S(\bar{{\cal R}}_{\bar{M}})$, and by direct calculation we have
\begin{align*}
\bar{K}_{ij}^{(\alpha)} \bar{r} = 
\frac{1-\alpha^2}{4}(\langle \partial_i\bar{r}\rangle_\theta 
\partial_j\bar{l}_n - 
\langle \partial_j\bar{r}\rangle_\theta \partial_i\bar{l}_n).
\end{align*}
Specifically for 
$\partial_k\bar{l}_n(x^n, \theta)\in S(\bar{{\cal T}}_{\bar{M}})$, 
the $\alpha$-enveloping curvature tensor of $\bar{M}$ is given by
\begin{align}
\label{eq:3-24}
\bar{K}_{ijkl}^{(\alpha)}(\theta) = 
\langle \bar{K}_{ij}^{(\alpha)}\partial_k\bar{l}_n \partial_l\bar{l}_n\rangle_\theta =  
\frac{1-\alpha^2}{4}(\bar{g}_{jk}\bar{g}_{il} - \bar{g}_{ik}\bar{g}_{jl}).
\end{align}
The intrinsic curvature of $\bar{M}$ is measured by the the $\alpha$-RC curvature tensor 
\begin{align}
\label{eq:3-25}
\bar{R}_{ijkl}^{(\alpha)}(\theta) = \partial_i\bar{\Gamma}_{jkl}^{(\alpha)} - \partial_j\bar{\Gamma}_{ikl}^{(\alpha)} + 
\bar{g}^{rs}(\bar{\Gamma}_{ikr}^{(\alpha)}\bar{\Gamma}_{jsl}^{(\alpha)} - 
\bar{\Gamma}_{jkr}^{(\alpha)}\bar{\Gamma}_{isl}^{(\alpha)}),
\end{align}
and from (\ref{eq:3-7}), (\ref{eq:3-23}) this is expressed as
\begin{align}
\label{eq:3-26}
&{\bar R}_{ijkl}^{(\alpha)} = 
\nu[R_{ijkl}^{(\alpha)} - g_{il}s_{jk}^{(\alpha)} + g_{jl}s_{ik}^{(\alpha)}
- g_{jk}s_{il}^{(-\alpha)} + g_{ik}s_{jl}^{(-\alpha)}],\\
\label{eq:3-27}
&s_{ij}^{(\alpha)} = \frac{1-\alpha}{2}
\left[\nabla_i^{(\alpha)}s_j - \frac{1-\alpha}{2}s_is_j + 
\frac{1+\alpha}{4}g_{ij}s_ks_lg^{kl}\right],\quad 
\nabla_i^{(\alpha)}s_j = \partial_is_j - \Gamma_{ij}^{(\alpha)k}s_k.
\end{align}
We note that under a conformal transformation the mutual duality of $\pm\alpha$-connections is preserved 
\begin{align*}
\partial_i{\bar g}_{jk} = {\bar \Gamma}_{ijk}^{(\alpha)} + 
{\bar \Gamma}_{ikj}^{(-\alpha)},
\end{align*}
and also the dual relation of the $\pm \alpha$-RC curvature tensors is preserved
\begin{align*}
{\bar R}_{ijkl}^{(\alpha)} = - {\bar R}_{ijlk}^{(-\alpha)}.
\end{align*}
These can be confirmed by direct calculations with 
(\ref{eq:3-17}), (\ref{eq:3-23}), and (\ref{eq:3-26}).

From a statistical viewpoint, the flatness in terms of the $\pm 1$-RC curvature tensors is also important (cf. Kumon et al., 2011). 
Thus we say that a statistical manifold $M$ is 
\textit{conformally mixture (exponential) flat} 
when there exists a gauge function $\nu(\theta) > 0$ such that 
${\bar R}_{ijkl}^{(-1)} \equiv 0\ ({\bar R}_{ijkl}^{(1)} \equiv 0)$ holds. 
Note that by (\ref{eq:3-31}) $M$ is conformally mixture flat if and only if $M$ is conformally exponential flat. 
For the sake of simplicity we hereafter express this notion as 
conformally $m(e)$-flat.

We can regard (\ref{eq:3-19}) as the partial differential equation for 
$\partial_j\bar{l}_n(x^n, \theta)$, and its integrability condition
\begin{align*}
\partial_k \partial_i \partial_j \bar{l}_n(x^n, \theta) = 
\partial_i \partial_k \partial_j \bar{l}_n(x^n, \theta)
\end{align*}
again yields the equation of Gauss
\begin{align}
\label{eq:3-28}
\bar{R}_{ijkl}^{(\alpha)}(\theta) = \bar{K}_{ijkl}^{(\alpha)}(\theta) + 
\langle \bar{h}_{jk}^{(\alpha)}\bar{h}_{il}^{(-\alpha)} - 
\bar{h}_{ik}^{(\alpha)}\bar{h}_{jl}^{(-\alpha)}\rangle_\theta.
\end{align}

The exponentiality of a statistical manifold is the key notion also in the sequential inferential procedure 
(cf. Proposition 2.1 in Kumon et al., 2011). 
In view of this fact, we call 
$\bar{M} = \{\bar{f}(x^n, n, \theta)\ |\ \theta \in \Theta \}$ an    
\textit{f.r.m.\ conformal exponential family} 
denoted by $\bar{M}_e^* = 
\{\bar{f}_e(z, \bar{\theta})\ |\ \bar{\theta} \in \bar{\Theta} \}$ when the density functions $\bar{f}(x^n, n, \theta)$ are rewritten as
\begin{align}
\label{eq:3-29}
\bar{f}(x^n, n, \theta) = \bar{f}_e(z, \bar{\theta}) = 
\exp\{\bar{\theta}^\alpha z_\alpha - \bar{\psi}(\bar{\theta})\},
\end{align}
where $z = z(x^n, n) = (z_1, \dots, z_m)^{\tp} \in {\mathbb R}^m$ serves as the canonical statistic, 
$\bar{\theta} = (\bar{\theta}^1, \dots, \bar{\theta}^m)\in \bar{\Theta}$ 
serves as the natural parameter with $\bar{\Theta}$ homeomorphic to ${\mathbb R}^m$, and $\bar{\eta} = \partial \bar{\psi}(\bar{\theta})/\partial \bar{\theta}$ serves as the expectation parameter with $\bar{\psi}(\bar{\theta})$ a smooth convex function of $\bar{\theta}$. 
For this notion we have the following criterion as in Theorem 3.1.

\begin{theorem}
An extended statistical manifold 
$\bar{M} = \{\bar{f}(x^n, n, \theta)\ |\ \theta \in \Theta \}$ is an  
f.r.m.\ conformal exponential family 
$\bar{M}_e^* = 
\{\bar{f}_e(z, \bar{\theta})\ |\ \bar{\theta} \in \bar{\Theta} \}$ 
if and only if 
\begin{align}
\label{eq:3-30}
&\bar{h}_{ij}^{(1)}(x^n, n, \theta) = 
\partial_i\partial_j\bar{l}_n(x^n, \theta) + 
\bar{g}_{ij}(\theta) - 
\bar{\Gamma}_{ij}^{(1)k}(\theta)\partial_k\bar{l}_n(x^n, \theta) 
\equiv 0.
\end{align}
\end{theorem}
\

\begin{proof}
Suppose that $\bar{M}$ is an f.r.m.\ conformal exponential family $\bar{M}_e^*$ with the natural parameter $\bar{\theta}$. 
From (\ref{eq:3-29}), we have 
$\bar{l}_e(z, \bar{\theta}) = \log \bar{f}_e(z, \bar{\theta}) = 
\bar{\theta}^\alpha z_\alpha - \bar{\psi}(\bar{\theta})$. 
Then a direct calculation and (\ref{eq:3-20}) with $\alpha = 1$ yield
\begin{align*}
&\bar{\partial}_\alpha \bar{\partial}_\beta \bar{l}_e(z, \bar{\theta}) = 
-\bar{\partial}_\alpha \bar{\partial}_\beta \bar{\psi}(\bar{\theta}) = 
-\bar{g}_{\alpha \beta}(\bar{\theta}),\quad 
\bar{\partial}_\alpha = \partial/\partial \bar{\theta}^\alpha,\quad
\bar{\Gamma}_{\alpha \beta \gamma}(\bar{\theta}) \equiv 0\
\Rightarrow\ \bar{h}_{\alpha \beta}^{(1)}(z, \bar{\theta}) \equiv 0,
\end{align*}
and since $\bar{h}_{\alpha \beta}^{(1)}(z, \bar{\theta})$ is the tensor variable, $\bar{h}_{ij}^{(1)}(x^n, n, \theta) \equiv 0$ also holds.

Conversely suppose that (\ref{eq:3-30}) holds with a certain 
parameterization $\theta$. 
Then from the equation of Gauss (\ref{eq:3-28}) with $\alpha = 1$, 
we have $\bar{R}_{ijkl}^{(1)}(\theta) \equiv 0$, so that there exists a coordinate system $(\bar{\theta}^\alpha)$ such that 
$\bar{\Gamma}_{\alpha \beta}^{(1)\gamma}(\bar{\theta}) \equiv 0$. 
Since 
$\bar{h}_{\alpha \beta}^{(1)}(x^n, n, \bar{\theta}) \equiv 0$ also holds, 
we have
\begin{align*}
&\bar{h}_{\alpha \beta}^{(1)}(x^n, n, \bar{\theta}) = 
\bar{\partial}_\alpha \bar{\partial}_\beta \bar{l}_n(x^n, \bar{\theta}) + \bar{g}_{\alpha\beta}(\bar{\theta}) \equiv 0\\
&\Rightarrow\ 
\bar{\partial}_\gamma \bar{\partial}_\alpha \bar{\partial}_\beta 
\bar{l}_n(x^n, \bar{\theta}) = 
-\bar{\partial}_\gamma \bar{g}_{\alpha \beta}(\bar{\theta}) = 
-\bar{\partial}_\alpha \bar{g}_{\gamma \beta}(\bar{\theta}) 
= \bar{\partial}_\alpha \bar{\partial}_\gamma \bar{\partial}_\beta 
\bar{l}_n(x^n, \bar{\theta})\\
&\Rightarrow\ \exists \bar{\psi}(\bar{\theta})\ 
\textrm{smooth and convex in $\bar{\theta}$ such that}\ 
\bar{g}_{\alpha \beta}(\bar{\theta}) = 
\bar{\partial}_\alpha \bar{\partial}_\beta \bar{\psi}(\bar{\theta})\\
&\Rightarrow\ \bar{\partial}_\alpha \bar{\partial}_\beta 
\bar{l}_n(x^n, \bar{\theta}) + 
\bar{\partial}_\alpha \bar{\partial}_\beta \bar{\psi}(\bar{\theta}) = 0\\
&\Rightarrow\ \bar{l}_n(x^n, \bar{\theta}) = 
\bar{l}_0(x^n, n) + \bar{\theta}^\alpha z_\alpha(x^n, n) - 
\bar{\psi}(\bar{\theta}),
\end{align*}
which shows that $\bar{M}$ is an f.r.m.\ conformal exponential family 
$\bar{M}_e^*$ with the natural parameter $\bar{\theta}$.
\end{proof}
\

The 1-ES curvature tensor variable 
${\bar h}_{ij}^{(1)}(x^n, n, \theta)$ under a conformal transformation is decomposed into
\begin{align}
\label{eq:3-31}
&{\bar h}_{ij}^{(1)}(x^n, n, \theta) = 
h_{ij}^{(1)}(x^n, n, \theta) - g_{ij}(\theta)n^{(1)}(x^n, \theta),\\ 
\label{eq:3-32}
&h_{ij}^{(1)}(x^n, n, \theta) = 
\partial_i\partial_j{\bar l}_n + n g_{ij} - 
\Gamma_{ij}^{(1)k}\partial_k{\bar l}_n = 
\sum_{t=1}^n h_{ij}^{(1)}(x(t), \theta)\in S(\bar{{\cal N}}_{\bar{M}}),\\ 
\label{eq:3-33}
&n^{(1)}(x^n, \theta) = n - \nu - s^i\partial_i{\bar l}_n
\in S(\bar{{\cal N}}_{\bar{M}}),\quad 
s_i = s_jg^{ji}.
\end{align}
Under a conformal transformation, the mean 1-ES curvature variable 
${\bar h}^{(1)}(x^n, n, \theta)$ is expressed as 
\begin{align}
\label{eq:3-34}
&{\bar h}^{(1)}(x^n, n, \theta) = 
\frac{1}{m}{\bar h}_{ij}^{(1)}(x^n, n, \theta){\bar g}^{ij}(\theta) = 
\frac{1}{\nu}[h^{(1)}(x^n, n, \theta) - n^{(1)}(x^n, \theta)]
\in S(\bar{{\cal N}}_{\bar{M}}),\\
\label{eq:3-35}
&h^{(1)}(x^n, n, \theta) = 
\frac{1}{m}h_{ij}(x^n, n, \theta)g^{ij}(\theta) = 
\sum_{t=1}^n h^{(1)}(x(t), \theta)\in 
S(\bar{{\cal N}}_{\bar{M}}),
\end{align}
and the 1-ES umbilic curvature tensor variable 
${\bar k}_{ij}^{(1)}(x^n, n, \theta)$ is expressed as
\begin{align}
\label{eq:3-36}
{\bar k}_{ij}^{(1)}(x^n, n, \theta) &= {\bar h}_{ij}^{(1)}(x^n, n, \theta) - {\bar g}_{ij}(\theta){\bar h}^{(1)}(x^n, n, \theta)\nonumber \\
& = h_{ij}^{(1)}(x^n, n, \theta) - 
g_{ij}(\theta)h^{(1)}(x^n, n, \theta)\nonumber \\
& = k_{ij}^{(1)}(x^n, n, \theta) = 
\sum_{t=1}^n k_{ij}^{(1)}(x(t), \theta)\in S(\bar{{\cal N}}_{\bar{M}}).
\end{align}
For the nonnegative scalar curvatures, we have
\begin{align}
\label{eq:3-37}
&\bar{\lambda}^2 = \langle \bar{h}_{ij}^{(1)}(x^n, n, \theta)
\bar{h}_{kl}^{(1)}(x^n, n, \theta)\rangle_\theta \bar{g}^{ik}(\theta)
\bar{g}^{jl}(\theta) = \bar{\kappa}^2 + \bar{\gamma}^2,\\
\label{eq:3-38}
&\bar{\kappa}^2 = \langle \bar{k}_{ij}^{(1)}(x^n, n, \theta)
\bar{k}_{kl}^{(1)}(x^n, n, \theta)\rangle_\theta \bar{g}^{ik}(\theta)
\bar{g}^{jl}(\theta) = \frac{1}{\nu}\kappa^2,\\
\label{eq:3-39}
&\bar{\gamma}^2 = m\langle \bar{h}^{(1)}(x^n, n, \theta)^2\rangle_\theta = 
m\langle (h^{(1)}(x^n, n, \theta) - n^{(1)}(x^n, \theta))^2\rangle_\theta.
\end{align}
Note that we can always make $\bar{\gamma}^2 = 0$ by setting
\begin{align}
\label{eq:3-40}
h^{(1)}(x^n, n, \theta) = n^{(1)}(x^n, \theta),
\end{align} 
which is realized by the stopping rule (cf. Okamoto et al., 1991)
\begin{align}
\label{eq:3-41}
-\frac{1}{m}\partial_i\partial_j\bar{l}_n(x^n, \hat{\theta}_{mle,n})
g^{ij}(\hat{\theta}_{mle,n}) = \nu(\hat{\theta}_{mle,n}),
\end{align}
where $\hat{\theta}_{mle,n}$ denotes the maximum likelihood estimator of $\theta$ at the stopping time $\tau = n$, that is, 
$\partial_k\bar{l}_n(x^n, \hat{\theta}_{mle,n}) = 0$.

Then we obtain the following result as to the possibility that $\bar{M}$ would be an f.r.m.\ conformal exponential family $\bar{M}_e^*$.

\begin{theorem}
An extended statistical manifold 
$\bar{M} = \{\bar{f}(x^n, n, \theta)\ |\ \theta \in \Theta \}$ can be an f.r.m.\ conformal exponential family 
$\bar{M}_e^* = 
\{\bar{f}_e(z, \bar{\theta})\ |\ \bar{\theta} \in \bar{\Theta} \}$ 
if and only if $M$ is totally exponential umbilic in ${\cal R}_M$, that is,
\begin{align}
\label{eq:3-42}
k_{ij}^{(1)}(x, \theta) = h_{ij}^{(1)}(x, \theta) - 
g_{ij}(\theta)h^{(1)}(x, \theta) \equiv 0.
\end{align}
\end{theorem}
\

\begin{proof}
From Theorem 3.2, (\ref{eq:3-37}), (\ref{eq:3-38}), (\ref{eq:3-39}), 
(\ref{eq:3-40}), we have
\begin{align*}
&{\bar h}_{ij}^{(1)}(x^n, n, \theta) \equiv 0\ 
\Leftrightarrow\ \bar{\lambda}^2 = \bar{\kappa}^2 + \bar{\gamma}^2 \equiv 0
\ \Leftrightarrow\ \bar{\kappa}^2 = \frac{1}{\nu}\kappa^2 \equiv 0,\quad
\bar{\gamma}^2 \equiv 0\\
&\Leftrightarrow\ \kappa^2 \equiv 0\ \Leftrightarrow\ 
k_{ij}^{(1)}(x, \theta) \equiv 0.
\end{align*}

This completes the proof of the theorem.
\end{proof}
\

As noted by (\ref{eq:3-13}), an f.r.m.\ exponential family $M_e$ is 
totally $e$-umbilic in ${\cal R}_M$, so that from Theorem 3.3, the transformed $\bar{M}_e$ can be an f.r.m.\ conformal exponential family $\bar{M}_e^*$. 
Concretely we obtain the following result.

\begin{corollary}
An f.r.m.\ exponential family $M_e$ is preserved under the conformal transformation satisfying
\begin{align}
\label{eq:3-43}
n^{(1)}(x^n, \theta) = n - \nu - s^i\partial_i{\bar l}_n \equiv 0,\ 
\textrm{that is},\ 
n - \nu \in S(\bar{{\cal T}}_{\bar{M}}).
\end{align}
For a stopping time $\tau_c \in {\cal C}$ introduced in Section 2, 
the above is asymptotically satisfied in the sense 
\begin{align}
\label{eq:3-44}
&\langle n_c^{(1)}(x^n, \theta)^2\rangle_\theta = o(\nu_c(\eta))\ \ 
\textrm{as}\ \ c\ \to\ \infty,\\
&n_c^{(1)}(x^n, \theta) = n_c - \nu(c, \theta) - 
s^i(c, \eta)\partial_i{\bar l}_n(x^n, \theta), \nonumber \\
&\nu(c, \theta) = \langle n_c\rangle_\theta,\quad 
s^i(c, \eta) = s_j(c, \theta)g^{ji}(\eta),\quad 
s_j(c, \theta) = \partial_j \log \nu(c, \theta). \nonumber
\end{align}
\end{corollary}
\

\begin{proof}
We first prove the former part. 
When $M_e$ is an f.r.m.\ exponential family, from Theorem 3.1 and 
(\ref{eq:3-32}), we have 
\begin{align*}
h^{(1)}_{ij}(x^n, n, \theta) = 
\sum_{t=1}^n h^{(1)}_{ij}(x(t), \theta) \equiv 0.
\end{align*}
Hence from Theorem 3.2 and (\ref{eq:3-31}), it follows for $M_e$ that
\begin{align*}  
&n^{(1)}(x^n, \theta) = n - \nu - s^i\partial_i{\bar l}_n \equiv 0\\ 
&\Leftrightarrow\
\bar{h}_{ij}^{(1)}(x^n, n, \theta) = h_{ij}^{(1)}(x^n, n, \theta) - 
g_{ij}(\theta)n^{(1)}(x^n, \theta) \equiv 0\\
&\Leftrightarrow\ 
\textrm{$\bar{M}_e$ is an f.r.m.\ conformal exponential family $\bar{M}_e^*$},
\end{align*}
which proves the former part.

We next prove the latter part. 
When $M_e$ is an f.r.m.\ exponential family with the natural parameter $\theta$ and the expectation parameter 
$\eta = \partial \psi(\theta)/\partial \theta$, by the direct decomposition 
in (\ref{eq:3-33})
\begin{align*}
n_c - \nu(c, \theta) = s^i(c, \eta)\partial_i\bar{l}_n(x^n, \theta) + 
n_c^{(1)}(x^n, \theta),\quad 
s^i(c, \eta)\partial_i\bar{l}_n(x^n, \theta) \in 
S(\bar{{\cal T}}_{\bar{M}}),\quad 
n_c^{(1)}(x^n, \theta) \in S(\bar{{\cal N}}_{\bar{M}}),
\end{align*}
we have
\begin{align*}
\textrm{Var}_\theta(\tau_c)& = 
\langle \left(n_c - \nu(c, \theta)\right)^2\rangle_\theta = 
\langle s^i(c, \eta)s^j(c, \eta)\partial_i \bar{l}_n(x^n, \theta) 
\partial_j \bar{l}_n(x^n, \theta)\rangle_\theta + 
\langle n_c^{(n)}(x^n, \theta)^2\rangle_\theta\\
&= \nu(c, \theta) g_{ij}(\theta)s^i(c, \eta)s^j(c, \eta) + 
\langle n_c^{(1)}(x, \theta)^2\rangle_\theta.
\end{align*}
From Proposition 2.1 with 
$\sigma_{ij}(\theta) = g_{ij}(\theta)$, it follows that
\begin{align*}
&\frac{\nu(c, \theta)}{\nu_c(\eta)}\ \to\ 1,\quad 
\frac{\textrm{Var}_\theta(\tau_c)}
{\nu_c(\eta) g_{ij}(\theta)s_c^i(\eta)s_c^j(\eta)}\ \to\ 1,\quad
\frac{\nu(c, \theta) g_{ij}(\theta)s^i(c, \eta)s^j(c, \eta)}
{\nu_c(\eta) g_{ij}(\theta)s_c^i(\eta)s_c^j(\eta)}\ \to\ 1
\ \ \textrm{as}\ \ c\ \to\ \infty \\
&\Rightarrow\ 
\frac{\textrm{Var}_\theta(\tau_c)}
{\nu(c, \theta) g_{ij}(\theta)s^i(c, \eta)s^j(c, \eta)}\ 
\to\ 1 \\
&\Rightarrow\ 
\langle n_c^{(1)}(x^n, \theta)^2\rangle_\theta = 
o(\textrm{Var}_\theta(\tau_c)) = 
o(\nu_c(\eta)),
\end{align*}
which proves the latter part.
\end{proof}
\

The relation (\ref{eq:3-43}) implies that the stopping time variable $n$ 
and the $m$ score functions 
$\partial_i\bar{l}_n = \theta^i \sum_{t=1}^n x_i(t) - n \partial_i\psi\ 
(i = 1, \dots, m)$ must be affinely dependent, in other words, there exist constants 
$c^i\ (i = 0, 1, \dots, m),\ d$ such that
\begin{align}
\label{eq:3-45}
c^0 n + c^i \sum_{t=1}^n x_i(t) = d,\quad 
\exists c^i \neq 0,\ i\in \{0, 1, \dots, m\},\ 
d\neq 0.
\end{align}
We note that the property (\ref{eq:3-45}) is the same as that given in Winkler and Franz (1979), which was derived from the statistical considerations of the efficient sequential estimators attaining the 
Cram\'er-Rao bound.

By definitions, the equation of Gauss (\ref{eq:3-28}) with $\alpha = 1$, Theorems 3.2 and 3.3, we also know the following implications.
\begin{align}
\label{eq:3-46}
&\textrm{$M$ is an f.r.m.\ exponential family $M_e$}\nonumber \\
&\Rightarrow\ \textrm{$M$ is totally $e$-umbilic in ${\cal R}_M$}\nonumber \\ 
&\Leftrightarrow\ 
\textrm{$\bar{M}$ can be an f.r.m.\ conformal exponential family $\bar{M}_e^*$}\nonumber \\ 
&\Rightarrow\ \textrm{$M$ is conformally $m(e)$-flat}. 
\end{align}

\section{SUMMARY OF ASYMPTOTICS ON NONSEQUENTIAL TEST}

In this section, we summarize preliminary results on the nonsequential test in a multidimensional curved exponential family which will be necessary for analyzing the sequential test.

\subsection{Curved exponential family}

We first introduce a curved exponential family. 
A family of probability densities $M_c = \{f_c(x, u)\ |\ u\in U\}$ parameterized by an $m$-dimensional vector parameter 
$u = (u^1, \dots, u^m)^{\tp}$ is said to be an $(n, m)$-curved exponential family when it is smoothly imbedded in an $n$-dimensional f.r.m.\ exponential family $M_e = \{f_e(x, \theta)\ |\ \theta\in \Theta\}$ in the sense 
\begin{align}
\label{eq:4-1}
f_c(x, u) = f_e[x, \theta(u)] = 
\exp\{\theta^i(u)x_i - \psi[\theta(u)]\},
\end{align}
where $U$ is homeomorphic to ${\mathbb R}^m\ (m < n)$ and the parametric representations 
$\theta(u) = [\theta^1(u), \dots, \theta^n(u)]^{\tp},\ 
\eta(u) = \partial \psi(\theta)(u) = [\eta_1(u), \dots,\eta_n(u)]^{\tp}$ are smooth functions of $u$ having full-rank Jacobian matrices. 
We use indices $i, j, k$ and so on to denote quantities in terms of the coordinate system $\theta$ (upper indices) or $\eta$ (lower indices) of $M_e$, and indices $a, b, c$ and so on to denote quantities in terms of the coordinate system $u$ of $M_c$.

\subsection{Asymptotic powers of tests}

In an $(n, m)$-curved exponential family $M_c$ we consider an unbiased test of the simple null hypothesis
\begin{align*}
H_0 : u = u_0\quad \textrm{vs}\quad H_1 : u \neq u_0.
\end{align*}
Let $x_1, \dots, x_N$ be $N$ independent observations from a distribution 
in $M_c$. 
The hypothesis $H_0$ is tested based on the sufficient statistic 
${\bar x} = \sum_{t=1}^N x(t)/N$. 
Let $R$ be the critical region of a test ${\cal T}$. 
Then the power of ${\cal T}$ at $u$ is written as
\begin{align*}
P_{{\cal T}}(u) = \int_R f({\bar x}, u)d\mu({\bar x}),
\end{align*}
where $f({\bar x}, u)$ denotes the density function of ${\bar x}$ when the true parameter is $u$.

In an asymptotic theory, we are treating a test sequence 
${\cal T}_N, N = 1, 2, \dots$, and to appropriately evaluate the power of $\{{\cal T}_N\}$, we define a set of alternatives $U_N(s)$ by
\begin{align*}
U_N(s) = \{u_{s,e} = u_0 + se/\sqrt{N}\in M_c,\ 
s \ge 0,\ g_{ab}(u_0)e^ae^b = 1 \},
\end{align*}
where $[g_{ab}]$ denotes the Riemannian metric on $M_c$ given by the Fisher information matrix, and $e = (e^a)$ denotes the $m$-dimensional unit tangent vector of $M_c$ at $u_0$. 
The set $U_N(s)$ consists of the points of $M_c$ which are separated from $u_0$ by the geodesic distance $s/\sqrt{N}$ (Figure 1).

\begin{figure}[htbp]
\begin{center}
\includegraphics[width=10cm,height=8cm]{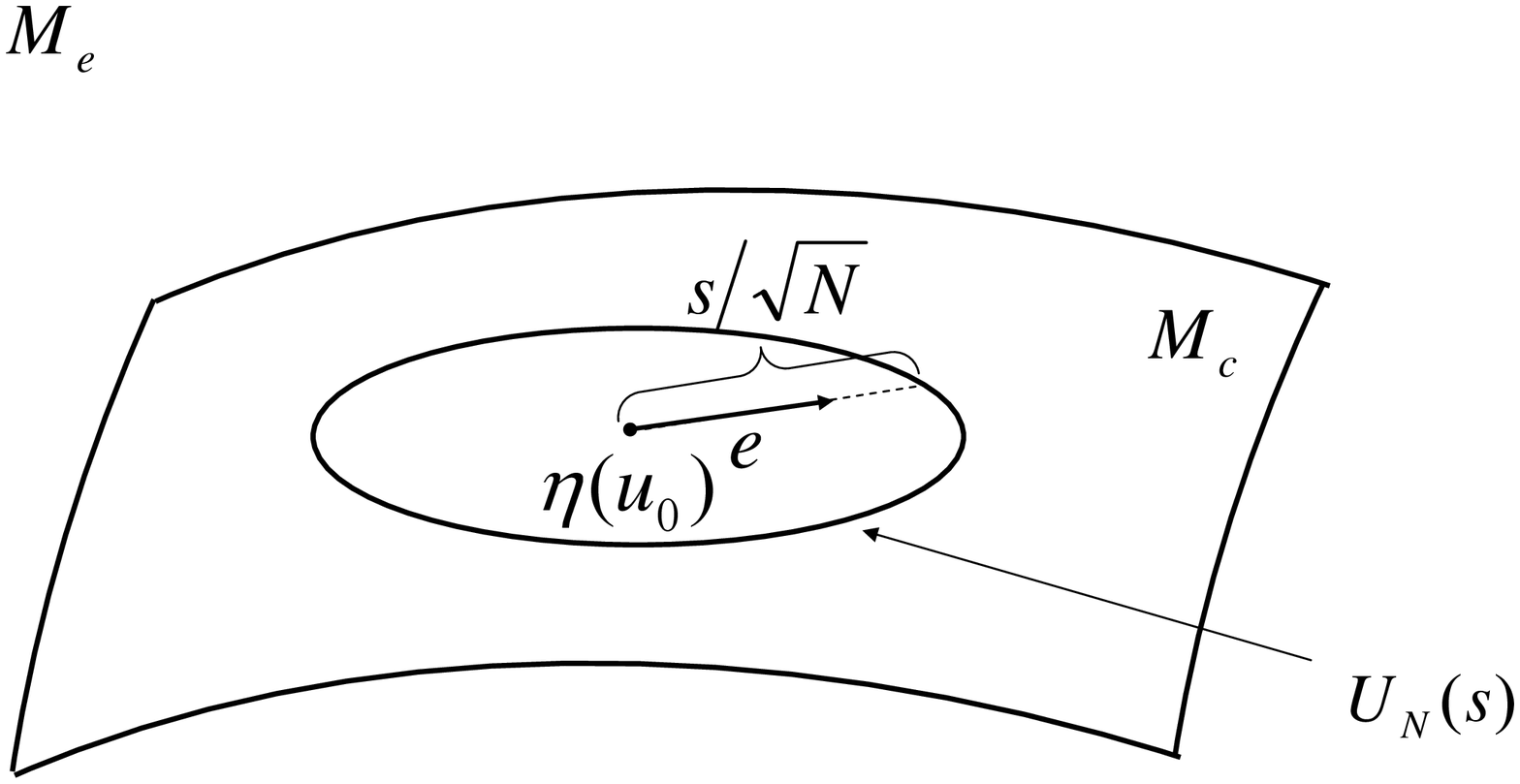}
\vspace*{-16mm}   
\caption{The set of alternatives $U_N(s)$ associated with an unbiased test of $u = u_0$.}
\label{fig:1-1}
\end{center}
\end{figure}

We compare tests based on the average power $P_{\cal T}(s)$ defined by
\begin{align*}
P_{{\cal T}}(s) = \langle P_{{\cal T}}(u)\rangle_e = 
\int_{u\in U_N(s)} P_{{\cal T}}(u)du/S_N(s),
\end{align*}
where $\langle\ \rangle_e$ means the average over all directions of $e$, so that $S_N(s)$ denotes the area of $U_N(s)$.

Tests are required to satisfy the level and the unbiasedness conditions. 
We say a test ${\cal T}$ is of $i$-th order level $\alpha$, when 
\begin{align*}
P_{{\cal T}}(0) = \alpha + O(N^{-i/2}),
\end{align*}
and we note that the unbiasedness condition 
\begin{align*}
\frac{d}{ds}P_{{\cal T}}(s)\Big|_{s=0} = 0
\end{align*}
is automatically satisfied, since the power $P_{{\cal T}}(s)$ is obtained by taking the average in all directions.

By expanding $P_{{\cal T}}(s)$ in the power series of $N^{-1/2}$, we have
\begin{align*}
P_{{\cal T}}(s) = P_{{\cal T}1}(s) + 
P_{{\cal T}2}(s)N^{-1/2} + P_{{\cal T}3}(s)N^{-1} + O(N^{-3/2}),
\end{align*}
and call $P_{{\cal T}i}(s)$ the $i$-th order power of a test ${\cal T}$ at $s$. 
When discussing the $i$-th order power, the tests are assumed to satisfy the level condition up to the same $i$-th order.

A test ${\cal T}$ is said to be first-order uniformly efficient or said simply to be efficient in short when
\begin{align*}
P_{{\cal T}1}(s) \ge P_{{\cal T}'1}(s),\quad 
\forall s \ge 0,\quad \forall {\cal T}'.
\end{align*}
It can be shown that an efficient test is automatically second-order uniformly efficient, and that there does not in general exist a third-order uniformly efficient test. 
Thus we say that an efficient test ${\cal T}$ is third-order 
$s_0$-efficient when
\begin{align*}
P_{{\cal T}3}(s_0) \ge P_{{\cal T}'3}(s_0),\quad  
\forall \textrm{efficient ${\cal T}'$}.
\end{align*}

In order to evaluate the power loss of a test, we next introduce the envelope power function by 
$P^*(s) = \sup_{{\cal T}} P_{{\cal T}}(s)$, where $\sup$ is taken at each $s$ separately. 
We expand it as 
\begin{align*}
P^*(s) = P_1^*(s) + P_2^*(s)N^{-1/2} + P_3^*(s)N^{-1} + O(N^{-3/2}),
\end{align*}
and call $P_i^*(s)$ the $i$-th order envelope power at $s$. 
Then the (third-order) power-loss function $\Delta P_{{\cal T}}(s)$ of an efficient test ${\cal T}$ is defined by
\begin{align}
\label{eq:4-2}
\Delta P_{{\cal T}}(s) = P_3^*(s) - P_{{\cal T}3}(s).
\end{align}

\subsection{Ancillary family associated with a test}

In calculating the power function $P_{{\cal T}}(s)$ of a test ${\cal T}$, it is convenient to associate an ancillary family with a test. 
We attach to each point $u\in M_c$ an $(n-m)$-dimensional smooth submanifold $A(u)$ of $M_e$ which transverses $M_c$ at $\theta(u)$ 
or equivalently at $\eta(u)$. 
Such an $A(u)$ is called an 
\textit{ancillary submanifold}, and the union 
$A = \{A(u)\ |\ u\in M_c\}$ is called an 
\textit{ancillary family}. 
Let $R_{M_c} = R\cap M_c$ be the intersection of $M_c$ and the critical region $R$ of a test ${\cal T}$. 
An ancillary family $A$ is said to be associated with a test ${\cal T}$ when
\begin{align}
\label{eq:4-3}
R = \underset{u\in R_{M_c}}\cup A(u) = 
\{\eta\ |\ \eta \in A(u),\ u \in R_{M_c}\}.
\end{align}

\begin{figure}[htbp]
\begin{center}
\includegraphics[width=10cm,height=8cm]{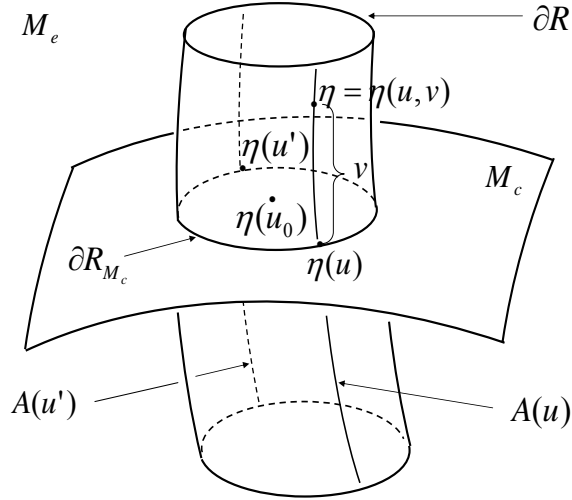}
\vspace*{-5mm}   
\caption{Ancillary family $A = \{A(u)\}$ and $w = (u, v)$-coordinates 
associated with an unbiased test of $u = u_0$.}
\label{fig:1-1}
\end{center}
\end{figure}

By introducing a coordinate system $v = (v^{m+1}, \dots, v^n)$ to each $A(u)$ such that the origin $v = 0$ is at the intersection of $A(u)$ and $M_c$, we have a new local coordinate system of $M_e$ (Figure 2)
\begin{align}
\label{eq:4-4}
w = (w^\alpha) = (u^a, v^\kappa),\quad \alpha = 1, \dots, n,\quad 
a = 1, \dots, m,\quad \kappa = m+1, \dots, n,
\end{align}
and we use indices $\alpha, \beta, \gamma$ and so on for quantities related to the coordinate system $w$, and indices $\kappa, \lambda, \mu$ and so on for quantities related to the coordinate system $v$.

The critical region $R$ is written as
\begin{align*}
R = \{w = (u, v)\ |\ u\in R_{M_c},\ v\ \textrm{is arbitrary}\}
\end{align*}
in the new coordinate system $w$. 
The sufficient statistic ${\bar x}$ is transformed to 
${\hat w} = ({\hat u}, {\hat v})$ by ${\bar x} = \eta({\hat w})$, 
and the power $P_{{\cal T}}(u)$ is written as
\begin{align*}
P_{{\cal T}}(u) = \int_R f({\hat w}, u)d{\hat w} = 
\int_{R_{M_c}}\int_{A(u)} f({\hat w}, u)d{\hat v}d{\hat u} 
= \int_{R_{M_c}} f({\hat u}, u)d{\hat u},
\end{align*}
where $f({\hat u}, u) = \int_{A(u)} f({\hat w}, u)d{\hat v}$ is the density function of ${\hat u}$ when the true parameter is $u$.

\subsection{Metric, connections and curvatures}

The basic tensors of $M_e$ are written as
\begin{align*}
g_{\alpha \beta}(w) = B_\alpha^iB_\beta^jg_{ij}(\theta),\quad 
T_{\alpha \beta \gamma}(w) = 
B_\alpha^iB_\beta^jB_\gamma^kT_{ijk}(\theta),\quad 
B_\alpha^i = \frac{\partial \theta^i}{\partial w^\alpha},
\end{align*}
and the $\alpha$-connection is given by 
\begin{align*}
\Gamma_{\beta \gamma \delta}^{(\alpha)}(w) &= 
B_\beta^iB_\gamma^jB_\delta^k\Gamma_{ijk}^{(\alpha)}(\theta) + 
B_\delta^i\partial_\beta B_\gamma^jg_{ij}(\theta) = 
\frac{1-\alpha}{2}T_{\beta \gamma \delta} + 
(\partial_\beta B_\gamma^j)B_{\delta j}\\
&= B_{\beta i}B_{\gamma j}B_{\delta k}\Gamma^{(\alpha)ijk}(\eta) + 
B_{\delta i}\partial_\beta B_{\gamma j}g^{ij}(\eta) = 
-\frac{1+\alpha}{2}T_{\beta \gamma \delta} + 
(\partial_\beta B_{\gamma j})B_\delta^j,\\
B_{\beta i} &= \frac{\partial \eta_i}{\partial w^\beta} = g_{ij}B_\beta^j,
\end{align*}
in the $w$-coordinate system. 
When we evaluate a quantity $q(u, v)$ on $M_c$, that is, at $v = 0$, we often denote it by $q(u)$ instead of by $q(u, 0)$ for brevity's sake. 
The metric tensors of $M_c$ and $A(u)$ are given by
\begin{align*}
g_{ab}(u) = B_a^iB_b^jg_{ij} = B_{ai}B_{bj}g^{ij},\quad 
g_{\kappa\lambda}(u) = B_\kappa^iB_\lambda^jg_{ij} = 
B_{\kappa i}B_{\lambda j}g^{ij},
\end{align*}
and then indices can be lowered or raised by using these metric tensors or their inverses 
$g^{ab}(u), g^{\kappa \lambda}(u)$. 
The mixed part
\begin{align*}
g_{a\kappa}(u) = B_a^iB_\kappa^jg_{ij} = B_{ai}B_{\kappa j}g^{ij}
\end{align*}
is the inner product of the tangent vectors $B_a$ of $M_c$ and the tangent vectors $B_\kappa$ of $A(u)$, which represents the angles between $M_c$ and $A(u)$. 
When $g_{a\kappa}(u_0) = 0$ holds at a point $u_0$, the ancillary family is said to be locally orthogonal at $u_0$. 
When $g_{a\kappa}(u) = 0$ holds at every $u$, the ancillary family is said to be an orthogonal family.

For a locally orthogonal ancillary family, the inner product 
$g_{a\kappa}(u_{s,e})$ at $u_{s,e} = u_0 + se/\sqrt{N}$ is expanded as
\begin{align}
\label{eq:4-5}
g_{a\kappa}(u_{s,e}) = Q_{ba\kappa}e^bs/\sqrt{N} + O(N^{-1}),\quad 
Q_{ba\kappa} = \partial _bg_{a\kappa}(u_0),
\end{align}
where $Q_{ba\kappa}$ is the tensor quantity characterizing the angles between $\cup_{u_{s,e}\in U_N(s)}A(u_{s,e})$ and $M_c$.

When the ancillary family is locally orthogonal, 
the mixed parts $\Gamma_{ab\kappa}^{(1)}(u_0)$ and 
$\Gamma_{\kappa\lambda a}^{(-1)}(u_0)$ of the $\pm1$-connections are tensor quantities, which are denoted by
\begin{align}
\label{eq:4-6}
H_{ab\kappa}^{(1)}(u_0) = \Gamma_{ab\kappa}^{(1)}(u_0) = 
(\partial_a B_b^i)B_{\kappa i},\quad
H_{\kappa\lambda a}^{(-1)}(u_0) = \Gamma_{\kappa\lambda a}^{(-1)}(u_0) = 
(\partial_\kappa B_{\lambda i})B_a^i.
\end{align}
These tensor quantities $H_{ab\kappa}^{(1)}(u_0)$ and 
$H_{\kappa\lambda a}^{(-1)}(u_0)$ are called the $\pm1$-ES curvature tensors of $M_c$ and $A(u_0)$ in $M_e$, respectively. 
When $H_{\kappa\lambda a}^{(-1)}(u_0) = 0$ holds, the ancillary family 
$A$ is said to be 
\textit{mixture flat ($m$-flat)} at $u_0$ in $M_e$.

Related to the $1$-ES curvature tensor $H_{ab\kappa}^{(1)}$ of $M_c$ in $M_e$, we also define
\begin{align}
\label{eq:4-7}
H_\kappa^{(1)} = \frac{1}{m}H_{ab\kappa}^{(1)}g^{ab},\quad
K_{ab\kappa}^{(1)} = H_{ab\kappa}^{(1)} - g_{ab}H_\kappa^{(1)}.
\end{align}
The first $H_\kappa^{(1)}$ is called the mean $1$-ES curvature tensor of $M_c$ in $M_e$, and the second $K_{ab\kappa}^{(1)}$ is called the $1$-ES umbilic tensor of 
$M_c$ in $M_e$, respectively. 
When $K_{ab\kappa}^{(1)} \equiv 0$ holds, $M_c$ is said to be totally exponential umbilic ($e$-umbilic) in $M_e$.

We remark that the 1-ES curvature tensor variable $h_{ab}^{(1)}(x, u)$ 
of $M_c$ in ${\cal R}_{M_c}$ can be also introduced and expressed as
\begin{align}
\label{eq:4-8}
h_{ab}^{(1)}(x, u) &= \partial_a \partial_b l_1(x, u) + g_{ab}(u) - 
\Gamma_{ab}^{(1)c}(u)\partial_c l_1(x, u) \nonumber \\
&= H_{ab\kappa}^{(1)}(u)g^{\kappa \lambda}(u)\partial_\lambda l_1(x, u),\quad l_1(x, u) = \log f_c(x, u).
\end{align}
This expression shows that the 1-ES curvature tensor 
$H_{ab\kappa}^{(1)}(u)$ of $M_c$ in $M_e$ and the 1-ES curvature tensor variable $h_{ab}^{(1)}(x, u)$ of $M_c$ in ${\cal R}_{M_c}$ are the equivalent notions.
Hence when $K_{ab\kappa}^{(1)} \equiv 0$ holds, $M_c$ can also said to be totally $e$-umbilic in $R_{M_c}$.

\subsection{First-order power}

The first-order analysis is based on the normal distribution on 
$f({\hat w}, u)$, and we have the following result 
(cf. Amari, 1985, Theorem 6.17).

\begin{proposition}
A test ${\cal T}$ is first-order uniformly efficient, if and only if the associated ancillary family is asymptotically orthogonal in the sense
\begin{align*}
g_{a\kappa}(u) = O(N^{-1/2}),\quad \forall u\in 
\partial R_{M_c} = \partial R \cap M_c,
\end{align*}
and the optimal intersection $R_{M_c}^* = R^* \cap M_c$ is given by
\begin{align*}
R_{M_c}^* = \{{\tilde u}_0\ |\ 
g_{ab}(u_0){\tilde u}_0^a{\tilde u}_0^b \ge c_0^2 \},\quad  
{\tilde u}_0 = \sqrt{N}({\hat u} - u_0),\quad
c_0^2 = \chi_{m,\alpha}^2,
\end{align*}
where $\chi_{m,\alpha}^2$ denotes the upper $100\alpha\%$ point of the $\chi^2$-distribution with $m$ degrees of freedom.
The first-order envelope power function $P_1^*(s)$ is expressed as
\begin{align*}
&P_1^*(s) = 1 - \int_0^{c_0} Z_m^{(0)}(s, r)dr,\\
&Z_m^{(0)}(s, r) = (2\pi)^{-m/2}r^{m-1}S_{m-2}\exp 
\bigg\{-\frac{1}{2}(s^2+r^2)\bigg\}A_k^{(0)}(s, r),\\
&A_k^{(0)}(s, r) = \int_0^\pi \sin^{m-2}\varphi \exp(sr\cos\varphi)d\varphi,
\end{align*}
where $S_{m-2}$ denotes the area of the $(m-2)$-dimensional unit sphere with $S_0 = 2$.
Note that $Z_m^{(0)}(s, r)$ is the density function of the 
$\chi^2$-distribution with $m$ degrees of freedom and non-centrality parameter $s^2$.
\end{proposition}

%\begin{figure}[htbp]
%\begin{center}
%\includegraphics[width=8cm,height=7cm]{ConTP1-2to5N.eps}
%\vspace*{-5mm}   
%\caption{First order envelope power function $P_1^*(s)$.}
%\label{fig:1-1}
%\end{center}
%\end{figure}
%
%
%Figure 3 shows the first order envelope power function $P_1^*(s)$ given by 
%(\ref{eq:4-22}) in the cases of $m = 2, 3, 4, 5$ with $\alpha = 0.05$.  
%

\subsection{Second- and third-order powers}

These powers can be calculated by taking account of the higher-order terms of $f({\hat w}, u)$, and by assuming the intersection $R_{M_c}$ as
\begin{align}
\label{eq:4-9}
R_{M_c} = \{{\tilde u}_0\ |\ 
g_{ab}(u_0){\tilde u}_0^{*a}{\tilde u}_0^{*b} 
\ge (c_0 + \epsilon)^2 \},\quad 
{\tilde u}_0^{*a} = {\tilde u}_0^a +  
\Gamma_{\beta \gamma}^{(-1)a}g^{\beta \gamma}({\hat u})/(2\sqrt{N}),
\end{align}
where ${\tilde u}_0^*$ is the bias-corrected version of ${\tilde u}_0$, 
and $\epsilon$ is introduced to adjust to the level condition up to $O(N^{-1})$. It is shown that $\epsilon = O(N^{-1})$, which depends on the ancillary family through $H_{\kappa\lambda a}^{(-1)}$ and $Q_{ab\kappa}$. 
Hence it follows that the second-order power function 
$P_{{\cal T}2}(s)$ is common to all the first-order efficient tests 
(see Appendices 1, 2).

Here we introduce a class of efficient tests whose ancillary family satisfies
\begin{align}
\label{eq:4-10}
Q_{(ab)\kappa} = k_1K_{ab\kappa}^{(1)} + k_2g_{ab}H_\kappa^{(1)},\quad
Q_{(ab)\kappa} = \frac{1}{2}(Q_{ab\kappa} + Q_{ba\kappa}),\quad
k_1, k_2 \in {\mathbb R}, 
\end{align}
and a test ${\cal T}$ with a pair of proportions $(k_1, k_2)$ is said to be the $k = (k_1, k_2)$-test.
Then after cumbersome calculations, we have the following result on the third-order efficiency of a $k = (k_1, k_2)$-test, which is the correction of Theorem 6.18 in Amari (1985).

\begin{proposition}
The third-order power loss function $\Delta P_{\cal K}(s)$ of a 
$k = (k_1, k_2)$-test ${\cal T}$ is given by
\begin{align}
\label{eq:4-11}
\Delta P_{\cal K}(s) &= \kappa^2 \Delta P_1(s) + \gamma^2 \Delta P_2(s)
+ H_A^{(-1)2}\xi_{0m}(s),
\end{align}
where
\begin{align}
\label{eq:4-12}
&\Delta P_1(s) = \xi_{1m}(s)\{K_1 - K_{1m}(s)\}^2,\quad
\Delta P_2(s) = \xi_{2m}(s)\{K_2 - K_{2m}(s)\}^2,\\
&\kappa^2 = K_{ab\kappa}^{(1)}K_{cd\lambda}^{(1)}
g^{ac}g^{bd}g^{\kappa \lambda},\quad 
\gamma^2 = mH_{\kappa}^{(1)}H_{\lambda}^{(1)}g^{\kappa \lambda}, \nonumber\\
&H_A^{(-1)2} = H_{\kappa \lambda a}^{(-1)}H_{\mu \nu b}^{(-1)}
g^{\kappa \mu}g^{\lambda \nu}g^{ab}. \nonumber
\end{align}
An efficient $k = (k_1, k_2)$-test ${\cal T}$ is third-order $s_0$-efficient if and only if
\begin{align}
\label{eq:4-13}
k_1 = K_{1m}(s_0),\quad k_2 = K_{2m}(s_0),\quad 
H_{\kappa \lambda a}^{(-1)}(u_0) = 0.
\end{align}
In the above expressions,
\begin{align}
\label{eq:4-14}
K_{1m}(s_0) = J_{1m}(s_0),\quad
K_{2m}(s_0) = J_{2m}(s_0) + 
\frac{\xi_{1m}(s_0)\{J_{1m}(s_0)-J_{2m}(s_0)\}}{\xi_{2m}(s_0)},
\end{align}
and
\begin{align*}
&\xi_{1m}(s) = \frac{1}{2m(m+2)}\left[2c_0^2sZ_m^{(1)}(s, c_0) - 
2c_0s^2\{Z_m^{(0)}(s, c_0)-Z_m^{(2)}(s, c_0)\} + 
\int_0^{c_0} f_{1m}(s, r) dr\right],\\
&f_{1m}(s, r) = -msrZ_m^{(1)}(s, r) + 
s^2\{(m + 2)Z_m^{(0)}(s, r) - 2kZ_k^{(2)}(s, r)\} + 
2s^3r\{Z_k^{(1)}(s, r) - Z_m^{(3)}(s, r) \},\\
&\xi_{2m}'(s) = \frac{1}{2m(m+2)}\left[c_0^2sZ_m^{(1)}(s, c_0) +  
c_0s^2\{Z_m^{(0)}(s, c_0)-Z_m^{(2)}(s, c_0)\} + 
\int_0^{c_0} f_{2m}(s, r) dr\right],\\
&f_{2m}(s, r) = srZ_m^{(1)}(s, r) - 
s^2\{2Z_m^{(0)}(s, r) - mZ_m^{(2)}(s, r)\} - 
s^3r\{Z_m^{(1)}(s, r) - Z_m^{(3)}(s, r) \},\\
&J_{1m}(s) = \frac{\xi_{3m}(s)}{2\xi_{1m}(s)},\quad 
J_{2m}(s) = \frac{\xi_{4m}(s)}{2\xi_{2m}'(s)},\quad
\xi_{2m}(s) = \xi_{1m}(s) + m\xi_{2m}'(s),\\
&\xi_{3m}(s) = \frac{1}{2m(m+2)}\left[4c_0^2sZ_m^{(1)}(s, c_0) - 
2c_0s^2Z_m^{(0)}(s, c_0) + 
\int_0^{c_0} f_{3m}(s, r) dr\right],\\
&f_{3m}(s, r) = -2msrZ_m^{(1)}(s, r) + 2ms^2Z_m^{(0)}(s, r),\\
&\xi_{4m}(s) = \frac{1}{2m(m+2)}\left[2c_0^2sZ_m^{(1)}(s, c_0) - 
c_0s^2Z_m^{(0)}(s, c_0) + 
\int_0^{c_0} f_{4m}(s, r) dr\right],\\
&f_{4m}(s, r) = 2srZ_m^{(1)}(s, r) - 2s^2Z_m^{(0)}(s, r),\\
&\xi_{0m}(s) = \frac{1}{4m}\int_0^{c_0}f_{0m}(s, r)dr,\\
&f_{0m}(s, r) = -srZ_m^{(1)}(s, r) + s^2Z_m^{(0)}(s, r),\\
&Z_m^{(l)}(s, r) = (2\pi)^{-m/2}r^{m-1}S_{m-2}\exp 
\left\{-\frac{1}{2}(s^2+r^2)\right\}A_m^{(l)}(s,r),\\
&A_m^{(l)}(s, r) = \int_0^\pi \sin^{m-2}\varphi \cos^l\varphi 
\exp(sr\cos \varphi)d\varphi.
\end{align*}
\end{proposition}
\

\clearpage

\begin{figure}[htbp]
\begin{center}
\includegraphics[width=8cm,height=7cm]{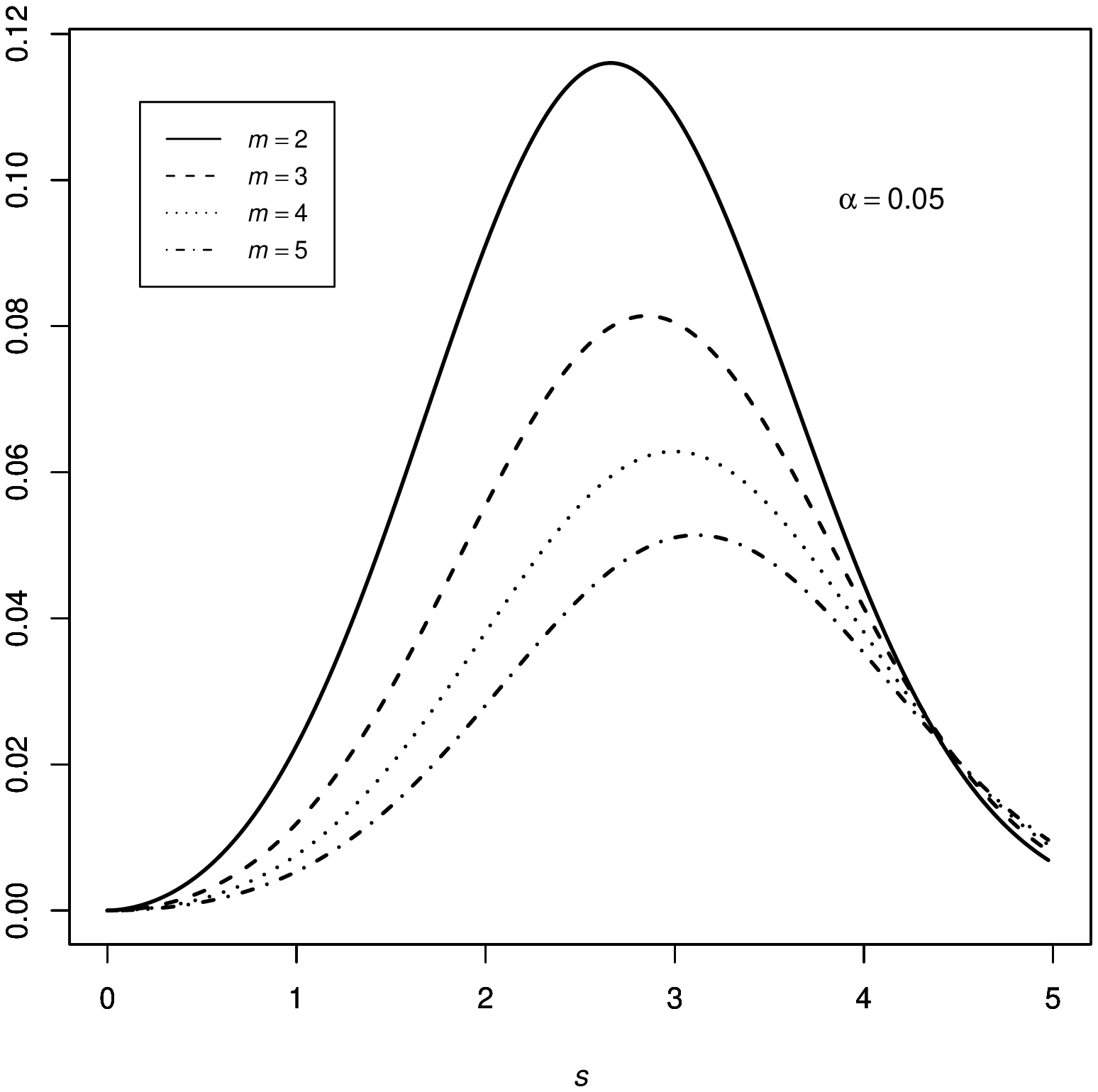}
\vspace*{-5mm}   
\caption{Power loss coefficient $\xi_{0m}(s)$.}
\label{fig:1-1}
\end{center}
\begin{minipage}{.5\linewidth}
\begin{center}
\includegraphics[width=8cm,height=7cm]{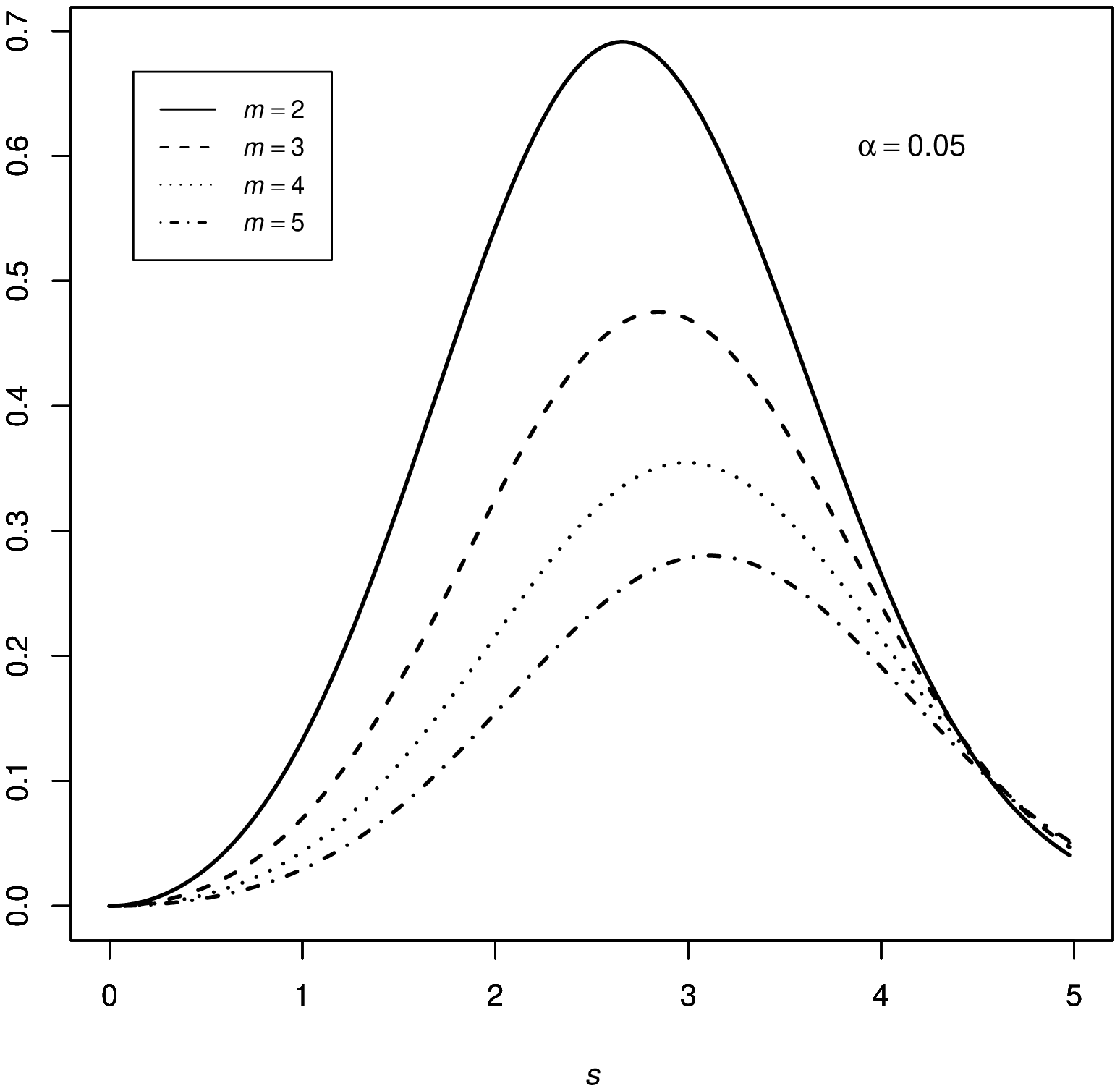}
\vspace*{-5mm}   
\caption{Power loss coefficient $\xi_{1m}(s)$.}
\label{fig:1-1}
\end{center}
\end{minipage}
\begin{minipage}{.5\linewidth}
\begin{center}
\includegraphics[width=8cm,height=7cm]{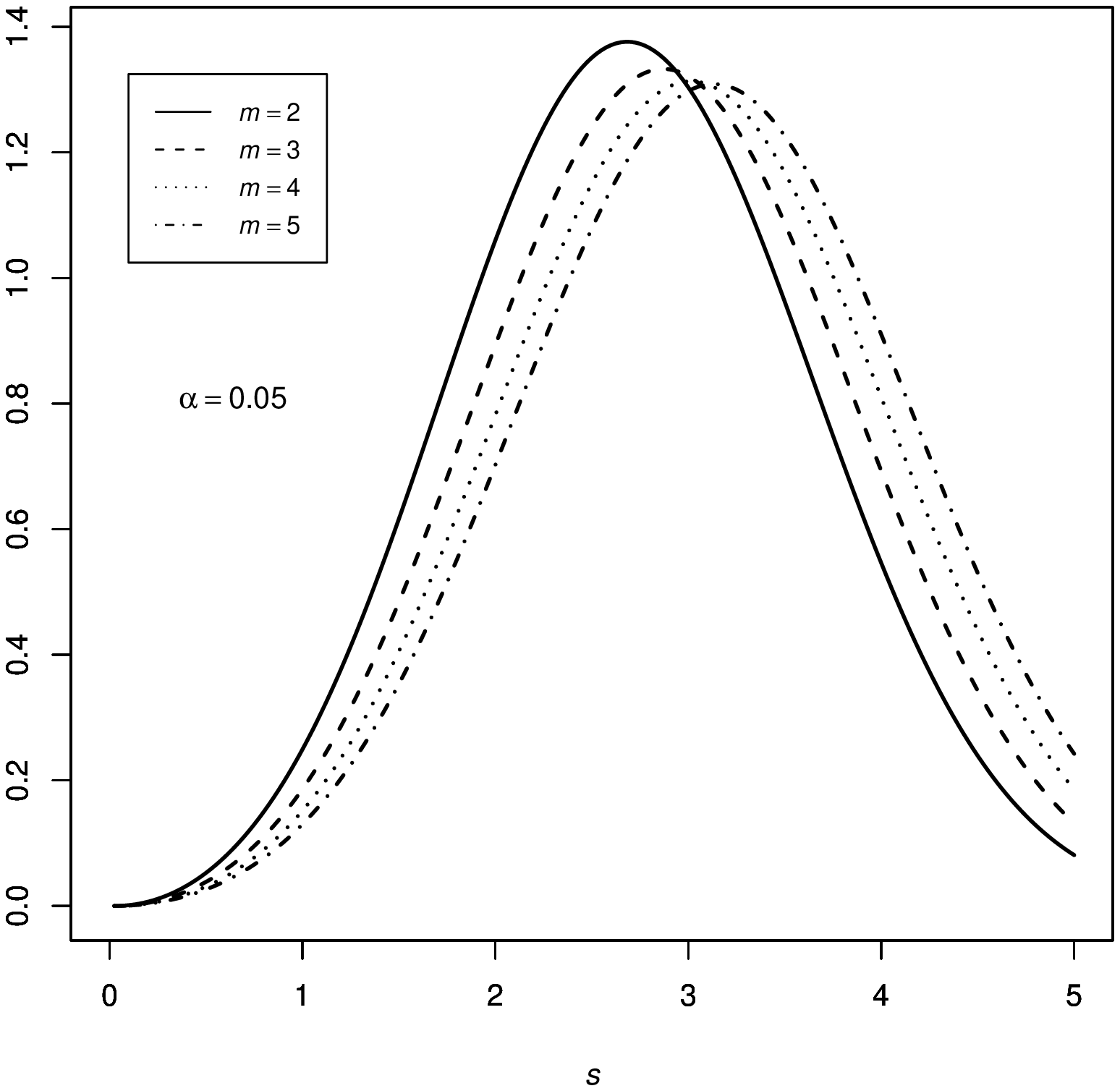}
\vspace*{-5mm}   
\caption{Power loss coefficient $\xi_{2m}(s)$.}
\label{fig:1-1}
\end{center}
\end{minipage}
\end{figure}

\begin{figure}[htbp]
\begin{minipage}{.5\linewidth}
\begin{center}
\includegraphics[width=8cm,height=7cm]{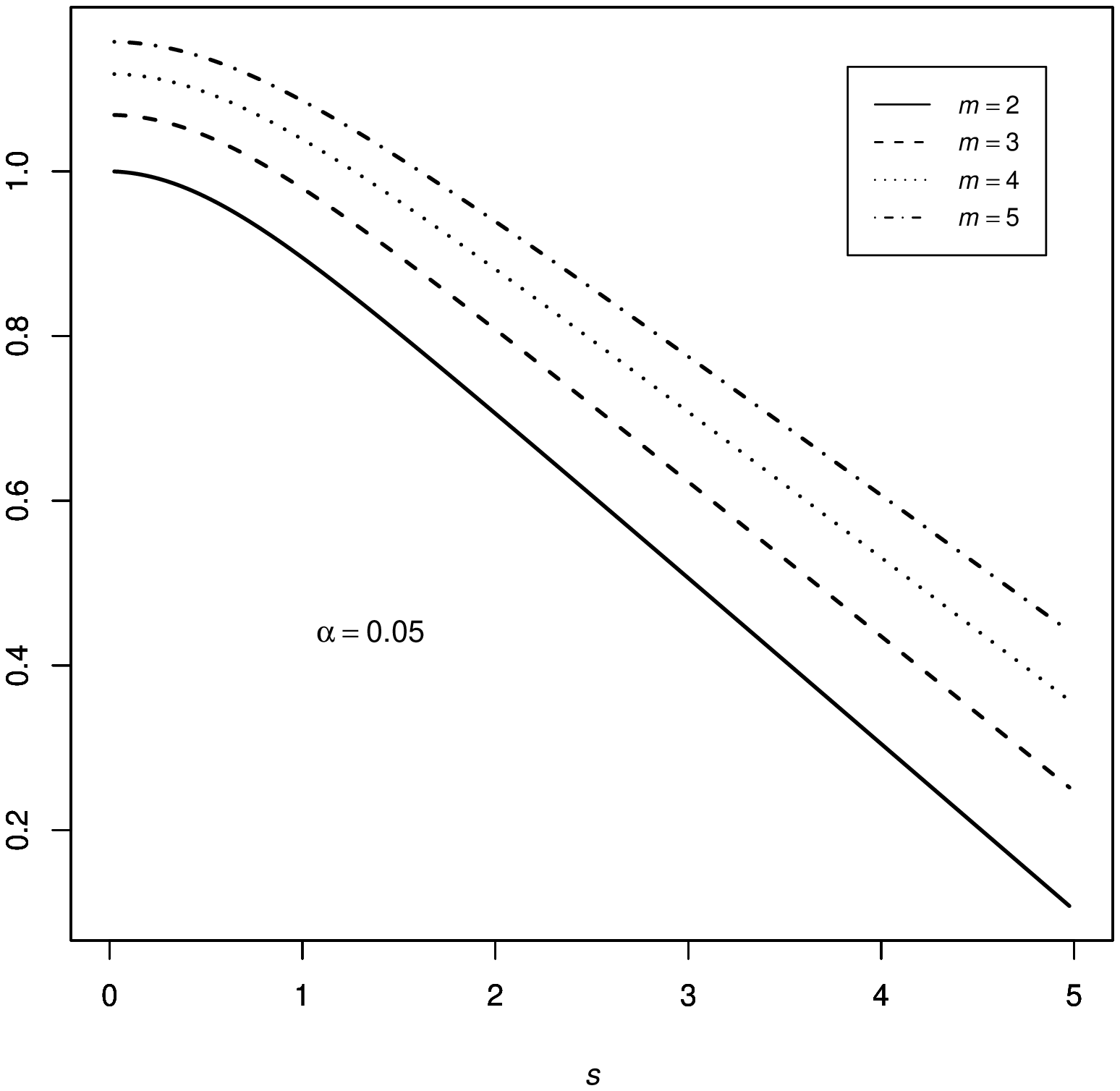}
\vspace*{-5mm}   
\caption{Optimal proportion $K_{1m}(s)$.}
\label{fig:1-1}
\end{center}
\end{minipage}
\begin{minipage}{.5\linewidth}
\begin{center}
\includegraphics[width=8cm,height=7cm]{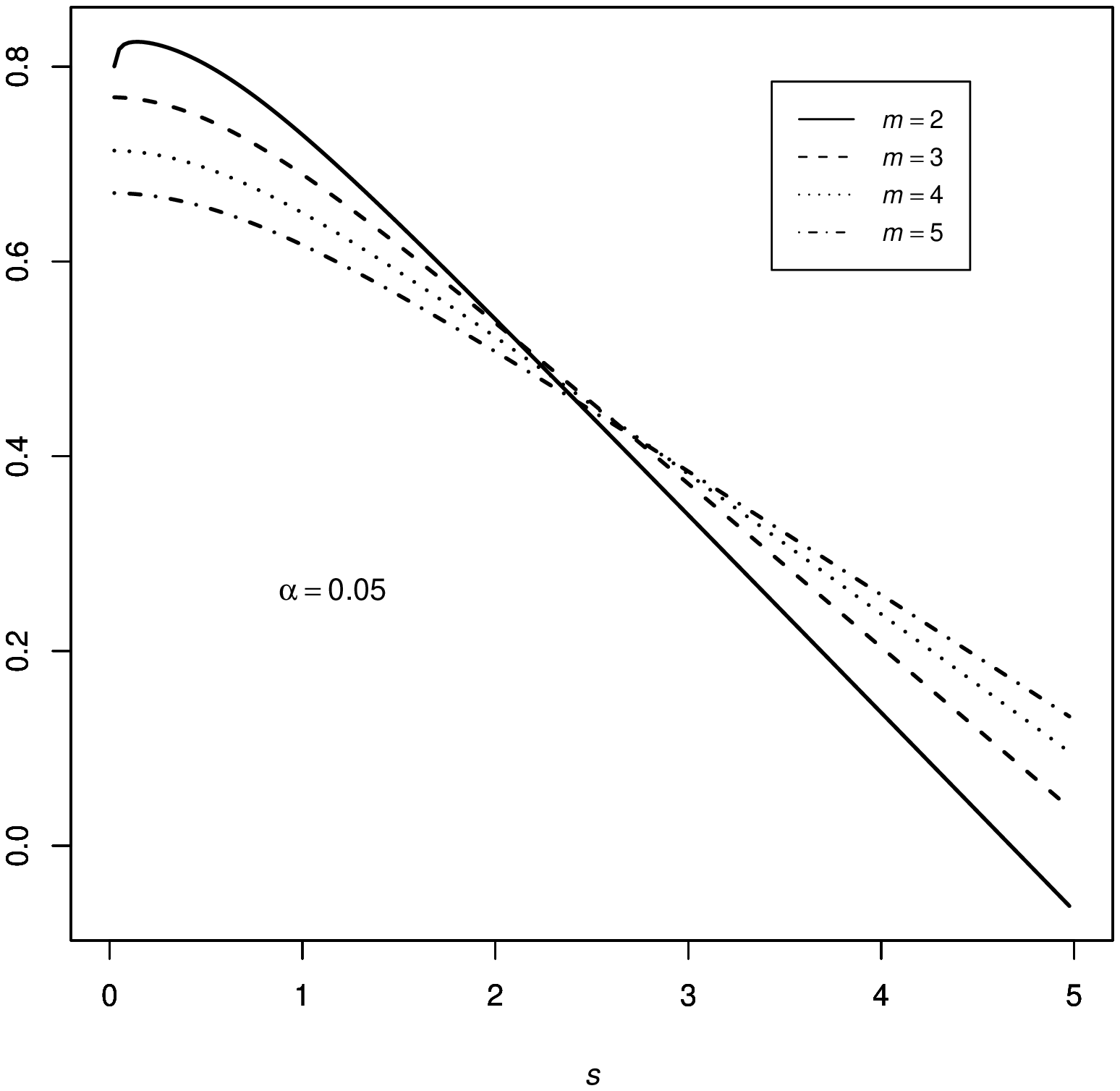}
\vspace*{-5mm}   
\caption{Optimal proportion $K_{2m}(s)$.}
\label{fig:1-1}
\end{center}
\end{minipage}
\end{figure}

The third-order power loss function $P_{\cal K}(s)$ of a $k$-test 
given by (\ref{eq:4-11}) depends on three kinds of nonnegative scalar curvatures $\kappa^2, \gamma^2$ and $H_A^{(-1)2}$. 
Among them, $\kappa^2$ measures the square of the anisotropy of the $1$-ES curvature of $M_c$ in $M_e$, and $\gamma^2$ measures the square of the mean $1$-ES curvature of $M_c$ in $M_e$. 
On the other hand, $H_A^{(-1)2}$ measures the square of the $(-1)$-ES curvature of the ancillary subspace $A(u_0)$ in $M_e$ associated with a 
$k$-test.

Figures 3, 4, 5 show the power loss coefficients 
$\xi_{0m}(s),\ \xi_{1m}(s),\ \xi_{2m}(s)$ given in (\ref{eq:4-11}), 
(\ref{eq:4-12}) as functions of $s$ in the cases of $m = 2, 3, 4, 5$ with $\alpha = 0.05$.

Figures 6, 7 show the optimal proportions $K_{1m}(s),\ K_{2m}(s)$ given by 
(\ref{eq:4-14}) as functions of $s$ in the cases of $m = 2, 3, 4, 5$ with $\alpha = 0.05$. 
\

\

\

Widely used typical tests can be regarded as the $m$-flat 
$k$-tests in the following manner (cf. Amari, 1985, Chapter 6).  
The m.l.e.\ test $(MLT)$ or the Wald test uses
\begin{align*}
g_{ab}(\hat{u}_{mle})(\hat{u}_{mle}^a-u_0^a)(\hat{u}_{mle}^b-u_0^b)\quad 
\textrm{or}\quad
g_{ab}(u_0)(\hat{u}_{mle}^a-u_0^a)(\hat{u}_{mle}^b-u_0^b)
\end{align*}
as the test statistic, where $\hat{u}_{mle}$ denotes the maximum likelihood estimator of $u$, and its characteristic is
\begin{align}
\label{eq:4-15}
k_1 = k_2 = 0,\quad H_{\kappa \lambda a}^{(-1)}(u) \equiv 0.
\end{align}
The likelihood ratio test $(LRT)$ uses 
\begin{align*}
-2\log\ \{f(\bar{x}, u_0)/f(\bar{x}, \hat{u}_{mle})\}
\end{align*}
as the test statistic, and its characteristic is
\begin{align}
\label{eq:4-16}
k_1 = k_2 = 0.5,\quad H_{\kappa \lambda a}^{(-1)}(u_0) = 0.
\end{align}
The efficient score test $(EST)$ or the Rao test uses 
\begin{align*}
g^{ab}(u_0)\partial_al(\bar{x}, u_0)\partial_bl(\bar{x}, u_0)
\end{align*}
as the test statistic, and its characteristic is
\begin{align}
\label{eq:4-17}
k_1 = k_2 = 1,\quad H_{\kappa \lambda a}^{(-1)}(u) \equiv 0.
\end{align}

\clearpage

\begin{figure}[htbp]
\begin{minipage}{.5\linewidth}
\begin{center}
\includegraphics[width=8cm,height=7cm]{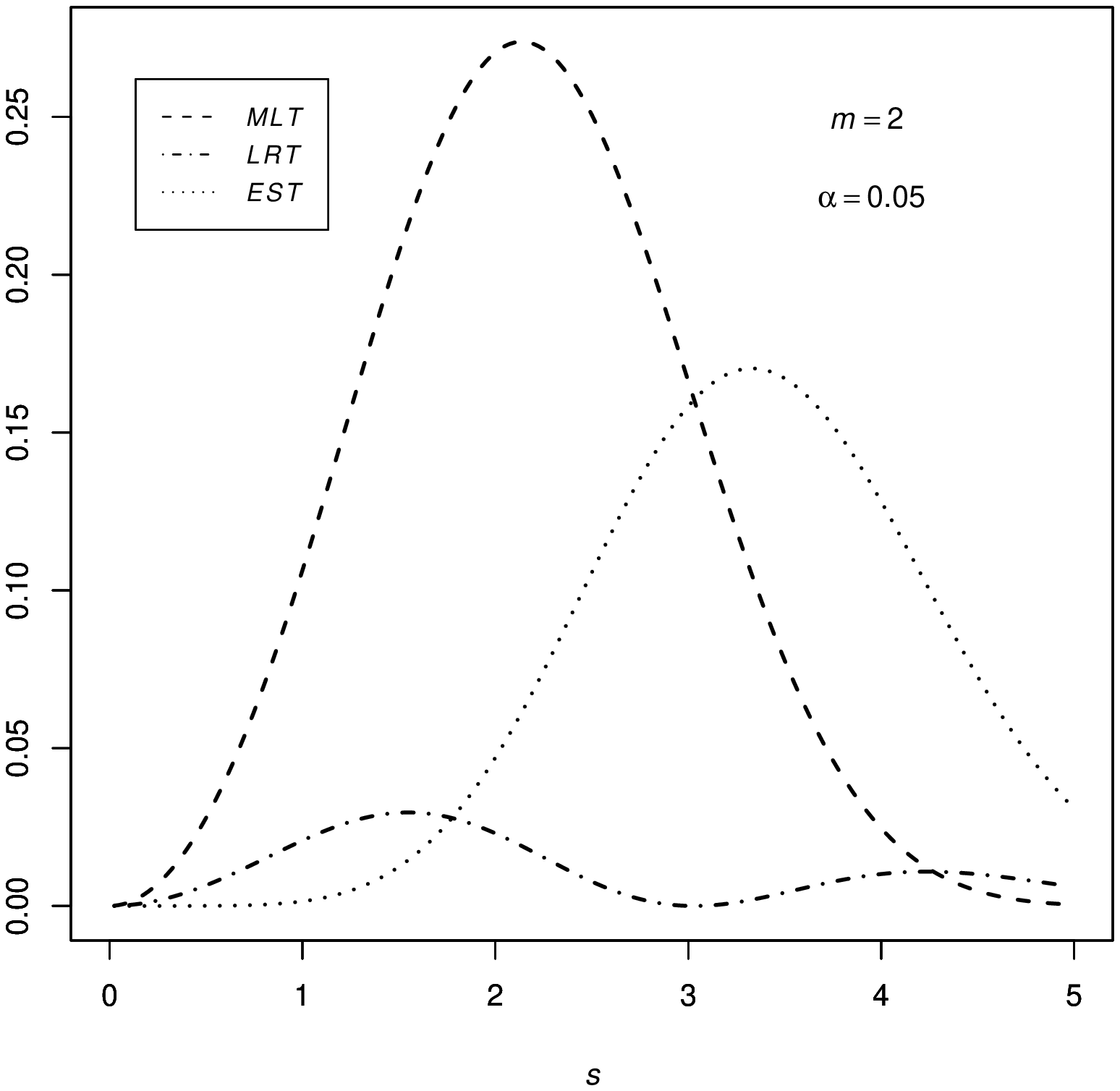}
\vspace*{-5mm}   
\caption{Power loss $\Delta P_1(s)$ of typical tests.}
\label{fig:1-1}
\end{center}
\end{minipage}
\begin{minipage}{.5\linewidth}
\begin{center}
\includegraphics[width=8cm,height=7cm]{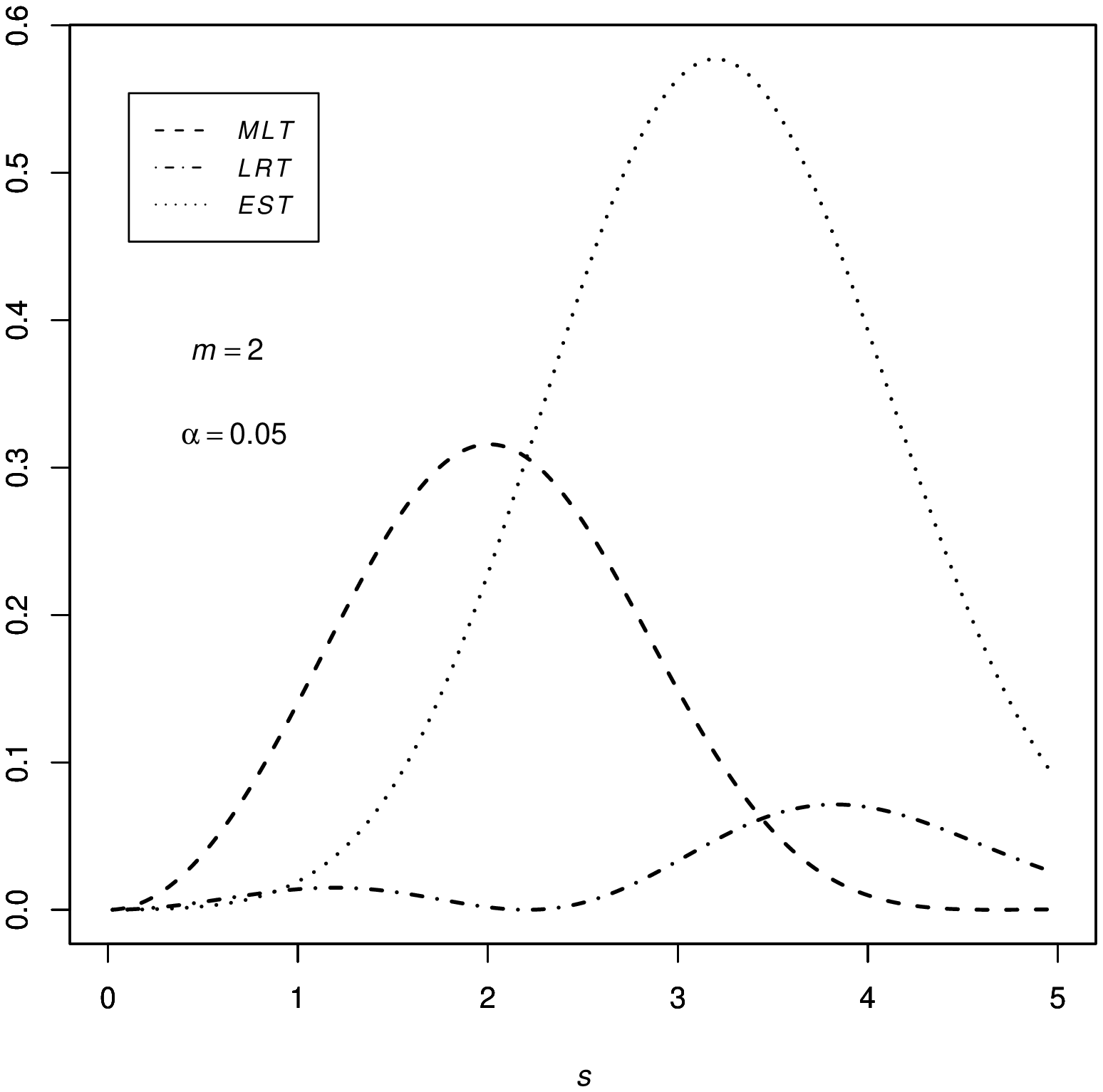}
\vspace*{-5mm}   
\caption{Power loss $\Delta P_2(s)$ of typical tests.}
\label{fig:1-1}
\end{center}
\end{minipage}
\end{figure}
\

Figures 8, 9 show the power loss functions 
$\Delta P_1(s),\ \Delta P_2(s)$ given by (\ref{eq:4-12}) 
for $MLT,\ LRT,\ EST$ in the case of $m = 2$ with $\alpha = 0. 05$.

Including these typical $k$-tests, we can design an efficient test whose ancillary family is $m$-flat and has a desired value of $Q_{(ab)k}$.  
When 
\begin{align*}
Q_{(ab)k} = k_1K_{ab\kappa}^{(1)} + k_2g_{ab}H_\kappa^{(1)}
\end{align*}
is specified, we modify the m.l.e.\ $\hat{u}_{mle}$ of $u$ to
\begin{align}
\label{eq:4-18}
&\hat{u}_{mle}'^a = \hat{u}_{mle}^a + 
g^{ab}(\hat{u}_{mle})\hat{A}_{bc}(\hat{u}_{mle}^c - u_0^c),\\
&\hat{A}_{bc} = k_1[\hat{a}_{bc} - g_{bc}(\hat{u}_{mle})\hat{a}] + 
k_2g_{bc}(\hat{u}_{mle})\hat{a},\nonumber \\
&\hat{a}_{bc} = g_{bc}(\hat{u}_{mle}) + 
\partial_b \partial_c l(\bar{x}, \hat{u}_{mle}),\quad
\hat{a} = \hat{a}_{bc}g^{bc}(\hat{u}_{mle}),\nonumber
\end{align}
and uses the quadratic form
\begin{align}
\label{eq:4-19}
g_{ab}(u_0)(\hat{u}_{mle}'^a - u_0^a)(\hat{u}_{mle}'^b - u_0^b)
\end{align}
as the test statistic. 
Then it gives the $m$-flat efficient test with the prescribed $Q_{(ab)k}$. 
In particular, by specifying
\begin{align*}
k_1 = K_{1m}(s_0),\quad k_2 = K_{2m}(s_0),
\end{align*}
we obtain the third-order $s_0$-efficient test for each $s_0$. 
Note that $\hat{a}_{bc}$ is the well-known asymptotic ancillary statistic, 
which is the expected Fisher information $g_{bc}(\hat{u}_{mle})$ 
minus the observed Fisher information 
$-\partial_b \partial_c l(\bar{x}, \hat{u}_{mle})$.

The proposed designs will be numerically examined in comparison with the sequential tests in Section 7.

\clearpage

\section{CONFORMAL TRANSFORMATION OF CURVED EXPONENTIAL FAMILY}

We reconsider the conformal transformation $M_e \mapsto {\bar M}_e$ of the 
f.r.m.\ exponential family $M_e$ by the gauge function $\nu(w) > 0$. 
As shown by (\ref{eq:3-17}) and (\ref{eq:3-22}), the metric tensor and the skewness tensor in terms of the $w$-coordinate system are changed into
\begin{align*}
\bar{g}_{\alpha \beta} = \nu g_{\alpha \beta},\quad
\bar{T}_{\alpha \beta \gamma} = \nu [T_{\alpha \beta \gamma} + 
3g_{(\alpha \beta}s_{\gamma)}].
\end{align*}
Dividing by the gauge factor $\nu$, we define
\begin{align}
\label{eq:5-1}
\tilde{g}_{\alpha \beta} = g_{\alpha \beta},\quad
\tilde{T}_{\alpha \beta \gamma} = T_{\alpha \beta \gamma} + 
3g_{(\alpha \beta}s_{\gamma)},
\end{align}
and introduce the scaled statistical manifold $\tilde{M}_e$ with the basic tensors $\tilde{g}_{\alpha \beta},\ \tilde{T}_{\alpha \beta \gamma}$.
Asymptotic properties of sequential inferences will be expressed by the geometrical quantities induced on $\tilde{M}_e$.
In view of (\ref{eq:3-23}), the $\alpha$-connection of $\tilde{M}_e$ 
in terms of the $w$-coordinate system is given by
\begin{align} 
\label{eq:5-2}
\tilde{\Gamma}_{\beta \gamma \delta}^{(\alpha)} = 
\Gamma_{\beta \gamma \delta}^{(\alpha)} + \frac{1-\alpha}{2}
(g_{\delta \beta}s_\gamma + g_{\delta \gamma}s_\beta)
- \frac{1+\alpha}{2}g_{\beta \gamma}s_\delta.
\end{align}
Then we can express the change of quantities related to the scaled statistical manifold $\tilde{M}_c$, that is, the $(-1)$-connection $\tilde{\Gamma}_{abc}^{(-1)}(u)$ of $\tilde{M}_c$, the $1$-ES curvature tensor $\tilde{H}_{ab\kappa}^{(1)}(u)$ of $\tilde{M}_c$ in $\tilde{M}_e$, and the $(-1)$-ES curvature tensor 
$\tilde{H}_{\kappa \lambda a}^{(-1)}(u)$ of $\tilde{A}(u)$ in $\tilde{M}_e$ are given by
\begin{align}
\label{eq:5-3}
\tilde{\Gamma}_{abc}^{(-1)} = \Gamma_{abc}^{(-1)} + g_{ca}s_b + g_{cb}s_a,\quad
\tilde{H}_{ab\kappa}^{(1)} = H_{ab\kappa}^{(1)} - g_{ab}s_\kappa,\quad
\tilde{H}_{\kappa \lambda a}^{(-1)} = H_{\kappa \lambda a}^{(-1)},
\end{align}
so that the mean 1-ES curvature tensor $\tilde{H}_\kappa^{(1)}(u)$ of $\tilde{M}_c$ in $\tilde{M}_e$, and the 1-ES umbilic tensor 
$\tilde{K}_{ab\kappa}^{(1)}(u)$ of $\tilde{M}_c$ in $\tilde{M}_e$ are
\begin{align}
\label{eq:5-4}
\tilde{H}_\kappa^{(1)} = H_\kappa^{(1)} - s_\kappa,\quad 
\tilde{K}_{ab\kappa}^{(1)} = K_{ab\kappa}^{(1)}.
\end{align}

We remark that the scaled 1-ES curvature tensor variable 
$\tilde{h}_{ab}^{(1)}(x^n, n, u)$ of 
$\tilde{M}_c$ in $\tilde{R}_{\tilde{M}_c}$ can be also introduced and expressed as
\begin{align}
\label{eq:5-5}
&\tilde{h}_{ab}^{(1)}(x^n, n, u) = \frac{1}{\sqrt{\nu}} 
\bar{h}_{ab}^{(1)}(x^n, n, u)
= \frac{1}{\sqrt{\nu}} 
\left[\partial_a \partial_b \bar{l}_n(x^n, u) + \bar{g}_{ab}(u) - 
\bar{\Gamma}_{ab}^{(1)c}(u)\partial_c \bar{l}_n(x^n, u)\right]\nonumber \\
&\qquad \qquad \quad \ \ = \frac{1}{\sqrt{\nu}} 
\left[\tilde{H}_{ab\kappa}^{(1)}(u)g^{\kappa \lambda}(u)
\partial_\lambda \bar{l}_n(x^n, u) - g_{ab}(u)n^{(1)}(x^n, w)\right],\\
\label{eq:5-6}
&n^{(1)}(x^n, w) = n - \nu(w) - s^\alpha \partial_\alpha \bar{l}_n(x^n, u),\quad s^\alpha = s_\beta g^{\beta \alpha},\quad 
\bar{l}_n(x^n, u) = \log \bar{f}_c(x^n, n, u). 
\end{align}

This expression shows that due to the additional term 
$g_{ab}(u)n^{(1)}(x^n, w)$, the $1$-ES curvature tensor $\tilde{H}_{ab\kappa}^{(1)}(u)$ of $\tilde{M}_c$ in $\tilde{M}_e$ and the scaled 1-ES curvature tensor variable 
$\tilde{h}_{ab}^{(1)}(x^n, n, u)$ of $\tilde{M}_c$ in $\tilde{R}_{\tilde{M}_c}$ are not in general the equivalent notions.

From Corollary 3.1 which states the condition for preserving the exponentiality of $M_e$ under the appropriate conformal transformation, we obtain the following result as to the above two curvature notions.

\begin{theorem}
For a scaled curved exponential family $\tilde{M}_c$, 
the squared scalar curvature $\tilde{\lambda}_\infty^2$ of the scaled 1-ES curvature tensor variable $\tilde{h}_{ab}^{(1)}(x^n, n, u)$ of 
$\tilde{M}_c$ in $\tilde{R}_{\tilde{M}_c}$ is not less than 
the squared scalar curvature $\tilde{\lambda}_e^2$ of the $1$-ES curvature tensor $\tilde{H}_{ab\kappa}^{(1)}(u)$ of 
$\tilde{M}_c$ in $\tilde{M}_e$ 
\begin{align}
\label{eq:5-7}
\tilde{\lambda}_\infty^2 = \langle \tilde{h}_{ab}^{(1)}(x^n, n, u)
\tilde{h}_{cd}^{(1)}(x^n, n, u)\rangle_u g^{ac}g^{bd} \ge 
\tilde{\lambda}_e^2 = 
\tilde{H}_{ab\kappa}^{(1)}\tilde{H}_{cd\kappa}^{(1)}
g^{ac}g^{bd}g^{\kappa \lambda}, 
\end{align}
where the equality holds if and only if the condition (\ref{eq:3-43}) in  Corollary 3.1 is satisfied.

For a stopping time $\tau_c\in {\cal C}$ introduced in Section 2, the equality asymptotically holds in the sense
\begin{align}
\label{eq:5-8}
\tilde{\lambda}_\infty^2 = \tilde{\lambda}_e^2 + o(1)\quad 
\textrm{as}\quad \nu_c(w)\ \to\ \infty.
\end{align}
\end{theorem}
\

\begin{proof}
From (\ref{eq:5-5}) since 
$\partial_\lambda \bar{l}_n(x^n, u)\in S(\bar{{\cal T}}_{\bar{M}_e}),\ 
n^{(1)}(x^n, w)\in S(\bar{{\cal N}}_{\bar{M}_e})$, we have
\begin{align*}
\tilde{\lambda}_\infty^2 = \tilde{\lambda}_e^2 + 
\frac{m}{\nu(w)} \langle n^{(1)}(x^n, w)^2\rangle_u \ge \tilde{\lambda}_e^2,
\end{align*}
and it follows that
\begin{align*}
\tilde{\lambda}_\infty^2 = \tilde{\lambda}_e^2\ \Leftrightarrow\ 
n^{(1)}(x^n, w) \equiv 0.
\end{align*}

For a stopping time $\tau_c\in {\cal C}$ introduced in Section 2, 
since (\ref{eq:3-44}) holds, we have
\begin{align*}
\tilde{\lambda}_\infty^2 = \tilde{\lambda}_e^2 + 
\frac{m}{\nu_c(w)} \langle n_c^{(1)}(x^n, w)^2\rangle_u = 
\tilde{\lambda}_e^2 + o(1)\quad 
\textrm{as}\quad \nu_c(w)\ \to\ \infty.
\end{align*}

This completes the proof of the theorem.
\end{proof}
\

From (\ref{eq:5-3}), (\ref{eq:5-4}) we know that under any conformal transformation 
$M_e \mapsto \tilde{M}_e$,
\begin{align}
\label{eq:5-9}
\tilde{\kappa}^2 = \kappa^2 = 
K_{ab\kappa}^{(1)}K_{cd\kappa}^{(1)}g^{ac}g^{bd}g^{\kappa \lambda},\quad 
\tilde{H}_A^{(-1)2} = H_A^{(-1)2} = 
H_{\kappa \lambda a}^{(-1)}H_{\mu \nu b}^{(-1)}g^{\kappa \mu}
g^{\lambda \nu}g^{ab},
\end{align}
and by setting $s_\kappa(u) = H_\kappa^{(1)}(u)$, we have
\begin{align}
\label{eq:5-10}
\tilde{\gamma}^2 = m\tilde{H}_\kappa^{(1)}
\tilde{H}_\lambda^{(1)}g^{\kappa \lambda}= 0.
\end{align}
When the ancillary family $A$ is $m$-flat and $M_c$ is totally 
$e$-umbilic in $M_e$, we further get 
\begin{align}
\label{eq:5-11}
\tilde{\kappa}^2 = \kappa^2 = 0,\quad
\tilde{H}_A^{(-1)2} = H_A^{(-1)2} = 0.
\end{align}
As noted after (\ref{eq:4-18}) $M_c$ is also totally $e$-umbilic in $R_{M_c}$, and then from (\ref{eq:3-46}) $M_c$ is conformally $m(e)$-flat, so that there exists a gauge function $\nu(u) > 0$ and a $(-1)$-affine coordinate system $\tilde{u}_\nu = (\tilde{u}_\nu^{\tilde{a}}),\ 
\tilde{a} = 1, \dots, m$ of $\tilde{M}_c$ satisfying 
\begin{align}
\label{eq:5-12}
\tilde{\Gamma}_{\tilde{a}\tilde{b}\tilde{c}}^{(-1)}(\tilde{u}_\nu) \equiv 0.
\end{align}
These possibilities will be realized in the subsequent sections.

\section{SEQUENTIAL TEST IN CURVED EXPONENTIAL FAMILY}

We deal with an unbiased sequential test of the simple null hypothesis 
\begin{align*}
H_0 : u = u_0\quad \textrm{vs}\quad H_1 : u \neq u_0,
\end{align*}
in an $(n, m)$-curved exponential family $M_c$. 
Let $K$ be a large number playing the role of the average time of observations, and let $\nu(\eta) > 0\ (\nu(w) > 0)$ be a smooth gauge function defined on $M_e$ in the $\eta$-($w$-)coordinate system.

The sample mean up to time $n$, $\bar{x}_n = \sum_{t=1}^n x(t)/n$ has the same value in the $\eta$-coordinate system, and its value 
in the $w$-coordinate system is denoted by 
${\hat w}_n = ({\hat u}_n, {\hat v}_n) = \eta^{-1}(\bar{x}_n)$. 
The random stopping time $\tau$ is assumed to satisfy 
(see Okamoto et al., 1991) 
\begin{align}
\label{eq:6-1}
&\tau = K\nu({\hat w}_\tau) + c({\hat u}_\tau) + \varepsilon,\\ 
&c(u) = -\frac{1}{2}(\partial_\alpha s_\beta - 
\Gamma_{\alpha \beta}^{(-1)\gamma}s_\gamma - 
s_\alpha s_\beta)g^{\alpha \beta},\quad 
s_\alpha = \partial_\alpha \log \nu(w),\nonumber \\ 
&\varepsilon = O_p(1),\quad E_u[\varepsilon] = o(1),\quad 
E_u[\tau] = K\nu(u),\quad \textrm{Var}_u[\tau] = O(K). \nonumber
\end{align}
The term $c(u)$ is due to the bias of ${\hat w}_\tau$ from the true 
$w = (u, 0)$, which is obtained by the requirement 
$E_u[\tau] = K\nu(u)$. The term $\varepsilon$ includes a rounding error 
and the ``overshooting'' at the stopping time $\tau$.

We note that a stopping time $\tau_c\in {\cal C}$ introduced in Section 2 satisfies the above conditions with $K = c^{1/(1-\beta)}$. 
On the other hand, by expanding (\ref{eq:6-1}) at the true $w = (u, 0)$ we have
\begin{align*}
\tau &= K\nu(w) + K\nu(w)s_\alpha (\hat{w}^\alpha - w^\alpha) + O_p(1)\\ 
&= K\nu(w) + s^\beta \partial_\beta \bar{l}_\tau(X^\tau, w) + O_p(1),\quad
s^\beta = s_\alpha g^{\alpha \beta},
\end{align*}
which implies that the condition (\ref{eq:3-43}) in Corollary 3.1 is asymptotically satisfied in the sense
\begin{align*}
&\langle n^{(1)}(x^n, w)^2\rangle_w = O(1)\quad \textrm{as}\quad 
K\ \to\ \infty,\\
&n^{(1)}(x^n, w) = n - K\nu(w) - s^\beta \partial_\beta \bar{l}_n(x^n, w).
\end{align*}

The formulation proceeds almost the same way as in the nonsequential case. 
Let $\tilde{R}$ be the critical region of a sequential test 
$\tilde{\cal T}$, and let $U_K(s)$ be a set of alternatives given by
\begin{align*}
&U_K(s) = \{u_{s,e} = u_0 + se/\sqrt{K\nu_0}\in \tilde{M}_c,\ 
s \ge 0,\ g_{ab}(u_0)e^ae^b = 1 \},\quad
\nu_0 = \nu(u_0).
\end{align*}
We compare sequential tests based on the average power 
$\tilde{P}_{\tilde{\cal T}}(s)$ defined by
\begin{align*}
\tilde{P}_{\tilde{\cal T}}(s) = \langle \tilde{P}_{\tilde{\cal T}}(u)\rangle_e = 
\int_{u\in U_K(s)} \tilde{P}_{\tilde{\cal T}}(u)du/S_K(s),
\end{align*}
where $S_K(s)$ denotes the area of $U_K(s)$.

By expanding $\tilde{P}_{\tilde{\cal T}}(s)$ in the power series of $(K\nu_0)^{-1/2}$, we have
\begin{align*}
\tilde{P}_{\tilde{\cal T}}(s) = \tilde{P}_{\tilde{\cal T}1}(s) + 
\tilde{P}_{\tilde{\cal T}2}(s)(K\nu_0)^{-1/2} + 
\tilde{P}_{\tilde{\cal T}3}(s)(K\nu_0)^{-1} + O(K^{-3/2}),
\end{align*}
and call $\tilde{P}_{\tilde{\cal T}i}(s)$ the $i$-th order power of a sequential test $\tilde{\cal T}$ at $s$. When discussing the $i$-th order power, the tests are assumed to satisfy the level condition 
\begin{align*}
\tilde{P}_{\tilde{\cal T}}(0) = \alpha + O(K^{-i/2}),
\end{align*}
up to the same $i$-th order.

A sequential test $\tilde{\cal T}$ is said to be (first-order uniformly) efficient when
\begin{align*}
\tilde{P}_{\tilde{\cal T}1}(s) \ge \tilde{P}_{\tilde{\cal T}'1}(s),\quad 
\forall s \ge 0,\ \forall \tilde{\cal T}'.
\end{align*}
It can be shown that an efficient sequential test is automatically second-order uniformly efficient, and that there does not in general exist a third-order uniformly efficient sequential test.

Introducing the envelope power function by 
$\tilde{P}^*(s) = \sup_{\tilde{\cal T}} \tilde{P}_{\tilde{\cal T}}(s)$,  
and expanding it as 
\begin{align*}
\tilde{P}^*(s) = \tilde{P}_1^*(s) + \tilde{P}_2^*(s)(K\nu_0)^{-1/2} + \tilde{P}_3^*(s)(K\nu_0)^{-1} + O(K^{-3/2}),
\end{align*}
we call $\tilde{P}_i^*(s)$ the $i$-th order envelope power at $s$. 
Then the (third-order) power-loss function 
$\Delta \tilde{P}_{\tilde{\cal T}}(s)$ of an efficient sequential test $\tilde{\cal T}$ is defined by
\begin{align}
\label{eq:6-2}
\Delta \tilde{P}_{\tilde{\cal T}}(s) = 
\tilde{P}_3^*(s) - \tilde{P}_{\tilde{\cal T}3}(s).
\end{align}

From Propositions 4.1, 4.2, Theorem 5.1 and Appendices 1, 2, 3, we can state the asymptotic results for the possibility of uniformly efficient sequential test.

\begin{theorem}
For an unbiased sequential test $\tilde{\cal T}$ of $D = \{u_0\}$ in a curved exponential family $M_c$, the following hold.\\
(i)\ A sequential test $\tilde{\cal T}$ is first-order uniformly efficient, if and only if the associated ancillary family is asymptotically orthogonal in the sense
\begin{align*}
g_{a\kappa}(u) = O(K^{-1/2}),\quad \forall u\in 
\partial \tilde{R}_{\tilde{M}_c},
\end{align*}
and the optimal intersection 
$\tilde{R}_{\tilde{M}_c}^* = \tilde{R}^*\cap \tilde{M}_c$ is given by
\begin{align*}
\tilde{R}_{\tilde{M}_c}^* = \{{\tilde u}_{n,0}\ |\ 
g_{ab}(u_0){\tilde u}_{n,0}^a{\tilde u}_{n,0}^b \ge c_0^2 \},\quad  
{\tilde u}_{n,0} = \sqrt{K\nu_0}({\hat u}_n - u_0),\quad 
c_0^2 = \chi_{m,\alpha}^2.
\end{align*}
The first-order envelope power function $\tilde{P}_1^*(s)$ is the same as 
$P_1^*(s)$ given in Proposition 4.1. 
\

\

\noindent
(ii)\ The second- and third-order powers are analyzed by assuming the intersection $\tilde{R}_{\tilde{M}_c}$ as
\begin{align}
\label{eq:6-3}
&\tilde{R}_{\tilde{M}_c} = \{\tilde{u}_{n,0}\ |\ 
g_{ab}(u_0)\tilde{u}_{n,0}^{*a}\tilde{u}_{n,0}^{*b} \ge 
(c_0 + \tilde{\epsilon})^2\},\\
&\tilde{u}_{n,0}^{*a} = \tilde{u}_{n,0}^a + 
\tilde{\Gamma}_{\alpha \beta}^{(-1)a}g^{\alpha \beta}(\hat{u}_n)/
(2\sqrt{K\nu_0}),\quad 
\tilde{\Gamma}_{\alpha \beta}^{(-1)a} = 
\Gamma_{\alpha \beta}^{(-1)a} + 
\delta_\alpha^a s_\beta + \delta_\beta^a s_\alpha. \nonumber
\end{align}
From Appendices 2, 3 we have $\tilde{\epsilon} = O(K^{-1})$, which depends on the ancillary family through $\tilde{Q}_{ab\kappa}$ and 
$H_{\kappa \lambda a}^{(-1)}$. 
Hence it follows that the second-order power function 
$\tilde{P}_{\tilde{{\cal T}}2}(s)$ is common to all the first-order uniformly efficient sequential tests.
\

\

\noindent
(iii)\ A class of efficient sequential tests is introduced whose ancillary family satisfies
\begin{align}
\label{eq:6-4}
\tilde{Q}_{(ab)\kappa} = k_1 K_{ab\kappa}^{(1)} + k_2 g_{ab}\tilde{H}_\kappa^{(1)},\quad k_1, k_2 \in {\mathbb R},
\end{align}
and a sequential test $\tilde{{\cal T}}$ with a pair of proportions 
$(k_1, k_2)$ is said to be the sequential $k = (k_1, k_2)$-test. 
The third-order power loss function $\Delta \tilde{P}_{\cal K}(s)$ of an efficient sequential $k = (k_1, k_2)$-test $\tilde{\cal T}$ is given by
\begin{align}
\label{eq:6-5}
\Delta \tilde{P}_{\cal K}(s) = \kappa^2 \Delta P_1(s) + 
\tilde{\gamma}^2 \Delta P_2(s) + H_A^{(-1)2}\xi_{0m}(s),\quad
\tilde{\gamma}^2 = m\tilde{H}_\kappa^{(1)}
\tilde{H}_\lambda^{(1)}g^{\kappa \lambda},
\end{align}
where functions $\Delta P_1(s),\ \Delta P_2(s),\ \xi_{0m}(s)$ are the same as given in Proposition 4.2.
\

\

\noindent
(iv)\ Using the stopping rule derived from (\ref{eq:3-41}) 
(cf. Okamoto et al., 1991),
\begin{align}
\label{eq:6-6}
\tau = \inf\left\{n\ \Big|\ -\frac{1}{m}\sum_{t=1}^n \partial_a\partial_b
{\bar l}_1(x(t), {\hat u}_{mle,n})g^{ab}({\hat u}_{mle,n}) \ge
K\nu({\hat u}_{mle,n}) + c(\hat{u}_{mle,n}) \right\},
\end{align}
we have $\tilde{\gamma}^2 = 0$, so that the third-order power loss function of a sequential $m$-flat $k = (k_1, k_2)$-test with 
$H_{\kappa \lambda a}^{(-1)}(u_0) = 0$ is given by
\begin{align}
\label{eq:6-7}
\Delta \tilde{P}_{\cal K}(s) = \kappa^2 \Delta P_1(s).
\end{align}
When $M_c$ is totally $e$-umbilic in $M_e$ or equivalently in $R_{M_c}$, 
we have $\kappa^2 = 0$, 
and thus a sequential $m$-flat $k = (k_1, k_2)$-test is third-order uniformly efficient in the sense 
\begin{align}
\label{eq:6-8}
\Delta \tilde{P}_{\cal K}(s) = 0.
\end{align}
\end{theorem}
\

\clearpage

\section{EXAMPLES}

\subsection{von Mises-Fisher model}

This is an $(m+1, m)$-curved exponential family, of which density function  with respect to the invariant measure on the $m$-dimensional unit sphere under rotational transformations is given by 
(cf. Barndorff-Nielsen et al., 1989, p.76)
\begin{align*}
&f_c(x, u) = \exp\{\theta(u)\cdot x - \psi[\theta(u)]\},\quad 
\theta\cdot x = \theta^1x_1 + \theta^2x_2 + \cdots + \theta^{m+1}x_{m+1},\\ 
&\ \theta = r\xi = (r\xi^i),\quad 
\xi \in S^m = \{\xi \in {\mathbb R}^{m+1}\ |\ 
\xi\cdot \xi = 1\},\quad 
x = (x_i) \in S^m,\\
&\ \psi(\theta) = -\log a_m(r),\quad
1/a_m(r) = (2\pi)^{(m+1)/2}r^{(1-m)/2}I_{(m-1)/2}(r),\quad r > 0,
\end{align*}
where $I_{(m-1)/2}(r)$ is the modified Bessel function of the first kind and of order $(m-1)/2$. 
The concentration parameter $r$ is assumed to be a given positive constant. The parametric representations $\theta = \theta(u)$ and 
$\eta = \eta(u)$ are given by\\

\begin{tabular}{ll}
$\theta^1(u) = r\cos u^1$ & $\eta_1(u) = r^\dagger \cos u^1$ \\
$\theta^2(u) = r\sin u^1 \cos u^2$ & 
$\eta_2(u) = r^\dagger \sin u^1 \cos u^2$ \\
$\theta^3(u) = r\sin u^1 \sin u^2 \cos u^3$ & 
$\eta_3(u) = r^\dagger \sin u^1 \sin u^2 \cos u^3$ \\
$\qquad \qquad \cdots$ & $\qquad \qquad \cdots$ \\
$\theta^{m+1}(u) = r\sin u^1 \sin u^2 \cdots \sin u^{m-1} \sin u^m$ & 
$\eta_{m+1}(u) = r^\dagger \sin u^1 \sin u^2 \cdots \sin u^{m-1} 
\sin u^m$  
\end{tabular}
\

\noindent
where $0 \le u^1, \dots, u^{m-1} \le \pi,\ 0 \le u^m < 2\pi$\ and\
$r^\dagger = -d\log a_m(r)/dr = I_{(m+1)/2}(r)/I_{(m-1)/2}(r)$. 
Note that $E[x] = r^\dagger \xi$, and $r^\dagger$ is a strictly increasing function of $r$ that maps $(0, \infty)$ onto $(0, 1)$.

From these representations the tangent vectors $B_a^i(u)$ and $B_{ai}(u)$ can be calculated, and then the unit normal vectors 
$B_\kappa^i(u)$ and $B_{\kappa i}(u)\ (\kappa = m+1)$ are derived from the relations 
$B_\kappa^i(u)B_{ai}(u) = 0$ and $B_{\kappa i}(u)B_a^i(u) = 0$ 
as follows. \\

\begin{tabular}{ll}
$B_\kappa^1(u) = \cos u^1$ & $B_{\kappa 1}(u) = \cos u^1$ \\
$B_\kappa^2(u) = \sin u^1 \cos u^2$ & 
$B_{\kappa 2}(u) = \sin u^1 \cos u^2$ \\
$B_\kappa^3(u) = \sin u^1 \sin u^2 \cos u^3$ & 
$B_{\kappa 3}(u) = \sin u^1 \sin u^2 \cos u^3$ \\
$\qquad \qquad \cdots$ & $\qquad \qquad \cdots$ \\
$B_\kappa^{m+1}(u) = \sin u^1 \sin u^2 \cdots \sin u^{m-1} \sin u^m$ & 
$B_{\kappa\ {m+1}}(u) = \sin u^1 \sin u^2 \cdots \sin u^{m-1} 
\sin u^m$  
\end{tabular}
\

\noindent
The related geometrical quantities are given below. 
\begin{align*}
&g_{ab}(u) = \delta_{ab}rr^\dagger \prod_{c=1}^a\sin^2 u^{c-1},\quad 
\sin^2 u^0 = 1,\\
&H_{ab\kappa}^{(1)}(u) = -\frac{1}{r^\dagger}g_{ab}(u),\quad
H_\kappa^{(1)}(u) = -\frac{1}{r^\dagger},\quad
K_{ab\kappa}^{(1)}(u) \equiv 0.
\end{align*}

This model is totally $e$-umbilic in $M_e$ or equivalently in $R_{M_c}$ with 
$\kappa^2 = 0,\ \gamma^2 = m/(r^\dagger)^2$,  
and since this $M_c$ is conformally $m$($e$)-flat,  
there exist a gauge function $\nu(u) > 0$ and a $(-1)$-affine coordinate system $\tilde{u}_\nu = (\tilde{u}_\nu^{\tilde{a}})$ of $\tilde{M}_c$ 
satisfying   
$\tilde{\Gamma}_{\tilde{a}\tilde{b}\tilde{c}}^{(-1)}(\tilde{u}_\nu) 
\equiv 0$. 
The concrete forms are (see Kumon et al., 2011)
\begin{align}
\label{eq:7-1}
\tilde{u}_\nu^{\tilde{a}} = \nu(u)D^{\tilde{a}i}\eta_i(u),\quad
\nu(u) = \frac{1}{\prod_{a=1}^m |\sin u^a|},
\end{align}
where $D^{\tilde{a}i}\ (\tilde{a} = 1, \dots, m,\ i = 1, \dots, m+1)$ are constants with $\textrm{rank}\ D^{\tilde{a}i} = m$.

\subsection{Hyperboloid model}

This is an $(m+1, m)$-curved exponential family, of which density function  with respect to the invariant measure on the $m$-dimensional unit hyperboloid under hyperbolic transformations is given by 
(cf. Barndorff-Nielsen et al., 1989, p.104)
\begin{align*}
&f_c(x, u) = \exp\{\theta(u)\cdot x - \psi[\theta(u)]\},\quad 
\theta^1 = -r\xi^1,\quad \theta^i = r\xi^i,\quad i = 2, \dots, m+1,\\ 
&\ \xi = (\xi^i)\in H^m = \{\xi \in {\mathbb R}^{m+1}\ |\ 
\xi*\xi = 1,\ \xi^1 > 0 \},\quad 
x = (x_i) \in H^m,\\
&\ \xi*\xi = (\xi^1)^2 - (\xi^2)^2 - \cdots - (\xi^{m+1})^2,\\
&\ \psi(\theta) = -\log a_m(r),\quad
1/a_m(r) = 2(2\pi)^{(m-1)/2}r^{(1-m)/2}K_{(m-1)/2}(r),\quad r > 0,
\end{align*}
where $K_{(m-1)/2}(r)$ is the modified Bessel function of the third kind and of order $(m-1)/2$. 
The concentration parameter $r$ is assumed to be a given positive constant. The parametric representations $\theta = \theta(u)$ and $\eta = \eta(u)$ are given by\\

\begin{tabular}{ll}
$\theta^1(u) = -r\cosh u^1$ & $\eta_1(u) = r^\dagger \cosh u^1$ \\
$\theta^2(u) = r\sinh u^1 \cos u^2$ & 
$\eta_2(u) = r^\dagger \sinh u^1 \cos u^2$ \\
$\theta^3(u) = r\sinh u^1 \sin u^2 \cos u^3$ & 
$\eta_3(u) = r^\dagger \sinh u^1 \sin u^2 \cos u^3$ \\
$\qquad \qquad \cdots$ & $\qquad \qquad \cdots$ \\
$\theta^{m+1}(u) = r\sinh u^1 \sin u^2 \cdots \sin u^{m-1} \sin u^m$ & 
$\eta_{m+1}(u) = r^\dagger \sinh u^1 \sin u^2 \cdots \sin u^{m-1} 
\sin u^m$  
\end{tabular}
\

\noindent
where $u^1 \in {\mathbb R},\ 0 \le u^2, \dots, u^{m-1} \le \pi,\ 
0 \le u^m < 2\pi$ and
$r^\dagger = d\log a_m(r)/dr = K_{(m+1)/2}(r)/K_{(m-1)/2}(r)$. 
Note that $E[x] = r^\dagger \xi$, and $r^\dagger$ is a strictly decreasing function of $r$ that maps $(0, \infty)$ onto $(1, \infty)$.

From these representations the tangent vectors $B_a^i(u)$ and $B_{ai}(u)$ can be calculated, and then the unit normal vectors 
$B_\kappa^i(u)$ and $B_{\kappa i}(u)\ (\kappa = m+1)$ are derived from the relations 
$B_\kappa^i(u)B_{ai}(u) = 0$ and $B_{\kappa i}(u)B_a^i(u) = 0$ 
as follows. \\

\begin{tabular}{ll}
$B_\kappa^1(u) = \cosh u^1$ & $B_{\kappa 1}(u) = \cosh u^1$ \\
$B_\kappa^2(u) = -\sinh u^1 \cos u^2$ & 
$B_{\kappa 2}(u) = \sinh u^1 \cos u^2$ \\
$B_\kappa^3(u) = -\sinh u^1 \sin u^2 \cos u^3$ & 
$B_{\kappa 3}(u) = \sinh u^1 \sin u^2 \cos u^3$ \\
$\qquad \qquad \cdots$ & $\qquad \qquad \cdots$ \\
$B_\kappa^{m+1}(u) = -\sinh u^1 \sin u^2 \cdots \sin u^{m-1} \sin u^m$ & 
$B_{\kappa\ {m+1}}(u) = \sinh u^1 \sin u^2 \cdots \sin u^{m-1} 
\sin u^m$  
\end{tabular}
\

\noindent
The related geometrical quantities are given below. 
\begin{align*}
&g_{11}(u) = rr^\dagger,\quad 
g_{ab}(u) = \delta_{ab}rr^\dagger \sinh^2 u^1 \prod_{c=2}^a\sin^2 u^{c-1},
\ \ \sin^2 u^1 = 1,\ \ a = 2, \dots, m,\\
&H_{ab\kappa}^{(1)}(u) = -\frac{1}{r^\dagger}g_{ab}(u),\quad
H_\kappa^{(1)}(u) = -\frac{1}{r^\dagger},\quad
K_{ab\kappa}^{(1)}(u) \equiv 0.
\end{align*}
This model is totally $e$-umbilic in $M_e$ or equivalently in $R_{M_c}$ with 
$\kappa^2 = 0,\ \gamma^2 = m/(r^\dagger)^2$, and since this $M_c$ is conformally $m$($e$)-flat,  
there exist a gauge function $\nu(u) > 0$ and a $(-1)$-affine coordinate system $\tilde{u}_\nu = (\tilde{u}_\nu^{\tilde{a}})$ of $\tilde{M}_c$ 
satisfying   
$\tilde{\Gamma}_{\tilde{a}\tilde{b}\tilde{c}}^{(-1)}(\tilde{u}_\nu) 
\equiv 0$. 
The concrete forms are (see Kumon et al., 2011)
\begin{align}
\label{eq:7-2}
\tilde{u}_\nu^{\tilde{a}} = \nu(u)D^{\tilde{a}i}\eta_i(u),\quad
\nu(u) = \frac{1}{|\sinh u^1| \prod_{a=2}^m |\sin u^a|},
\end{align}
where $D^{\tilde{a}i}\ (\tilde{a} = 1, \dots, m,\ i = 1, \dots, m+1)$ are constants with $\textrm{rank}\ D^{\tilde{a}i} = m$.

\subsection{Numerical results}

We examine our theoretical results numerically by using the von Mises-Fisher and the hyperboloid models.
It will be shown that the sequential $m$-flat test is uniformly superior to the nonsequential $m$-flat tests. 
At first the third-order power loss functions in the nonsequential tests 
are evaluated in the two dimensional case $m=2$. 
We take 20 values of the distance 
$s$ from 0 to 5, and for each value of $s$, we select the number of 
$H_1$ equally spaced alternatives. 
At each point of alternatives, the number of $N$ random simulated data are generated. 
The empirical average power $P_{{\cal K}e}(s)$ of an $m$-flat $k$-test is  calculated as the number of rejections relative to the number of $H_1$. 
Tests are compared by the empirical third-order power loss function
\begin{align*}
\Delta P_{{\cal K}e}(s) = N[P_{{\cal K}^*e}(s) - P_{{\cal K}e}(s)]/\gamma^2,
\end{align*}
where $P_{{\cal K}^*e}(s)$ denotes the empirical average power of the third-order $s$-efficient test for each value of $s$.

As for the von Mises-Fisher model, numerical results are based on the following set of values
\begin{align*}
&m = 2,\quad \alpha = 0.05,\quad r = 0.1,\quad
H_0\ :\ (u_0^1, u_0^2) = (\pi/2, \pi/2),\quad
N = 2000,\quad H_1 = 1000, 
\end{align*}
and for the hyperboloid model, numerical results are based on the following set of values
\begin{align*}
&m = 2,\quad \alpha = 0.05,\quad r = 2,\quad
H_0\ :\ (u_0^1, u_0^2) = (1, \pi/2),\quad
N = 50,\quad H_1 = 5000.
\end{align*}

Figure 10 shows the von Mises-Fisher model, and Figure 11 shows the hyperboloid model. 
In each figure, $\Delta P_{{\cal K}e}(s)$ are depicted for 
$MLT,\ LRT,\ EST$. 
The notations in the figures indicate the following tests.
\begin{align*}
DMLT\ :\ k = (0, 0),\quad
DLRT\ :\ k = (0.5, 0.5),\quad
DEST\ :\ k = (1, 1).
\end{align*}
We remark that these tests are designed by (\ref{eq:4-19}), which are the approximate versions of the exact $m$-flat $k=(k_1, k_2)$-tests.
Thus the theoretical third-order power loss functions 
$\Delta P_2(s) = \Delta P_{{\cal K}}(s)/\gamma^2$ shown by Figure 9 are not precisely reflected in these figures.

\begin{figure}[htbp]
\begin{minipage}{.5\linewidth}
\begin{center}
\includegraphics[width=8cm,height=7cm]{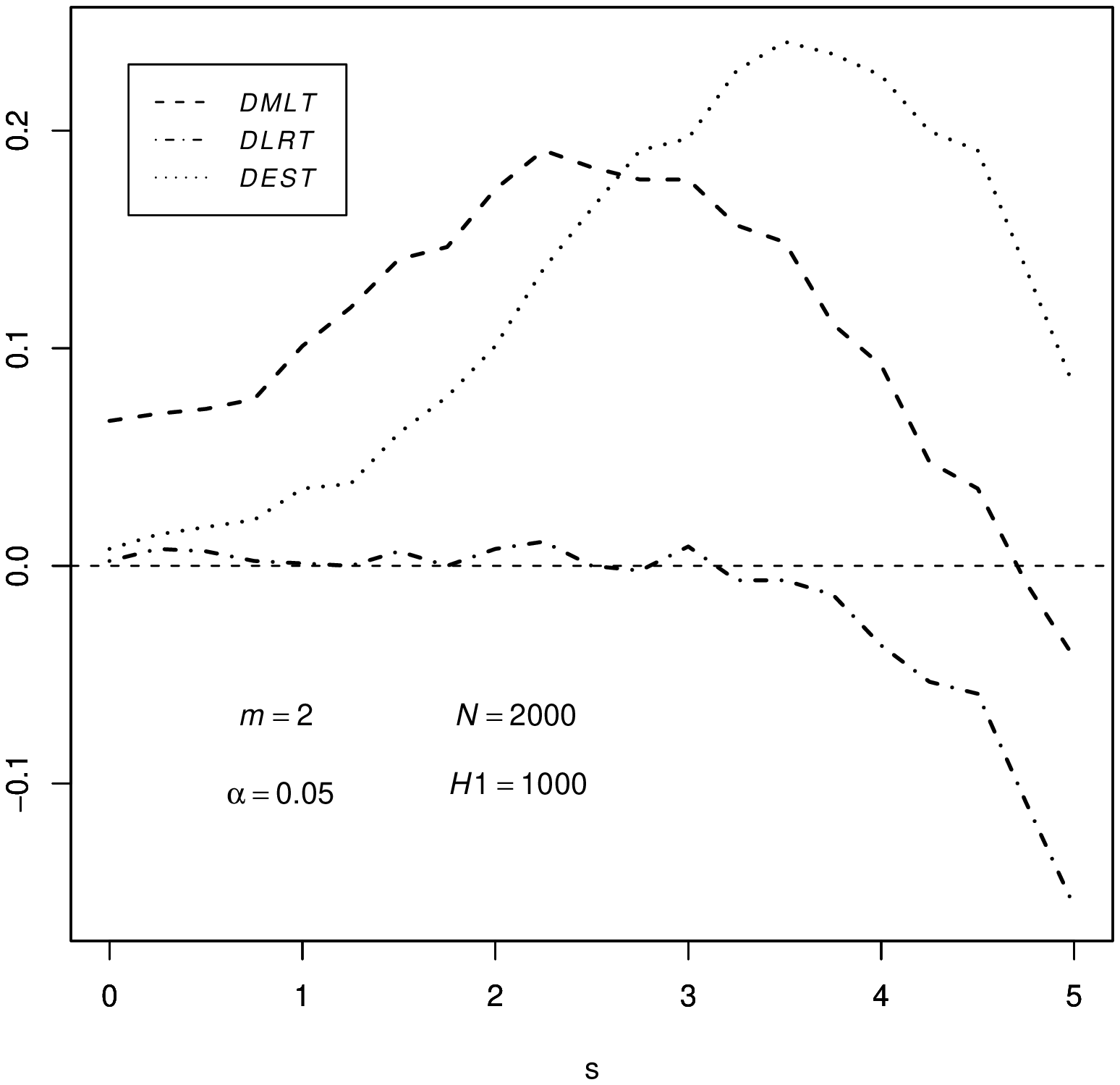}
\vspace*{-11mm}   
\caption{$\Delta P_{{\cal K}e}(s)$ \ von Mises-Fisher}
\label{fig:1-1}
\end{center}
\end{minipage}
\begin{minipage}{.5\linewidth}
\begin{center}
\includegraphics[width=8cm,height=7cm]{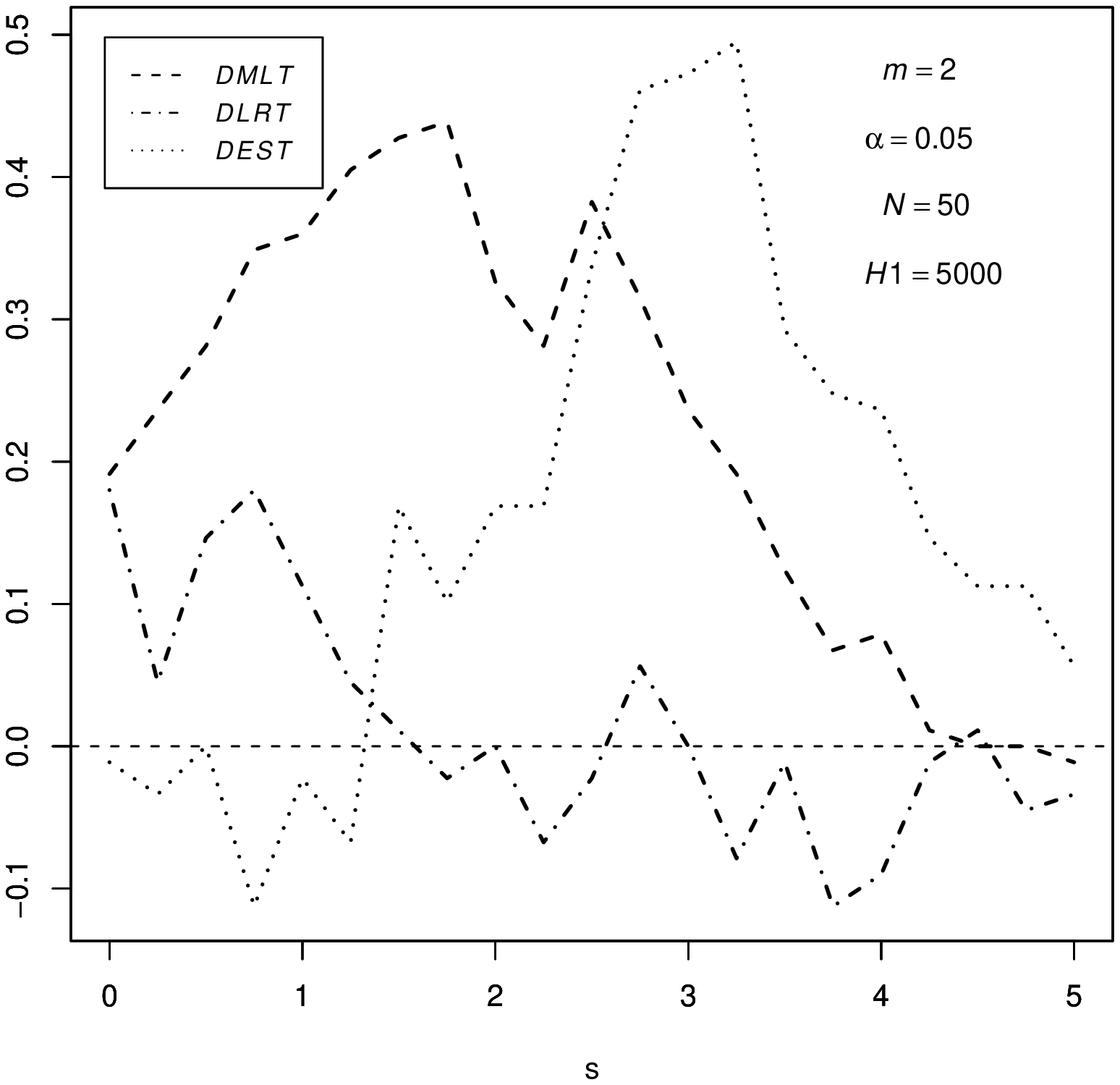}
\vspace*{-11mm}   
\caption{$\Delta P_{{\cal K}e}(s)$\ \ hyperboloid}
\label{fig:1-1}
\end{center}
\end{minipage}
\end{figure}

Next the power function of the sequential $m$-flat $k$-test is evaluated in the two dimensional case $m=2$. 
Again we take 20 values of the distance 
$s$ from 0 to 5, and for each value of $s$, we select the number of 
$H_1$ equally spaced alternatives. 
The empirical average power $\tilde{P}_{{\cal K}^*e}(s)$ of the sequential $m$-flat $k^* = (0, 0)$-test is calculated as the number of rejections relative to the number of $H_1$. 
This empirical average power $\tilde{P}_{{\cal K}^*e}(s)$ is compared with the nonsequential empirical average powers $P_{{\cal K}e}(s)$ of 
$MLT,\ LRT,\ EST$. 
The stopping time $\tau$ for the sequential test is determined by 
(\ref{eq:6-6}) 
\begin{align*}
&\tau = \inf \left\{n\ \Big|\ -\frac{1}{m}\sum_{t=1}^n \partial_a\partial_b 
l_1(x(t), {\hat u}_{mle,n})g^{ab}({\hat u}_{mle,n}) \ge 
K\nu({\hat u}_{mle,n}) + c({\hat u}_{mle,n}) \right\},\\ 
&c(u) = -\frac{1}{2}\Big(\frac{m}{rr^\dagger} - \frac{1}{r^{\dagger 2}} \Big)\ :\ \textrm{von Mises-Fisher}\quad 
c(u) = -\frac{1}{2}\Big(-\frac{m}{rr^\dagger} - \frac{1}{r^{\dagger 2}} \Big)\ :\ \textrm{hyperboloid}.
\end{align*}
Then the rejection region of the sequential $m$-flat $k^* = (0, 0)$-test is given by (\ref{eq:6-3})
\begin{align*}
&\tilde{R}_{\tilde{M}_c} = \left\{\tilde{u}_\nu\ \big|\ 
g_{\tilde{a}\tilde{b}}(u_0)
\tilde{u}_{\nu0}^{\tilde{a}}\tilde{u}_{\nu0}^{\tilde{b}} \ge 
(c_0 + \tilde{\epsilon})^2\right\},\quad
\tilde{u}_{\nu0}^{\tilde{a}} = 
\sqrt{K\nu_0}\left[\tilde{u}^{\tilde{a}}_{\nu,mle,\tau} - 
\tilde{u}_{\nu}^{\tilde{a}}(u_0)\right].
\end{align*}

As for the von Mises-Fisher model, numerical results are based on the following set of values (cf. (\ref{eq:7-1}))
\begin{align*}
&m = 2,\quad \alpha = 0.05,\quad r = 0.2,\quad
H_0\ :\ (u_0^1, u_0^2) = (\pi/2, \pi/2),\quad
N = 1000,\quad H_1 = 500,\\ 
&\tilde{u}_\nu^1 = \nu(u)\eta_1(u),\quad 
\tilde{u}_\nu^2 = \nu(u)\eta_2(u),\quad 
\nu(u) = 1/|\sin u^1|| \sin u^2|,\qquad 
K = 1000,
\end{align*}
and for the hyperboloid model, numerical results are based on the following set of values (cf. (\ref{eq:7-2}))
\begin{align*}
&m = 2,\quad \alpha = 0.05,\quad r = 2,\quad
H_0\ :\ (u_0^1, u_0^2) = (1, \pi/2),\quad
N = 50,\quad H_1 = 500,\\
&\tilde{u}_\nu^1 = \nu(u)\eta_1(u),\quad 
\tilde{u}_\nu^2 = \nu(u)\eta_2(u),\quad 
\nu(u) = 1/|\sinh u^1|| \sin u^2|,\quad
K = 60
\end{align*}

Figures 12, 13 show the von Mises-Fisher model, and Figures 14, 15 show the hyperboloid model. 
Figures 12, 14 depict 
$\tilde{P}_{{\cal K}^*e}(s),\ P_{{\cal K}e}(s)$, and the notations in the figures indicate the following tests.
\begin{align*}
&CMLT (sequential)\ :\ k^* = (0, 0),\quad
OMLT(nonsequential)\ :\ k = (0, 0),\\
&OLRT(nonsequential)\ :\ k = (0.5, 0.5),\quad
OEST(nonsequential)\ :\ k = (1, 1).
\end{align*}

Figures 13, 15 depict the differences 
$DP_{{\cal K}e}(s) = \tilde{P}_{{\cal K}^*e}(s) - P_{{\cal K}e}(s)$, 
and the notations in the figures indicate the following tests.
\begin{align*}
DMLT\ :\ k = (0, 0),\quad
DLRT\ :\ k = (0.5, 0.5),\quad
DEST\ :\ k = (1, 1).
\end{align*}

We see that in each model, the empirical average power $\tilde{P}_{{\cal K}^*e}(s)$ of the sequential $m$-flat $k^* = (0, 0)$-test $(CMLT)$ surpasses the empirical average powers $P_{{\cal K}e}(s)$ of typical nonsequential $m$-flat $k$-tests $(OMLT,\ OLRT,\ OEST)$ almost uniformly in the distance $s > 0$.

\clearpage

\begin{figure}[htbp]
\begin{minipage}{.5\linewidth}
\begin{center}
\includegraphics[width=8cm,height=7cm]{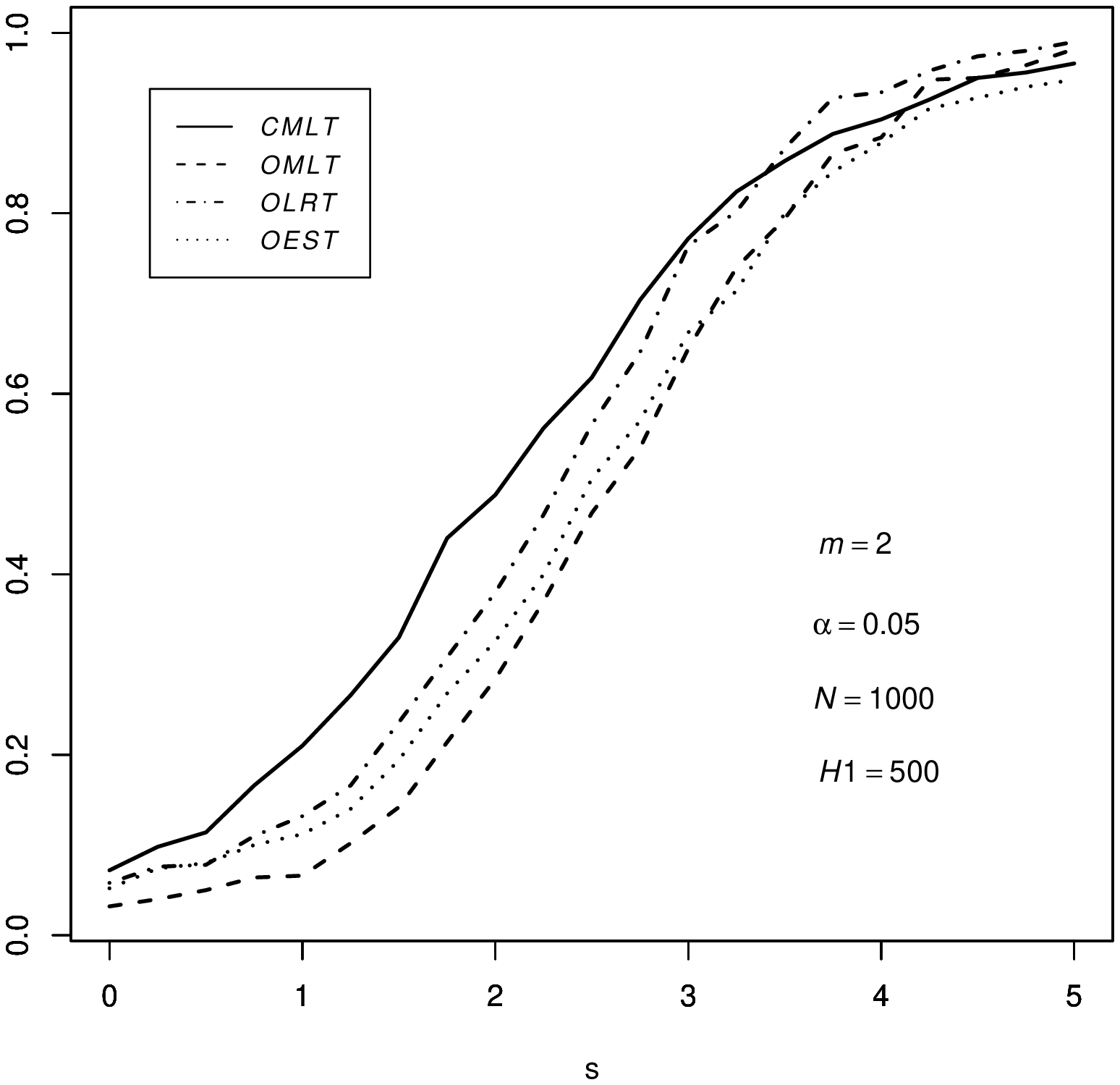}
\vspace*{-11mm}   
\caption{$\tilde{P}_{{\cal K}^*e}(s), P_{{\cal K}e}(s)$\ \ 
von Mises-Fisher}
\label{fig:1-2}
\end{center}
\end{minipage}
\begin{minipage}{.5\linewidth}
\begin{center}
\includegraphics[width=8cm,height=7cm]{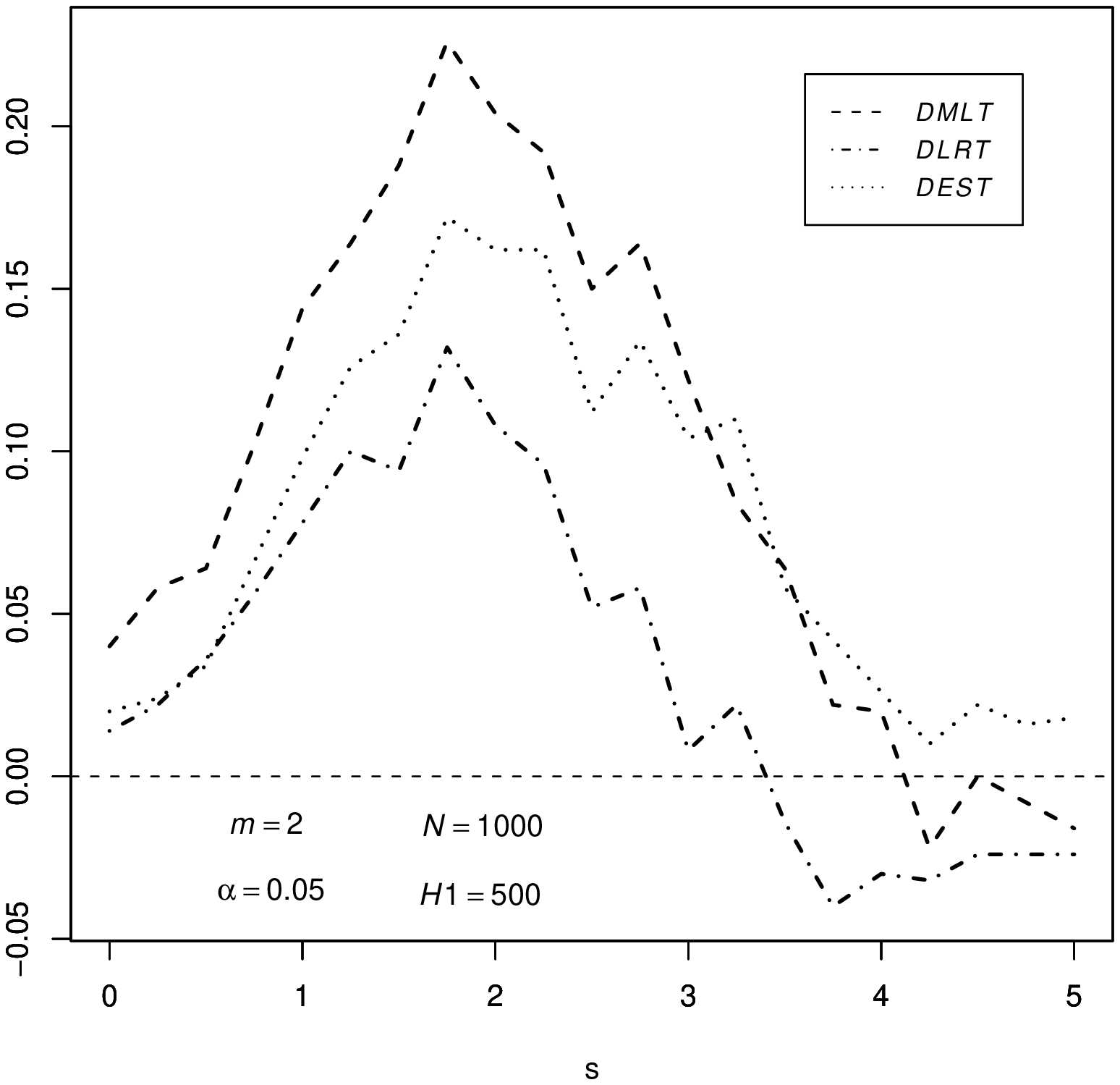}
\vspace*{-11mm}   
\caption{$DP_{{\cal K}e}(s)$\ \ von Mises-Fisher}
\label{fig:1-3}
\end{center}
\end{minipage}
\end{figure}

\begin{figure}[htbp]
\begin{minipage}{.5\linewidth}
\begin{center}
\includegraphics[width=8cm,height=7cm]{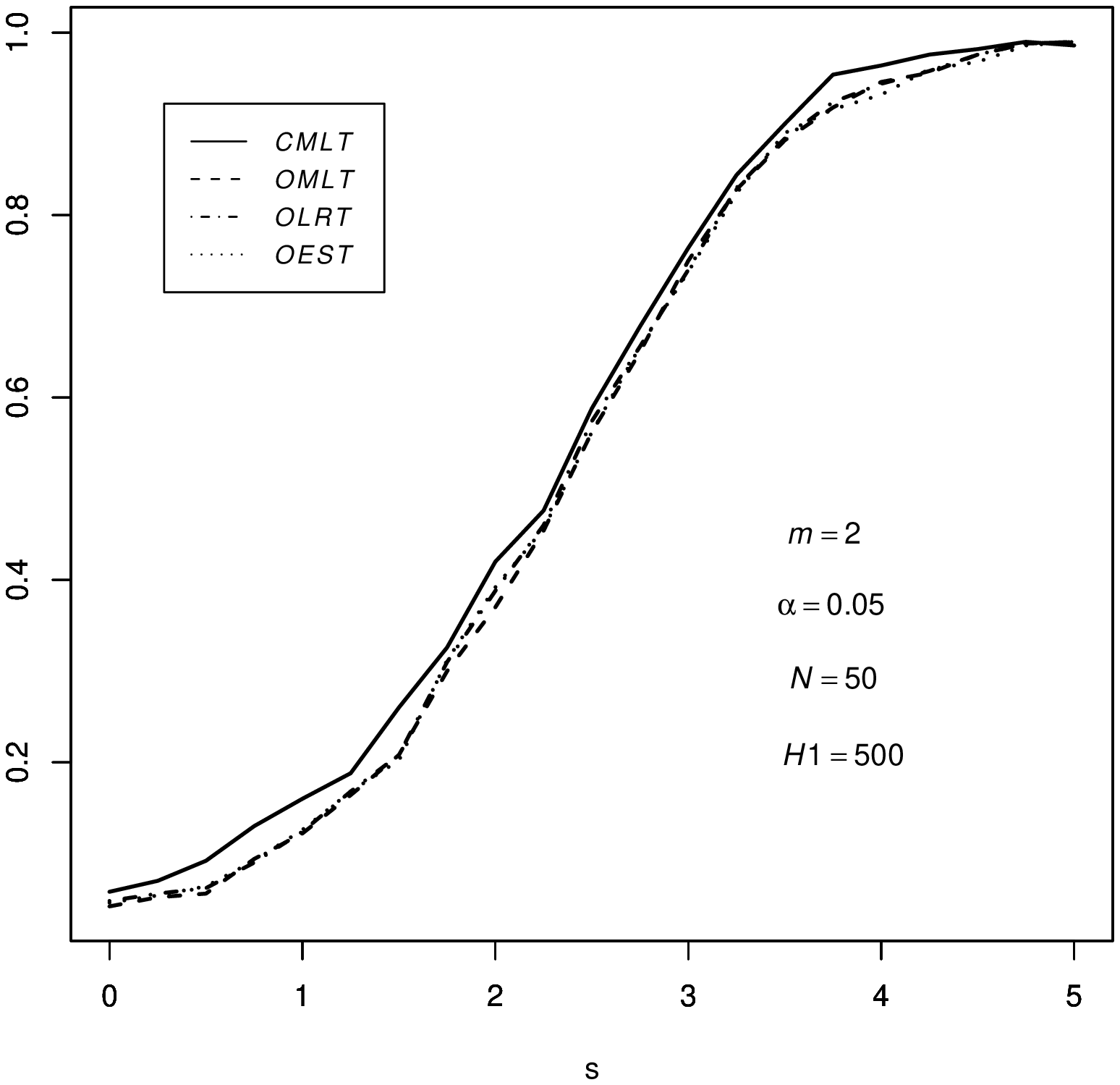}
\vspace*{-11mm}   
\caption{$\tilde{P}_{{\cal K}^*e}(s), P_{{\cal K}e}(s)$\ \ hyperboloid}
\label{fig:1-4}
\end{center}
\end{minipage}
\begin{minipage}{.5\linewidth}
\begin{center}
\includegraphics[width=8cm,height=7cm]{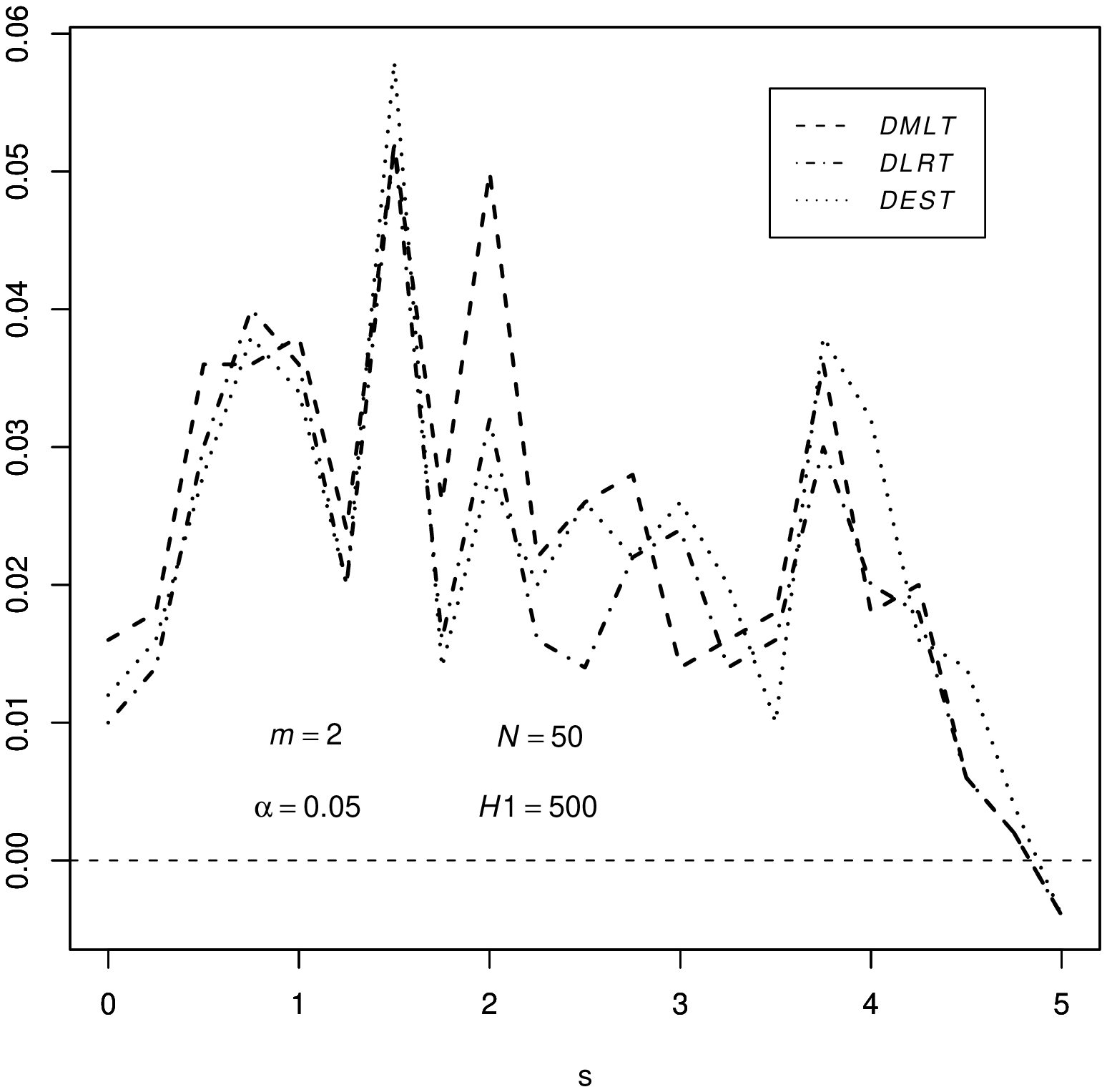}
\vspace*{-11mm}   
\caption{$DP_{{\cal K}e}(s)$\ \ hyperboloid}
\label{fig:1-5}
\end{center}
\end{minipage}
\end{figure}

\section{DISCUSSION}

Conformal geometry of statistical manifold has produced the concrete third-order efficient sequential test procedure in a multidimensional curved exponential family $M_c$. 
The method is divided into two stages: the first is to choose the stopping rule which reduces the mean 1-ES curvature of the conformally transformed $\tilde{M}_c$ to zero, and the second is to determine the significance of the null hypothesis based on the mixture-flat $k$-test statistic at the stopping time.

This method was also designed to verify that the totally exponential umbilic statistical manifold (TEU for short) guarantees the third-order uniformly efficient sequential test procedure. 
TEU implies the conformal mixture (exponential) flatness, and this flatness leads to a gauge function $\nu(u) > 0$ and a dual affine coordinate system of $\tilde{M}_c$ which reduces the $(-1)$-connection of $\tilde{M}_c$ to zero. 
This coordinate system is not pertinent to the power function itself but is effective for the simplification of the sequential mixture-flat $k$-test statistic  without the bias-correction. 
Therefore TEU can be regarded as an f.r.m.\ exponential family in the appropriate sequential inferential procedure.  
The von-Mises Fisher model and the hyperboloid model are the examples of the dual quadric hypersurface which was coined in the previous work 
(Kumon et al., 2011).
This type of hypersurface belongs to the class of TEU, and the complete characterization of TEU remains a future subject.

The present and the previous works exhibit one feature of active information geometry which encompasses geometries not only as interpretations but as strategic actions. 
It has been recognized that the sequential inferential procedures 
are rather restrictive in the sense that they cannot reduce the anisotoropy of statistical manifold. 
Another feature of the proposed spirit can be found in the field of dynamical systems, where dynamical actions such as feedback and feedforward may produce more flexible geometrical actions on system manifolds. 
This is a promising theme which is to be studied in a future work.

\clearpage

\noindent
{\textbf{\Large{APPENDICES}}\\
\

\

\noindent
\textbf{\Large{}1\ \ Edgeworth expansion of $f(\hat{w}, u)$}\\
\

\noindent
The observed sufficient statistic $\bar{x}$ is decomposed into 
$(\hat{u}, \hat{v})$ by $\bar{x} = \eta(\hat{u}, \hat{v})$, implying that 
$\bar{x}$ is in the ancillary subspace $A(\hat{u})$ and its coordinates are $\hat{v}$ in $A(\hat{u})$.
Then we define the normalized and bias-corrected statistics 
$\tilde{u}_{s,e}^*$ and $\tilde{v}$ by
\begin{align*}
\tilde{u}_{s,e} = \sqrt{N}(\hat{u} - u_{s,e}),\quad 
\tilde{u}_{s,e}^{*a} = \tilde{u}_{s,e}^a + 
\frac{1}{2\sqrt{N}}\Gamma_{\beta \gamma}^{(-1)a}g^{\beta \gamma}(\hat{u}),\quad \tilde{v} = \sqrt{N}\hat{v},
\end{align*}
and denote them by $\tilde{w}_{s,e}^* = (\tilde{u}_{s,e}^*, \tilde{v})$. 
The density functions of the statistics $\tilde{w}_{s,e}^*$ and 
$\tilde{u}_{s,e}^*$ are given in Amari and Kumon (1983), and we only enumerate the results.

When the true parameter is $u_{s,e}$, the density function 
$f(\tilde{w}_{s,e}^* ; u_{s,e})$ is expanded as
\begin{align}
\label{eq:8-1}
f(\tilde{w}_{s,e}^* ; u_{s,e}) = 
\frac{|g_{\alpha \beta}(u_{s,e})|^{1/2}}{(2\pi)^{n/2}}
\exp\left\{-\frac{1}{2}g_{\alpha \beta}(u_{s,e})
\tilde{w}_{s,e}^{*\alpha}\tilde{w}_{s,e}^{*\beta}\right\} + O(N^{-1/2}).
\end{align}
The distribution of $\tilde{u}_{s,e}^*$ is obtained by integrating 
(\ref{eq:8-1}) with respect to $\tilde{v}$,
\begin{align}
\label{eq:8-2}
f(\tilde{u}_{s,e}^* ; u_{s,e}) = 
\frac{|\bar{g}_{ab}(u_{s,e})|^{1/2}}{(2\pi)^{m/2}} 
\exp\left\{-\frac{1}{2}\bar{g}_{ab}(u_{s,e})
\tilde{u}_{s,e}^{*a}\tilde{u}_{s,e}^{*b}\right\} + O(N^{-1/2}),
\end{align}
where 
\begin{align}
\label{eq:8-3}
\bar{g}_{ab}(u_{s,e}) = g_{ab}(u_{s,e}) - 
g_{a\kappa}g_{b\lambda}g^{\kappa \lambda}.
\end{align}
The matrix $g_{ab} - \bar{g}_{ab}$ is positive semi-definite, and 
$\bar{g}_{ab}(u_{s,e}) = g_{ab}(u_{s,e})$ holds if and only if 
$g_{a\kappa}(u_{s,e}) = 0$, that is, $A$ is orthogonal at $u_{s,e}$.

When $A$ is locally orthogonal at $u_0$, the density function of 
$\tilde{u}_{s,e}^*$ is further expanded as
\begin{align}
\label{eq:8-4}
f(\tilde{u}_{s,e}^* ; u_{s,e}) = 
\varphi(\tilde{u}_{s,e}^* ; u_{s,e})\left[1 + A_N(\tilde{u}_{s,e}^* ; u_{s,e}) + B_N(\tilde{u}_{s,e}^* ; u_{s,e})\right] + O(N^{-3/2}),
\end{align}
where
\begin{align}
\label{eq:8-5}
&\varphi(\tilde{u}_{s,e}^* ; u_{s,e}) = 
\frac{|g_{ab}(u_{s,e})|^{1/2}}{(2\pi)^{m/2}} 
\exp \left\{-\frac{1}{2}g_{ab}(u_{s,e})
\tilde{u}_{s,e}^{*a}\tilde{u}_{s,e}^{*b}\right\},\\
\label{eq:8-6}
&A_N(\tilde{u}_{s,e}^* ; u_{s,e}) = 
\frac{1}{6\sqrt{N}}K_{abc}h^{abc} \nonumber \\
&\qquad \qquad \qquad \ \ + \frac{1}{4N}K_{ab}h^{ab} + 
\frac{1}{24N}K_{abcd}h^{abcd} + 
\frac{1}{72N}K_{abc}K_{def}h^{abcdef},\\
&K_{abc} = T_{abc} - 3C_{abc},\nonumber \\
&K_{ab} = C_{cda}C_{efd}g^{ce}g^{df} + 
2H_{ac\kappa}^{(1)}H_{bd\lambda}^{(1)}g^{cd}g^{\kappa \lambda},\nonumber \\
&K_{abcd} = 
S_{abcd} - 4D_{abcd} + 12(C_{eab} + C_{abe} - T_{abe})C_{fcd}g^{ef},\nonumber \\
&T_{abc} = T_{ijk}B_a^iB_b^jB_c^k,\quad 
C_{abc} = \partial_bB_{ai}B_c^i = \Gamma_{abc}^{(-1)},\nonumber \\
&S_{abcd} = \partial_i\partial_j\partial_k\partial_l\psi B_a^iB_b^jB_c^kB_d^l,\quad
D_{abcd} = \partial_c\partial_bB_{ai}B_d^i,\nonumber 
\end{align}
\begin{align}
\label{eq:8-7}
&2NB_N(\tilde{u}_{s,e}^* ; u_{s,e}) = 
\left[Q_{ac\kappa}\left\{\left(g^{ac} + s^2e^ae^c\right)Q_{bd\lambda} - 
2g^{ac}H_{bd\lambda}^{(1)}\right\}g^{\kappa \lambda} + 
\frac{1}{2}H_{ab}^{(-1)2}\right]h^{ab} \nonumber \\
&\qquad \qquad \qquad \qquad \ +sQ_{cd\lambda}
\left(2Q_{ab\kappa} - H_{ab\kappa}^{(1)}\right)g^{\kappa \lambda}
e^dh^{abc} + 
Q_{ab\kappa}\left(Q_{cd\lambda} - H_{cd\lambda}^{(1)}\right)
g^{\kappa \lambda}h^{abcd},\\
&H_{ab}^{(-1)2} = H_{\kappa \lambda a}^{(-1)}H_{\mu \nu b}^{(-1)}
g^{\kappa \mu}g^{\lambda \nu},\nonumber
\end{align}
and $h^{ab}$ etc. are the tensorial Hermite polynomials in 
$\tilde{u}_{s,e}$ as follows.
\begin{align*}
&h^{ab}(u) = u^au^b - g^{ab},\quad
h^{abc}(u) = u^au^bu^c - 3g^{(ab}u^{c)},\\
&h^{abcd}(u) = u^au^bu^cu^d - 6g^{(ab}u^cu^{d)} + 3g^{(ab}g^{cd)},\\
&h^{abcdef}(u) = u^au^bu^cu^du^eu^f - 15g^{(ab}u^cu^du^eu^{f)} + 45g^{(ab}g^{cd}u^eu^{f)} - 15g^{(ab}g^{cd}g^{ef)}.
\end{align*}

The term $A_N(\tilde{u}_{s,e}^* ; u_{s,e})$ does not depend on the ancillary family $A$, so that it is common to all the efficient tests. 
The term $B_N(\tilde{u}_{s,e}^* ; u_{s,e})$ shows that the density function 
$f(\tilde{u}_{s,e}^* ; u_{s,e})$ depends on $A$ through the two geometrical quantities $Q_{ab\kappa}$ and $H_{\kappa \lambda a}^{(-1)}$. 
Thus when $A$ is associated with a test, the third-order characteristics 
of an efficient test are determined by the angles $Q_{ab\kappa}$ between the boundary $\partial R$ of the critical region $R$ and $M_c$, and the $(-1)$-ES curvature $H_{\kappa \lambda a}^{(-1)}$ of $\partial R$.

\clearpage

\noindent
\textbf{\Large{}2\ \ Outline on the second- and third-order powers}\\
\

\noindent
From the third-order level condition 
$P_{{\cal T}}(0) = \alpha + O(N^{-3/2})$, we have
\begin{align*}
&\int_{\bar{R}_{M_c}^*} f(\tilde{u}_0; u_0)\ d\tilde{u}_0 + 
\epsilon \langle f(\tilde{u}_0; u_0)\rangle_{\partial \bar{R}_{M_c}^*} = 
1 - \alpha,\quad 
\bar{R}_{M_c}^* = \{\tilde{u_0}\ |\ 
g_{ab}(u_0)\tilde{u}_0^a \tilde{u}_0^b \le c_0^2\},\\
&\langle f(\tilde{u}_0; u_0)\rangle_{\partial \bar{R}_{M_c}^*} = 
\int_{\partial \bar{R}_{M_c}^*} f(\tilde{u}_0; u_0)\ d\tilde{u}_0/
S_N(\partial \bar{R}_{M_c}^*),\quad
\partial \bar{R}_{M_c}^* = \{\tilde{u_0}\ |\ 
g_{ab}(u_0)\tilde{u}_0^a \tilde{u}_0^b = c_0^2\},
\end{align*}
where $S_N(\partial \bar{R}_{M_c}^*)$ denotes the area of the $(m-1)$-dimensional sphere $\partial \bar{R}_{M_c}^*$ with the radius $c_0$.

We denote by $\epsilon_0$ for the m.l.e.\ test, and for a general efficient test, we put $\epsilon = \epsilon_0 + \epsilon_1$. 
Then from the expression of $A_N(\tilde{u}_0^* ; u_0)$ given by 
(\ref{eq:8-6}) we have for $\epsilon_0$ 
\begin{align*}
&N\epsilon_0 = \frac{1}{4m}K_{ab}g^{ab}h_1(c_0) + 
\frac{1}{24m(m+2)}K_{abcd}3g^{(ab}g^{cd)}h_3(c_0)\\ 
&\qquad \quad
+ \frac{1}{72m(m+2)(m+4)}K_{abc}K_{def}15g^{(ab}g^{cd}g^{ef)}h_5(c_0),\\
&h_1(c_0) = c_0,\quad 
h_3(c_0) = c_0^3 - (m+2)c_0,\quad
h_5(c_0) = c_0^5 - 2(m+4)c_0^3 + (m+2)(m+4)c_0. 
\end{align*}
On the other hand, from the expression of $B_N(\tilde{u}_0^* ; u_0)$ 
given by (\ref{eq:8-7}) we have for $\epsilon_1$ 
\begin{align*}
&N\epsilon_1 = \left[\frac{1}{4m}H_A^{(-1)2} + 
\frac{1}{2m}Q_{ab\kappa}\left(Q_{cd\lambda} - 2H_{cd\lambda}^{(1)}\right)
g^{ac}g^{bd}g^{\kappa \lambda}\right]h_1(c_0)\\
&\qquad \quad + \left[\frac{1}{2m(m+2)}
Q_{ab\kappa}\left(Q_{cd\lambda} - H_{cd\lambda}^{(1)}\right)
g^{(ab}g^{cd)}g^{\kappa \lambda}\right]h_3(c_0).
\end{align*}

Since it needs long and complicated calculations to derive Proposition 3.2, we only enumerate key formulae on the integrations of the tensorial Hermite polynomials which were also used to get $\epsilon_0, \epsilon_1$. 
\begin{align*}
&\int_{\bar{R}_{M_c}^*} \langle h^{ab}\varphi(\tilde{u}_{s,e}^*; u_{s,e})\ 
d\tilde{u}_{s,e}^* \rangle_e = -\frac{1}{m}g^{ab}h_1(c_0),\\
&\int_{\bar{R}_{M_c}^*} \langle h^{abcd}\varphi(\tilde{u}_{s,e}^*; u_{s,e})\ d\tilde{u}_{s,e}^* \rangle_e = 
-\frac{1}{m(m+2)}3g^{(ab}g^{cd)}h_3(c_0),\\
&\int_{\bar{R}_{M_c}^*} \langle h^{abcdef}\varphi(\tilde{u}_{s,e}^*; u_{s,e})\ d\tilde{u}_{s,e}^* \rangle_e = 
-\frac{1}{m(m+2)(m+4)}15g^{(ab}g^{cd}g^{ef)}h_5(c_0).
\end{align*}

\clearpage

\noindent
\textbf{\Large{}3\ \ Edgeworth expansion of $\tilde{f}(\hat{w}, u)$}\\
\

\noindent
We define the normalized and bias-corrected statistics 
$\tilde{u}_{n,s,e}^*$ and $\tilde{v}_n$ by
\begin{align*}
\tilde{u}_{n,s,e} = \sqrt{K\nu_0}(\hat{u}_n - u_{s,e}),\quad 
\tilde{u}_{n,s,e}^{*a} = \tilde{u}_{n,s,e}^a + 
\frac{1}{2\sqrt{K\nu_0}}\tilde{\Gamma}_{\beta \gamma}^{(-1)a}
g^{\beta \gamma}(\hat{u}_n),\quad 
\tilde{v}_n = \sqrt{K\nu_0}\hat{v},
\end{align*}
and denote them by 
$\tilde{w}_{n,s,e}^* = (\tilde{u}_{n,s,e}^*, \tilde{v}_n)$.

When the true parameter is $u_{s,e}$, the density function 
$\tilde{f}(\tilde{w}_{n,s,e}^* ; u_{s,e})$ of $\tilde{w}_{n,s,e}^*$ is expanded as
\begin{align}
\label{eq:8-8}
\tilde{f}(\tilde{w}_{n,s,e}^* ; u_{s,e}) = 
\frac{|g_{\alpha \beta}(u_{s,e})|^{1/2}}{(2\pi)^{n/2}}
\exp\left\{-\frac{1}{2}g_{\alpha \beta}(u_{s,e})
\tilde{w}_{n,s,e}^{*\alpha}\tilde{w}_{n,s,e}^{*\beta}\right\} + O(K^{-1/2}).
\end{align}
The distribution of $\tilde{u}_{n,s,e}^*$ is obtained by integrating 
(\ref{eq:8-8}) with respect to $\tilde{v}_n$,
\begin{align}
\label{eq:8-9}
\tilde{f}(\tilde{u}_{n,s,e}^* ; u_{s,e}) = 
\frac{|\bar{g}_{ab}(u_{s,e})|^{1/2}}{(2\pi)^{m/2}} 
\exp\left\{-\frac{1}{2}\bar{g}_{ab}(u_{s,e})
\tilde{u}_{n,s,e}^{*a}\tilde{u}_{n,s,e}^{*b}\right\} + O(K^{-1/2}),
\end{align}
where $\bar{g}_{ab}(u_{s,e})$ is the same as (\ref{eq:8-3}).

When $A$ is locally orthogonal at $u_0$, the density function of 
$\tilde{u}_{n,s,e}^*$ is further expanded as
\begin{align}
\label{eq:8-10}
\tilde{f}(\tilde{u}_{n,s,e}^* ; u_{s,e}) = 
\varphi(\tilde{u}_{n,s,e}^* ; u_{s,e})\left[1 + 
\tilde{A}_K(\tilde{u}_{n,s,e}^* ; u_{s,e}) + 
\tilde{B}_K(\tilde{u}_{n,s,e}^* ; u_{s,e})\right] + O(K^{-3/2}),
\end{align}
where
\begin{align}
\label{eq:8-11}
&\varphi(\tilde{u}_{n,s,e}^* ; u_{s,e}) = 
\frac{|g_{ab}(u_{s,e})|^{1/2}}{(2\pi)^{m/2}} 
\exp \left\{-\frac{1}{2}g_{ab}(u_{s,e})
\tilde{u}_{n,s,e}^{*a}\tilde{u}_{n,s,e}^{*b}\right\},\\
\label{eq:8-12}
&\tilde{A}_K(\tilde{u}_{n,s,e}^* ; u_{s,e}) = 
\frac{1}{6\sqrt{K\nu_0}}\tilde{K}_{abc}h^{abc} \nonumber \\
&\qquad \qquad \qquad \ \ \ + \frac{1}{4K\nu_0}\tilde{K}_{ab}h^{ab} + 
\frac{1}{24K\nu_0}\tilde{K}_{abcd}h^{abcd} + 
\frac{1}{72K\nu_0}\tilde{K}_{abc}\tilde{K}_{def}h^{abcdef},\\
&\tilde{K}_{abc} = \tilde{T}_{abc} - 3\tilde{C}_{abc},\nonumber \\
&\tilde{K}_{ab} = \tilde{C}_{cda}\tilde{C}_{efd}g^{ce}g^{df} + 
2\tilde{H}_{ac\kappa}^{(1)}\tilde{H}_{bd\lambda}^{(1)}g^{cd}
g^{\kappa \lambda},\nonumber \\
&\tilde{K}_{abcd} = 
\tilde{S}_{abcd} - 4\tilde{D}_{abcd} + 
12(\tilde{C}_{eab} + \tilde{C}_{abe} - \tilde{T}_{abe})
\tilde{C}_{fcd}g^{ef},\nonumber \\
&\tilde{C}_{abc} = \tilde{\Gamma}_{abc}^{(-1)},\quad
\tilde{D}_{abcd} = \partial_c\tilde{C}_{abd} - \tilde{C}_{abe}
(\tilde{C}_{cdf} - \tilde{T}_{cdf})g^{ef},\nonumber \\ 
&\tilde{S}_{abcd} = S_{abcd} + 6s_{(ab}g_{cd)} + 4s_{(a}T_{bcd)} + 
(1 - s_\alpha s_\beta g^{\alpha \beta})3g_{(ab}g_{cd)},\quad
s_{ab} = \partial_a s_b + s_as_b - \tilde{\Gamma}_{ab}^{(1)c}s_c, \nonumber
\end{align}
\begin{align}
\label{eq:8-13}
&2K\nu_0\tilde{B}_K(\tilde{u}_{n,s,e}^* ; u_{s,e}) = 
\left[\tilde{Q}_{ac\kappa}\left\{\left(g^{ac} + s^2e^ae^c\right)
\tilde{Q}_{bd\lambda} - 
2g^{ac}\tilde{H}_{bd\lambda}^{(1)}\right\}g^{\kappa \lambda} + 
\frac{1}{2}H_{ab}^{(-1)2}\right]h^{ab} \nonumber \\
&\qquad \qquad \qquad \qquad \ +s\tilde{Q}_{cd\lambda}
\left(2\tilde{Q}_{ab\kappa} - \tilde{H}_{ab\kappa}^{(1)}\right)
g^{\kappa \lambda}
e^dh^{abc} + 
\tilde{Q}_{ab\kappa}
\left(\tilde{Q}_{cd\lambda} - \tilde{H}_{cd\lambda}^{(1)}\right)
g^{\kappa \lambda}h^{abcd}.
\end{align}
The term $\tilde{A}_K(\tilde{u}_{n,s,e}^* ; u_{s,e})$ does not depend on the ancillary family $\tilde{A}$, so that it is common to all the efficient sequential tests. 
The term $\tilde{B}_K(\tilde{u}_{n,s,e}^* ; u_{s,e})$ shows that the density function 
$\tilde{f}(\tilde{u}_{n,s,e}^* ; u_{s,e})$ depends on $\tilde{A}$ through the two geometrical quantities $\tilde{Q}_{ab\kappa}$ and 
$H_{\kappa \lambda a}^{(-1)}$. 
Thus when $\tilde{A}$ is associated with a sequential test, the third-order characteristics of an efficient sequential test are determined by the angles $\tilde{Q}_{ab\kappa}$ between the boundary $\partial \tilde{R}$ of the critical region $\tilde{R}$ and 
$\tilde{M}_c$, and the $(-1)$-ES curvature $H_{\kappa \lambda a}^{(-1)}$ of $\partial \tilde{R}$.


\clearpage



\begin{thebibliography}{99}


\bibitem{Akahira and Takeuchi}
Akahira, M. and Takeuchi, K. (1989). 
Third Order Asymptotic Efficiency of the Sequential Maximum Likelihood Estimation Procedure,  
{\it Sequential Analysis} 8: 333-359.


\bibitem{Amari a}
Amari, S. (1985). 
{\it Differential Geometrical Methods in Statistics,   
Lecture Notes in Statistics} 28, 2nd print, New York: Springer.


\bibitem{Amari b}
Amari, S., Barndorff-Nielsen, O. E., Kass, R. E., Lauritzen, S. L.  
and Rao, C. R. (1987). 
{\it Differential Geometry in Statistical Inference,}
Hayward, Calif.: Institute of Mathematical Statistics.


\bibitem{Amari and Kumon a}
Amari, S. and Kumon, M. (1983). 
Differential Geometry of Edgeworth Expansions in Curved Exponential Family, 
{\it Annals of the Institute of Statistical Mathematics} 35: 1-24.


\bibitem{Amari and Kumon b}
Amari, S. and Kumon, M. (1988). 
Estimation in the Presence of Infinitely Many Nuisance 
Parameters-Geometry of Estimating Functions, 
{\it Annals of Statistics} 16: 1044-1068.


\bibitem{Amari and Nagaoka}
Amari, S. and Nagaoka, H. (2000). 
{\it Methods of Information Geometry,} 
Providence: American Mathematical Society and Oxford University Press.


\bibitem{Barndorff-Nielsen et al.}
Barndorff-Nielsen, O. E., Bl\ae sild, P. and Eriksen, P. S. (1989). 
{\it Decomposition and Invariance of Measures, and Statistical Transformation Models,   
Lecture Notes in Statistics} 58, New York: Springer.


\bibitem{Gut}
Gut, A. (2009). 
{\it Stopped Random Walks, Limit Theorems and Applications,} 
2nd edition, New York: Springer.


\bibitem{Kumon a}
Kumon, M. (2009). 
On the Conditions for the Existence of Ancillary Statistics in a Curved Exponential Family,  
{\it Statistical Methodology} 6: 320-335.


\bibitem{Kumon b}
Kumon, M. (2010). 
Studies of Information Quantities and Information Geometry of Higher Order Cumulant Spaces,  
{\it Statistical Methodology} 7: 152-172.


\bibitem{Kumon and Amari}
Kumon, M. and Amari, S. (1983). 
Geometrical Theory of Higher-Order Asymptotics of Test, Interval Estimator and Conditional Inference, 
{\it Proceedings of Royal Society of London} A 387: 429-458.


\bibitem{Kumon et al.}
Kumon, M., Takemura, A and Takeuchi, K. (2011). 
Conformal Geometry of Statistical Manifold with Application to Sequential Estimation, {\it Sequential Analysis,} 30: 308-337.


\bibitem{Mariera}
Magiera, R. (1974). 
On the Inequality of Cram\'er-Rao Type in Sequential Estimation Theory, 
{\it Zastosowania Matematyki (Applicationes Mathematicae)} 14: 227-235.


\bibitem{Okamoto, Amari and Takeuchi}
Okamoto, I., Amari, S. and Takeuchi, K. (1991). 
Asymptotic Theory of Sequential Estimation: Differential Geometrical Approach, 
{\it Annals of Statistics} 19: 961-981.


\bibitem{Takeuchi and Akahira}
Takeuchi, K. and Akahira, M.(1988). 
Second Order Asymptotic Efficiency in Terms of Asymptotic Variances of the Sequential Maximum Likelihood Estimation Procedures, In 
{\it Statistical Theory and Data Analysis II: Proceedings of the Second Pacific Area Statistical Conference,} 
Matusita, K. ed., 191-196, Amsterdam: North-Holland.


\bibitem{Vos}
Vos, P. W. (1989). 
Fundamental Equations for Statistical Submanifolds with Applications to the Bartlett Correction,  
{\it Annals of the Institute of Statistical Mathematics} 41: 429-450.


\bibitem{Winkler and Franz}
Winkler, W. and Franz, J. (1979). 
Sequential Estimation Problems for the Exponential Class of Processes with Independent Increments,  
{\it Scandinavian Journal of Statistics} 6: 129-139.


\end{thebibliography}
\end{document}